\documentclass[letterpaper]{amsart}
\usepackage{amsfonts,amssymb}
\usepackage{xspace}
\usepackage[dvips]{epsfig}
\usepackage[matrix,arrow,curve]{xy}
\CompileMatrices
\usepackage{eucal}
\usepackage{afterpage}
\usepackage{amscd}

\theoremstyle{plain}

\newtheorem{thm}{Theorem}[section]
\newtheorem{theorem}[thm]{Theorem}

\newtheorem{corollary}[thm]{Corollary}
\newtheorem{lemma}[thm]{Lemma}

\newtheorem{prop}[thm]{Proposition}
\newtheorem{proposition}[thm]{Proposition}

\newtheorem*{conjecture*}{Conjecture}

\newtheorem*{question*}{Question}

\theoremstyle{definition}

\newtheorem{definition}[thm]{Definition}

\newtheorem*{definitions*}{Definitions}

\newtheorem*{rem*}{Remark}
\newtheorem{remark}[thm]{Remark}
\newtheorem*{remark*}{Remark}

\newtheorem*{remarks*}{Remarks}
\newtheorem*{example*}{Example}
\newtheorem{example}[thm]{Example}
\newtheorem*{examples*}{Examples}

\newtheorem*{notation*}{Notation}
\newtheorem*{convention*}{Convention}
\newtheorem*{conventions*}{Conventions}
\newtheorem*{note*}{Note}

\newcommand{\zeroindent}{\parindent0cm \parskip1ex}

\newenvironment{theoremlist}%
{\begin{list}{{\rm(\alph{enumi}) }}{\usecounter{enumi}

\leftmargin0cm \labelsep0cm \rightmargin0cm \parsep1em \listparindent0em \itemsep-1em
\topsep1em \parskip1em \setlength{\labelwidth}{\fill}}} {\parskip0em \end{list}}

{

\begin{enumerate}}{\end{enumerate}}

{

\begin{enumerate}}{\end{enumerate}}

{

\begin{enumerate} \parsep0cm
\itemsep0cm \parskip0cm}{\end{enumerate}}

{

\begin{enumerate}}{\end{enumerate}}

{

\begin{enumerate}}{\end{enumerate}}

{\begin{list}{(\alph{enumi}) }{\usecounter{enumi}

\leftmargin0cm \labelsep0cm
\itemsep1em \rightmargin0cm
\setlength{\labelwidth}{\fill}}} {\end{list}}
\newcommand{\includefigure}[3]{%
\begin{figure}[#3]
\begin{center}
\epsfig{file=#2}\\ \caption{\label{fig:#1}}
\end{center}
\end{figure}}

\newcommand{\R}{\mathbb{R}}
\newcommand{\Z}{\mathbb{Z}}

\newcommand{\C}{\mathbb{C}}
\newcommand{\N}{\mathbb{N}}
\newcommand{\half}{{\textstyle\frac{1}{2}}}

\newcommand{\iso}{\cong}           
\newcommand{\smooth}{C^\infty}
\newcommand{\CP}[1]{\C {\mathrm P}^{#1}}
\newcommand{\RP}[1]{\R {\mathrm P}^{#1}}
\newcommand{\leftsc}{\langle}
\newcommand{\rightsc}{\rangle}

\newcommand{\Rleq}{\R^{\scriptscriptstyle \leq 0}}
\newcommand{\suchthat}{\; | \;}

\newcommand{\id}{\mathrm{id}}

\newcommand{\im}{\mathrm{im}}
\renewcommand{\ker}{\mathrm{ker}}

\newcommand{\Hom}{\mathrm{Hom}}
\newcommand{\End}{\mathrm{End}}

\newcommand{\mo}{(M,\omega)}
\renewcommand{\o}{\omega}
\renewcommand{\O}{\Omega}

\newcommand{\Symp}{\mathrm{Sp}}
\newcommand{\Diff}{\mathrm{Diff}}

\newenvironment{Blist}%
{

\begin{enumerate}}{\end{enumerate}}

\newenvironment{Jlist}%
{

\begin{enumerate}}{\end{enumerate}}

\newenvironment{Slist}%
{

\begin{enumerate}}{\end{enumerate}}

\newenvironment{Sprimelist}%
{

\begin{enumerate}}{\end{enumerate}}

\newenvironment{Plist}%
{

\begin{enumerate}}{\end{enumerate}}

{

\begin{enumerate}}{\end{enumerate}}

\newenvironment{Flist}%
{

\begin{enumerate}}{\end{enumerate}}

\newenvironment{Alist}%
{

\begin{enumerate}}{\end{enumerate}}

\newenvironment{Olist}%
{

\begin{enumerate}}{\end{enumerate}}

\newcommand{\cat}{\mathcal Q}
\newcommand{\Amod}{A_m\text{-}\mathrm{mod}}
\newcommand{\amod}{A_m\mathrm{-mod}}
\newcommand{\Ob}{\mathrm{Ob}}

\newcommand{\funct}{\mathcal F}
\newcommand{\action}{\mathcal R}
\newcommand{\iP}{{}_iP}

\newcommand{\Diffeo}{\mathcal{G}}
\newcommand{\Conf}{\mathrm{Conf}}
\newcommand{\tD}{\widetilde{D}}
\newcommand{\isotopic}{\simeq}
\newcommand{\ta}{\tilde{a}}
\newcommand{\tb}{\tilde{b}}
\newcommand{\tc}{\tilde{c}}
\newcommand{\td}{\tilde{d}}
\newcommand{\tf}{\tilde{f}}
\newcommand{\tg}{\tilde{g}}

\newcommand{\tP}{\widetilde{P}}
\newcommand{\Pchi}{\chi}

\newcommand{\Ibigr}{I^{\mathrm{bigr}}}
\newcommand{\mubigr}{\mu^{\mathrm{bigr}}}

\renewcommand{\Symp}{\mathrm{Symp}}
\newcommand{\Lag}{\mathrm{Lag}}
\newcommand{\reg}{\mathrm{reg}}

\newcommand{\J}{\mathbf{J}}
\newcommand{\JJ}{\mathcal{J}}
\renewcommand{\P}{\mathcal{P}}
\newcommand{\moduli}{\mathcal{M}}
\newcommand{\fhs}[1]{\cite[#1]{floer-hofer-salamon94}}
\newcommand{\gen}[1]{\leftsc #1 \rightsc}
\newcommand{\W}{\mathcal{W}}
\newcommand{\Y}{\mathbf{Y}}
\newcommand{\T}{\mathcal{T}}
\newcommand{\LL}{\mathcal{L}}
\newcommand{\muabs}{\mu^{\mathrm{abs}}}
\newcommand{\mupaths}{\mu^{\mathrm{paths}}}
\newcommand{\barB}{\overline{B}{}}
\newcommand{\Rg}{\R^{\scriptscriptstyle >0}}
\newcommand{\tLL}{{\widetilde{\LL}}}
\newcommand{\tL}{\widetilde{L}}

\newcommand{\oplusop}[1]{{\mathop{\oplus}\limits_{#1}}}

\newcommand{\sumop}[1]{{\mathop{\sum}\limits_{#1}}}
\newcommand{\bigoplusop}[1]{{\mathop{\bigoplus}\limits_{#1}}}

\newcommand\F{\mathbb F}

\newcommand\oo{\otimes}

\newcommand\mc{\mathcal}
\newcommand\mf{\mathfrak}
\newcommand\cR{{\mathcal{R}}}
\renewcommand\sl{\mathfrak{sl}}
\newcommand\U{\mathcal U}

\newcommand\hsm{\hspace{0.03in}}
\newcommand\vsp{\vspace{0.15in}}

\newcommand\lra{\longrightarrow}

\newcommand\ab{\mathcal{Q}}
\newcommand\db{D^b(\amod)}
\newcommand\Ccat{{\mathcal{C}}_m}
\newcommand\Oo{\mathcal{O}}

\newcommand\simp{S}        
\newcommand\nd{\mathrm{cr}}
\newcommand{\cross}{\mathrm{cr}}

\newcommand\dga{\widehat{A}}

\zeroindent
\numberwithin{equation}{section}

\begin{document}
\title[Quivers and Floer cohomology]{Quivers, Floer cohomology, and\\ braid group actions}
\author{Mikhail Khovanov and Paul Seidel}
\date{\today}
\begin{abstract}
We consider the derived categories of modules over a certain family $A_m$ ($m
\geq 1$) of graded rings, and Floer cohomology of Lagrangian intersections in
the symplectic manifolds which are the Milnor fibres of simple singularities of
type $A_m.$ We show that each of these two rather different objects encodes the
topology of curves on an $(m+1)$-punctured disc. We prove that the braid group
$B_{m+1}$ acts faithfully on the derived category of $A_m$-modules, and that it
injects into the symplectic mapping class group of the Milnor fibers. The
philosophy behind our results is as follows. Using Floer cohomology, one should
be able to associate to the Milnor fibre a triangulated category (its
construction has not been carried out in detail yet). This triangulated
category should contain a full subcategory which is equivalent, up to a slight
difference in the grading, to the derived category of $A_m$-modules. The full
embedding would connect the two occurrences of the braid group, thus explaining
the similarity between them.
\end{abstract}

\maketitle

%
%
%
%
%
%

\section{Introduction}
\subsection{Generalities\label{subsec:generalities}}

This paper investigates the connection between symplectic geometry and those parts of
representation theory which revolve around the notion of categorification. The existence of
such a connection, in an abstract sense, follows from simple general ideas. The difficult
thing is to make it explicit. On the symplectic side, the tools needed for a systematic study
of this question are not yet fully available. Therefore we concentrate on a single example,
which is just complicated enough to indicate the depth of the relationship. The results can be
understood by themselves, but a glimpse of the big picture certainly helps to explain them,
and that is what the present section is for.

Let $\cat$ be a category. An {\em action of a group $G$ on $\cat$} is a family $(\funct_g)_{g
\in G}$ of functors from $\cat$ to itself, such that $\funct_e \iso \id_\cat$ and
$\funct_{g_1} \funct_{g_2} \iso \funct_{g_1g_2}$ for all $g_1,g_2 \in G$; here $\iso$ denotes
isomorphism of functors\footnote{Strictly speaking, this should be called a \emph{weak}
action of $G$ on $\cat;$ a full-fledged action comes with preferred isomorphisms
$\funct_{g_1} \funct_{g_2} \iso \funct_{g_1g_2}$ that satisfy obvious compatibility
relations.}.
We will not distinguish between two actions $(\funct_g)$ and
$(\widetilde{\funct}_g)$ such that $\funct_g \iso \widetilde{\funct}_g$ for all $g$. A
particularly nice situation is when $\cat$ is triangulated and the $\funct_g$ are exact
functors. Then the action induces a linear representation of $G$ on the Grothendieck group
$K(\cat)$. The inverse process, in which one lifts a given linear representation to a
group action on a triangulated category, is called categorification (of group
representations).

The connection with symplectic geometry is based on an idea of Donaldson. He
proposed (in talks circa $1994$) to associate to a compact symplectic manifold
$(M^{2n},\o)$ a category $\Lag\mo$ whose objects are Lagrangian submanifolds $L
\subset M$, and whose morphisms are the Floer cohomology groups\footnote{The
definition of Floer cohomology in general involves difficult analytic and
algebraic questions. To simplify the exposition, we tacitly ignore them here.}
$HF(L_0,L_1)$. The composition of morphisms would be given by products
$HF(L_1,L_2) \times HF(L_0,L_1) \longrightarrow HF(L_0,L_2),$ which are
defined, for example, in  \cite{desilva98}. Let $\Symp\mo$ be the group of
symplectic automorphisms of $M$. Any $\phi \in \Symp\mo$ determines a family of
isomorphisms $HF(L_0,L_1) \iso HF(\phi L_0,\phi L_1)$ for $L_0,L_1 \in
\Ob(\Lag\mo)$ which are compatible with the products. In other words, $\phi$
induces an equivalence $\funct_\phi$ from $\Lag\mo$ to itself. This is just a
consequence of the fact that $\Lag\mo$ is an object of symplectic geometry,
hence natural under symplectic maps. Assume for simplicity that $H^1(M;\R) =
0$, so that all symplectic vector fields are Hamiltonian. Then a smooth isotopy
of symplectic automorphisms $(\phi_t)_{0 \leq t \leq 1}$ gives rise to
distinguished elements in $HF(\phi_0 L,\phi_1 L)$ for all $L.$ This means that
the ${\mathcal F}_\phi$ define a canonical action of the symplectic mapping
class group $\pi_0(\Symp\mo)$ on $\Lag\mo$. As a consequence, any symplectic
fibre bundle with fibre $\mo$ and base $B$ gives rise to a $\pi_1(B)$-action on
$\Lag\mo$, through the monodromy map $\pi_1(B) \rightarrow \pi_0(\Symp\mo)$.
Interesting examples can be obtained from families of smooth complex projective
varieties.

While it is thus easy to construct potentially interesting group actions on the
categories $\Lag\mo$, both the group actions and the categories themselves are
difficult to understand. A particularly intriguing question is whether, in any
given case, one can relate them to objects defined in a purely algebraic way.
In homological algebra there is a standard technique for approaching similar
comparison problems (an example is Beilinson's work \cite{beilinson78} on
coherent sheaves on $\CP{k}$). A very crude attempt to adapt this technique to
our situation goes like this: pick a finite number of Lagrangian submanifolds
$L_0, \dots, L_m \subset M$ which, for some reason, appear to be particularly
important. Using the product on Floer cohomology, turn
\begin{equation} \label{eq:hom-algebra}
A = \bigoplus_{i,j = 0}^m HF(L_i,L_j)
\end{equation}
into a ring. Now associate to an arbitrary Lagrangian submanifold $L \subset M$
the $A$-module $\bigoplus_i HF(L,L_i)$. This defines a functor from $\Lag\mo$
to the category of $A$-modules. As it stands the functor is not particularly
useful, since it does not take into account all the available structure.

To begin with, Floer cohomology groups are graded, and hence $A$ is a graded
ring (depending on $M$, this may be only a $\Z/N$-grading for some finite $N$).
To make more substantial progress one needs to refine the category $\Lag\mo$.
This theory is as yet under construction, and we can only give a vague outline
of it. As pointed out by Fukaya \cite{fukaya93}, working directly with the
Floer cochain complexes should enable one to construct an $A_\infty$-category
underlying $\Lag\mo$. Kontsevich \cite{kontsevich94} suggested to consider the
derived category of this $A_\infty$-category. This ``derived Fukaya category''
is expected to be triangulated, and to contain $\Lag\mo$ as a full subcategory;
we denote it by $D^b\Lag\mo$ (which is an abuse of notation). The action of
$\pi_0(\Symp\mo)$ on $\Lag\mo$ should extend to an action on $D^b\Lag\mo$ by
exact functors. Moreover, it seems natural to suppose that the Grothendieck
group of $D^b\Lag\mo$ is related to $H_n(M;\Z)$; that would mean that the group
actions coming from families of projective varieties could be considered as
categorifications of the classical monodromy representations. Returning to the
rings \eqref{eq:hom-algebra}, one expects to get from Fukaya's construction a
canonical (up to quasi-isomorphism) $A_\infty$-algebra ${\mathcal A}$ with
cohomology $A$. The functor introduced above would lift to an
$A_\infty$-functor from the $A_\infty$-category underlying $\Lag\mo$ to the
$A_\infty$-category of $A_\infty$-modules over ${\mathcal A}$. In a second
step, this would induce an exact functor
\[
D^b\Lag\mo \longrightarrow D({\mathcal A}),
\]
where $D({\mathcal A})$ is the derived category of $A_\infty$-modules. A
standard argument based on exactness indicates that this functor, when
restricted to the triangulated subcategory of $D^b\Lag\mo$ generated by
$L_0,\dots, L_m$, would be full and faithful. This is a much stronger
comparison theorem than one could get with the primitive approach which we had
mentioned first. We will not say more about this; what the reader should keep in
mind is the expected importance of the rings $A$.

It should also be mentioned that there is a special reason to study actions of
braid groups and mapping class groups on categories. In 1992--1994 a number of
mathematicians independently suggested extending the formalism of TQFTs to
manifolds with corners. While three-dimensional extended TQFTs are rather well
understood, there has been much less progress in the four-dimensional case.
Part of the data which make up a four-dimensional extended TQFT are categories
${\mathcal Q}(S)$ associated to closed surfaces $S$; these would govern the
process of cutting and pasting for the vector spaces associated to closed
three-manifolds. The mapping class group of $S$ should act on ${\mathcal Q}(S)$
in a natural way, and for interesting TQFTs this action is likely to be highly
nontrivial. As a matter of historical interest, we note that the categories
$\Lag\mo$ were invented by Donaldson in an attempt to make gauge theory fit
into the extended TQFT formalism; see \cite{fukaya97} for more recent
developments.

There is a version of this for invariants of two-knots, where the categories
are associated to punctured discs and carry actions of the braid groups. The
papers \cite{bernstein-frenkel-khovanov98} and \cite{khovanov98} go some way
towards constructing such extended TQFTs in a representation-theoretic way.

\subsection{The results\label{subsec:results}}

In the algebraic part of the paper, we study the derived categories $D^b(\Amod)$ for a
particular family $A_m$ ($m \geq 1$) of graded rings. The definition of $A_m$ starts with the
quiver (oriented graph) $\Gamma_m$ shown in Figure\ \ref{fig:quiver}, whose vertices are
labelled $0,1,\dots,m$. Recall that the path ring of an arbitrary quiver
$\Gamma$ is the abelian group freely generated by the set of all paths in $\Gamma$, with
multiplication given by composition. Paths of length $l$ in  $\Gamma_m$ correspond
to $(l+1)$-tuples $(i_1|i_2|\dots|i_{l+1})$ of numbers $i_k \in \{0,1,\dots,m\}$ such that the
difference of any two consecutive ones is $\pm 1$. Therefore, the path ring of $\Gamma_m$
is the abelian
group freely generated by such $(l+1)$-tuples for all $l \geq 0$, with the multiplication
\[
(i_1|\dots|i_{l+1})(i'_1|\dots|i'_{l'+1}) = \begin{cases}
 (i_1|\dots|i_l|i'_1|\dots|i'_{l'+1}) & \text{if } i_{l+1} = i'_1, \\
 0 & \text{otherwise.}
\end{cases}
\]
The paths $(i)$ of length zero are mutually orthogonal idempotents, and their
sum is the unit element. We make the path ring of $\Gamma_m$ into a graded ring
by setting $\deg(i) = \deg(i|i+1) = 0$ and $\deg(i+1|i) = 1$ for all $i$, and
extending this in the obvious way (this is \emph{not} the same as the grading
by lengths of paths). Finally, $A_m$ is the quotient of the path ring of
$\Gamma_m$ by the relations
\[
(i-1|i|i+1) = (i+1|i|i-1) = 0, \quad (i|i+1|i) = (i|i-1|i), \quad (0|1|0) = 0
\]
for all $0 < i < m$. These relations are homogeneous with respect to the above grading,
so that $A_m$ is a graded ring. As an
abelian group $A_m$ is free and of finite rank; a basis is given by the $4m+1$ elements $(0),
\dots, (m), (0|1), \dots (m-1|m), (1|0), \dots, (m|m-1), (1|0|1), \dots, (m|m-1|m)$.
\includefigure{quiver}{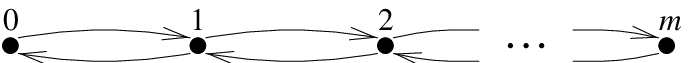}{h}

Let $\Amod$ be the category of finitely generated graded left $A_m$-modules, and $D^b(\Amod)$
its bounded derived category. We will write down explicitly exact self-equivalences ${\mathcal
R}_1, \dots, {\mathcal R}_m$ of $D^b(\Amod)$ which satisfy
\begin{align*}
 \mathcal{R}_i\mathcal{R}_{i+1}\mathcal{R}_i &\iso
 \mathcal{R}_{i+1}\mathcal{R}_i\mathcal{R}_{i+1}
 &&\text{for $1 \leq i < m$, and} \\
 \mathcal{R}_j\mathcal{R}_k &\iso \mathcal{R}_k\mathcal{R}_j
 &&\text{for $|j-k| \geq 2$}.
\end{align*}
These are the defining relations of the braid group $B_{m+1}$, which means that the functors
$\mathcal{R}_i$ generate an action $(\action_\sigma)_{\sigma \in B_{m+1}}$ of $B_{m+1}$ on
$D^b(\Amod)$. The first step in analyzing this action is to look at the induced linear
representation on the Grothendieck group $K(D^b(\Amod)) \cong K(\Amod)$. The category $\Amod$
carries a self-equivalence $\{1\}$ which raises the grading of a module by one (this should
not be confused with the translation functor $[1]$ in the derived category). The action of
$\{1\}$ makes $K(\Amod)$ into a module over $\Z[q,q^{-1}]$. In fact
\[ K(\Amod) \iso \Z[q,q^{-1}] \otimes \Z^{m+1};
\]
a basis is given by the indecomposable projective modules $P_i= A_m(i) \subset A_m$.
The functors ${\mathcal R}_i$ commute with $\{1\}.$
Therefore the induced representation is a homomorphism $B_{m+1} \longrightarrow
GL_{m+1}(\Z[q,q^{-1}])$. A computation shows that this is the well-known Burau representation.
\includefigure{basic-curves}{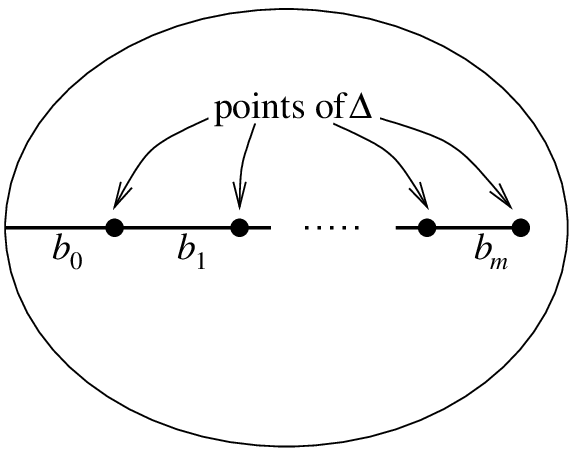}{ht}

At this point, it is helpful to recall the topological meaning of the Burau representation.
Take a closed disc $D$ and a set $\Delta \subset D \setminus \partial D$ of $(m+1)$ marked
points on it. Let $\tD$ be the infinite cyclic cover of $D \setminus \Delta$ whose restriction
to a small positively oriented loop around any point of $\Delta$ is isomorphic to $\R
\longrightarrow S^1$. Fix a point $z$ on $\partial D$, and let $\tilde{z} \subset \tD$ be its
preimage. Let $\Diffeo = \Diff(D,\partial D;\Delta)$ be the group of diffeomorphisms $f$ of
$D$ which satisfy $f|\partial D = \id$ and $f(\Delta) = \Delta$. Any $f \in \Diffeo$ can be
lifted in a unique way to a $\Z$-equivariant diffeomorphism of $\tD$ which acts trivially on
$\tilde{z}$. The Burau representation can be defined as the induced action of $\pi_0(\Diffeo)
\iso B_{m+1}$ on $H_1(\tD, \tilde{z};\Z) \iso \Z[q,q^{-1}] \otimes \Z^{m+1}$. Our first result
will show that the action $(\action_\sigma)$ itself, or rather certain numbers attached to it,
admits a similar topological interpretation. For any $\sigma \in B_{m+1}$ and $0 \leq i,j \leq
m$ consider the bigraded abelian group
\begin{equation} \label{eq:bigraded-group}
\bigoplus_{r_1,r_2 \in \Z} \Hom_{D^b(\Amod)}(P_i,\action_\sigma P_j[r_1]\{-r_2\}).
\end{equation}
Fix a collection of curves $b_0,\dots,b_m$ on $(D,\Delta)$ as in Figure\ \ref{fig:basic-curves};
this determines a preferred isomorphism $\pi_0(\Diffeo) \iso B_{m+1}$. Take some $\sigma \in
B_{m+1}$ and a diffeomorphism $f_\sigma \in \Diffeo$ representing it. The geometric
intersection number $I(b_i,f_\sigma(b_j)) \geq 0$ counts the number of essential intersection
points between the curves $b_i$ and $f_\sigma(b_j)$. The exact definition, which is a slight
variation of the usual one, will be given in the body of the paper.

\begin{theorem} \label{th:gin}
The total rank of the bigraded group \eqref{eq:bigraded-group} is $2\,I(b_i,f_\sigma(b_j))$.
\end{theorem}

We will actually prove a slightly stronger result, which describes \eqref{eq:bigraded-group}
up to isomorphism in terms of {\em bigraded intersection numbers} invented for that purpose.
Theorem \ref{th:gin}, together with standard properties of geometric intersection numbers, has
an interesting implication. Call a group action $(\funct_g)_{g \in G}$ on a category $\cat$
faithful if $\funct_g \not\iso \id_{\cat}$ for all $g \neq e$.

\begin{corollary} \label{th:faithful}
$(\action_\sigma)_{\sigma \in B_{m+1}}$ is faithful for all $m \geq 1$.
\end{corollary}

Note that this cannot be proved by considerations on the level of Grothendieck groups, since
the Burau representation is not faithful (in the usual sense of the word) for $m \gg 0$
\cite{moody91}. Corollary \ref{th:faithful} can be used to derive the faithfulness of certain
braid group actions on categories which occur in algebraic geometry \cite{seidel-thomas99}.

On the symplectic side we will study the Milnor fibres of certain singularities. The
$n$-dimensional singularity of type $(A_m)$, for $m,n \geq 1$, is the singular point $x = 0$
of the hypersurface $\{x_0^2 + \dots + x_{n-1}^2 + x_n^{m+1} = 0\}$ in $\C^{n+1}$. To define
the Milnor fibre one perturbs the defining equation, so as to smooth out the singular point,
and then intersects the outcome with a ball around the origin. In the present case one can
take
\[
M = \{ x_0^2 + \dots + x_{n-1}^2 + h(x_n) = 0, \; |x| \leq 1\}
\]
where $h(z) = z^{m+1} + w_mz^m + \dots + w_1z + w_0 \in \C[z]$ is a polynomial with no
multiple zeros, with $w = (w_0,\dots,w_m) \in \C^{m+1}$ small. Equip $M$ with the restriction
$\o$ of the standard symplectic form on $\C^{n+1}$, and its boundary with the restriction
$\alpha$ of the standard contact one-form on $S^{2n+1}$. Assume that the disc $D$ considered
above is in fact embedded in $\C$ as the subset $\{z \in \C,\; |z|^2 + |h(z)| \leq 1\}$, and
that the set of marked points is $\Delta = h^{-1}(0)$. Then one can associate to any curve $c$
in $(D,\Delta)$ a Lagrangian submanifold $L_c$ of $(M,\o,\alpha)$. We postpone the precise
definition, and just mention that $L_c$ lies over $c$ with respect to the projection to the
$x_n$-variable. In particular, from curves $b_0,\dots,b_m$ as in Figure\ \ref{fig:basic-curves}
one obtains Lagrangian submanifolds $L_0,\dots,L_m$. These are all $n$-spheres except for
$L_0$, which is an $n$-ball with boundary in $\partial M$; they intersect each other transversally,
and
\begin{equation} \label{eq:a-chain}
L_i \cap L_j = \begin{cases} \text{one point} & \text{if
} |i-j| = 1,
\\ \emptyset & \text{if } |i-j| \geq 2. \end{cases}
\end{equation}
The case $m = 3$ and $n = 1$ is shown in Figure\ \ref{fig:vanish}.
\includefigure{vanish}{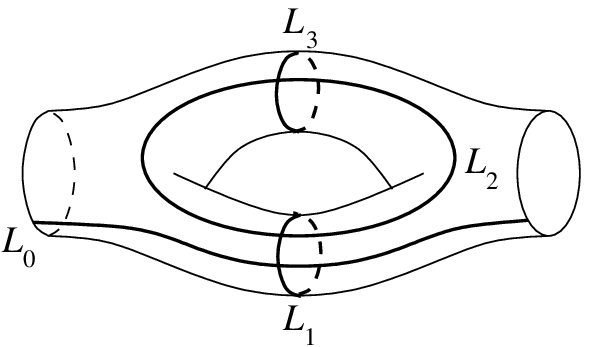}{hb}%

By varying the polynomial $h$ one gets a family of Milnor fibres parametrized by an open
subset $W \subset \C^{m+1}$. Parallel transport in this family, for a suitable class of
connections, defines a {\em symplectic monodromy map}
\[
\rho_s: \pi_1(W,w) \longrightarrow \pi_0(\Symp(M,\partial M,\o)).
\]
The actual construction is slightly more complicated than this rough description may suggest.
Similar maps can be defined for general hypersurface singularities. In our particular case,
there are isomorphisms
\begin{equation} \label{eq:isos}
\pi_1(W,w) \iso \pi_1(\Conf_{m+1}(D \setminus \partial D),\Delta) \iso \pi_0(\Diff(D,\partial
D;\Delta)) \iso B_{m+1};
\end{equation}
the first one comes from a canonical embedding $W \subset \Conf_{m+1}(D \setminus \partial
D)$, the second one is well-known and again canonical, and the third one is determined by the
choice of curves $b_0,\dots,b_m$ as before. Using these isomorphisms, one can consider
$\rho_s$ as a homomorphism from the braid group to $\pi_0(\Symp(M,\partial M,\o))$. We
mention, although this has no relevance for the present paper, that there is a more
direct definition of
$\rho_s$ in terms of generalized Dehn twists along $L_1,\dots,L_m$. This is explained in
\cite[Appendix]{seidel98b} for $n = 2$, and the general case is similar.

Take some $\sigma \in B_{m+1}$ and a map $\phi_\sigma \in \Symp(M,\partial M,\o)$ which
represents $\rho_s(\sigma)$ in the sense of (\ref{eq:isos}).
If $n = 1$ we also assume that $\phi_\sigma$ has a particular
property called $\theta$-exactness, which will be defined in the body of the paper; this
assumption is necessary because in that dimension $H^1(M,\partial M;\R) \neq 0$. The Floer
cohomology $HF(L_i,\phi_\sigma(L_j))$ is a finite-dimensional vector space over $\Z/2$
constructed from the intersection points of $L_i$ and $\phi_\sigma(L_j)$, and it is
independent of the choice of $\phi_\sigma$.

\begin{theorem} \label{th:gin-two}
The dimension of $HF(L_i,\phi_\sigma(L_j))$ is $2\,I(b_i,f_\sigma(b_j))$.
\end{theorem}

Here $f_\sigma \in \Diff(D,\partial D;\Delta)$ is as in Theorem \ref{th:gin}.
There is also a faithfulness result analogous to Corollary \ref{th:faithful}:

\begin{corollary} \label{th:faithful-two}
$\rho_s$ is injective for all $m,n \geq 1$.
\end{corollary}

It is instructive to compare this with the outcome of a purely topological consideration. The
classical geometric monodromy $\rho_g$ is the composition
\[
\xymatrix{&&
 {\pi_0(\Symp(M,\partial M,\o))} \ar[d] \\
 {B_{m+1}} \ar@/^/[rru]^-{\rho_s} \ar[rr]^-{\rho_g}
 && {\pi_0(\Diff(M,\partial M))}
}
\]
where the vertical arrow is induced by inclusion. For $n = 1$ this arrow is an
isomorphism, and the injectivity of $\rho_g$ is an old result of Birman and
Hilden \cite{birman-hilden73} (see Proposition \ref{th:birman-hilden}). The
situation changes drastically for $n = 2$, where Brieskorn's simultaneous
resolution \cite{brieskorn66} implies that $\rho_g$ factors through the
symmetric group $S_{m+1}$. The difference between $\rho_s$ and $\rho_g$ in this
case shows that the symplectic structure of the Milnor fibre contains essential
information about the singularity, which is lost when one considers it only as
a smooth manifold. The situation for $n \geq 3$ is less clear-cut. One can use
unknottedness results of Haefliger \cite{haefliger62}, together with an easy
homotopy computation, to show that $\rho_g$ is not injective for all $m \geq 2$
if $n>2$ is even, and for all $m \geq 3$ if $n \equiv 1 \,\text{mod}\, 4$, $n
\neq 1$ (with more care, these bounds could probably be improved).
%
%

By combining Theorems \ref{th:gin} and \ref{th:gin-two} one sees that the dimension of
$HF(L_i,\phi_\sigma(L_j))$ is equal to the total rank of \eqref{eq:bigraded-group}. One can
refine the comparison by taking into account the graded structure of Floer cohomology, as
follows:

\begin{corollary} \label{th:comparison}
For all $\sigma \in B_{m+1}$, $0 \leq i,j \leq m$, and $r \in \Z$,
\begin{equation} \label{eq:comparison-eq}
 HF^r(L_i,\phi_\sigma(L_j)) \iso \bigoplus_{r_1+n r_2 = r}
 \Hom_{D^b(\Amod)}(P_i,\action_\sigma P_j[r_1]\{-r_2\}) \otimes \Z/2.
\end{equation}
\end{corollary}

Strictly speaking, the grading in Floer theory is canonical only up to an
overall shift (in fact, for $n = 1$ there is an even bigger ambiguity), and the
isomorphism holds for a particular choice. One can cure this problem by using
Kontsevich's idea of graded Lagrangian submanifolds. We will not state the
corresponding improved version of Corollary \ref{th:comparison} explicitly, but
it can easily be extracted from our proof. Note that if one looks at the graded
groups $HF^*(L_i,\phi_\sigma(L_j))$ in isolation, they appear to depend on $n$
in a complicated way; whereas the connection with the bigraded group
\eqref{eq:bigraded-group} makes this dependence quite transparent.

One can see this collapsing of the bigrading on a much simpler level, which
corresponds to passing to Euler characteristics on both sides in Corollary
\ref{th:comparison}. Namely, the representation of the braid group on
$H_n(M,\partial L_0;\Z) \iso \Z^{m+1}$ induced by $\rho_s$ is the
specialization $q = (-1)^n$ of the Burau representation.

While our proof of Corollary \ref{th:comparison} is by dimension counting, it
is natural to ask whether there are canonical homomorphisms between the two
graded groups; these should be compatible with the product structures existing
on both sides. Unfortunately, the techniques used here are not really suitable
for that; while an inspection of our argument suggests a possible definition of
a homomorphism in \eqref{eq:comparison-eq}, showing that it is independent of
the various choices involved would be extremely complicated. And there is no
easy way to adapt the proof of Theorem \ref{th:gin-two} to include a
description of the product on Floer cohomology. It seems that to overcome these
problems, one would have to use a more abstract approach, in line with the
general picture of Section \ref{subsec:generalities}. While the details of this
are far from clear, an outline can easily be provided, and we will do that now.

It is convenient to modify the algebraic setup slightly. Let $A_{m,n}$ be the
graded algebra obtained from $A_m$ by multiplying the grading with $n$ and
changing the coefficients to $\Z/2$; regard this as a differential graded
algebra $\dga_{m,n}$ with zero differential, and take the derived category of
differential graded modules $D(\dga_{m,n})$. From the way in which
$L_0,\dots,L_m \subset M$ intersect each other, and standard properties of
Floer cohomology, one can see that
\[
A_{m,n} \iso \bigoplus_{i,j = 0}^m HF^*(L_i,L_j)
\]
(actually, this was the starting point of our work). It seems reasonable to
assume that the $A_\infty$-algebra of Floer cochain complexes underlying the
right hand side should be formal (see \cite{seidel-thomas99} for much more
about this), which would make its derived category of $A_\infty$-modules
equivalent to $D(\dga_{m,n})$. As mentioned in Section
\ref{subsec:generalities}, the general features of the construction should
imply that one has an exact functor $D^b\Lag\mo \rightarrow D(\dga_{m,n})$,
which would be full and faithful on the triangulated subcategory of
$D^b\Lag\mo$ generated by $L_0,\dots,L_m$. $\rho_s$ defines a braid group
action on $\Lag\mo$, and there should be an induced action on $D^b\Lag\mo$. On
the other hand, one can define a braid group action on $D(\dga_{m,n})$ in a
purely algebraic way, in the same way as we have done for the closely related
category $D^b(\Amod)$. It is natural to conjecture that the exact functor
should intertwine the two actions. This would imply Corollary
\ref{th:comparison}, in a stronger form which repairs the deficiencies
mentioned above.

Finally, the reader may be wondering whether the objects which we have considered can be used
to define invariants of two-knots. The category $D^b(\Amod)$ is certainly not rich enough for
this. However, it has been conjectured \cite{bernstein-frenkel-khovanov98} that one can
construct invariants of two-knots from a certain family of categories, of which $D^b(\Amod)$
is the simplest nontrivial example. The homology groups for classical knots introduced in
\cite{khovanov98} are part of this program.

\emph{Note.}
During the final stage of work on this paper we learned of a preprint \cite{RZ} by Rouquier
and Zimmermann. They construct a braid group action in the derived category of modules over
multiplicity one Brauer tree algebras, which are close relatives of $A_m$ (the idempotent
$(0)$ of $A_m$ is set to zero; these smaller algebras also appear in \cite{seidel-thomas99}
and \cite{huerfano-khovanov}).
Our action essentially coincides with the one of Rouquier and Zimmermann.

{\em Acknowledgements.} M.K.\ would like to thank Maxim Vybornov who introduced him to the
algebras $A_m$ by pointing out that they describe perverse sheaves on projective spaces. The
construction of the Lagrangian submanifolds $L_c$ was suggested to P.S.\ by Simon Donaldson.
Most of this work was done at the Institute for Advanced Study, during the year 97/98 and two
shorter visits later on. We are very grateful to the Institute, and to Robert MacPherson in
particular, for making this collaboration possible. M.K.\ was supported by NSF grants DMS
96-27351 and DMS 97-29992 and, later on, by the University of
California at Davis. P.S.\ was supported by NSF grant DMS-9304580 and by the Institut
Universitaire de France.

%
%
%
%
%
%
%

\section {The braid group action on $\db$}
\label{der-cat}

\subsection{The category of $A_m$-modules}
\label{Am-mod}

Fix $m \geq 1$. We recall some of the objects introduced in Section
\ref{subsec:results}:
\begin{itemize}
\item the graded ring $A_m$. In future all left modules,
right modules, or bimodules over $A_m$ are understood to be graded;

\item the category $\Amod$, whose objects are finitely generated
$A_m$-modules, and whose morphisms are grading-preserving module homomorphisms.
In future, we call the objects of $\Amod$ just
$A_m$-modules, and write $\Hom(M,N)$ instead of $\Hom_{\Amod}(M,N)$;
Throughout the paper, all homomorphisms are assumed grading-preserving unless
otherwise mentioned.

\item the self-equivalence $\{1\}$ of $\Amod$ which shifts the
grading upwards by one;

\item the $A_m$-modules $P_i = A_m(i).$ These are indecomposable and projective.
In fact, if one considers $A_m$ as a left
module over itself, it splits into the direct sum $\bigoplus_{i=0}^m
P_i$. Conversely, any indecomposable projective $A_m$-module
is isomorphic to $P_i\{k\}$ for some $i$ and $k$.
\end{itemize}

Similarly, if one considers $A_m$ as a right module over itself,
it splits into the direct sum of the indecomposable projective right
modules $\iP = (i)A_m$.

\begin{prop}\label{prop:finite-hom-dim} $A_m$ has finite homological dimension.
\end{prop}
\proof
Introduce  $A_m$-modules $\simp_0,\simp_1,\dots ,\simp_m$ as
follows. As a graded abelian group $\simp_i$ is isomorphic to $\Z,$
placed in degree $0.$ The idempotent $(i)$ acts as the identity on $\simp_i,$
and all other paths in the quiver act as zero. Any $A_m$-module
has a finite length composition series with subsequent quotients isomorphic,
up to shifts in the grading, to $\simp_i$ or $\simp_i/p\simp_i$ for various
$i$ and primes $p.$ Moreover, $\simp_i/p\simp_i$ has a resolution by
modules $\simp_i$:
  \begin{equation*}
  0 \lra \simp_i \lra \simp_i \lra \simp_i/p\simp_i \lra 0.
  \end{equation*}
Therefore, the proposition will follow if we construct a finite length
projective resolution of $\simp_i,$ for each $i.$
Consider the commutative diagram
\[
 \begin{array}{ccccccc}
     P_m\{m-i\} & \to & \dots & \to & P_{i+1}\{1\} & \to & P_i \\
     \uparrow &  & &    & \uparrow  &  &  \uparrow \\
     P_{m-1}\{m-i\} & \to & \dots & \to & P_i\{1\} & \to & P_{i-1} \\
     \uparrow &  &    &  &  \uparrow &  & \uparrow \\
     \dots    &  & \dots   &  &   \dots   &  & \dots \\
     \uparrow &  &    &  &  \uparrow &  & \uparrow \\
     P_{m-i}\{m-i\} & \to   & \dots  & \to & P_1\{1\}  &  \to  & P_0 \\
     \uparrow &   &  \dots  &   &  \uparrow &  &  \\
     \dots   & \to   &\dots    & \to  & P_0\{1\}  &   &   \\
     \dots   &  \to  &  \dots  &  &            &   &  \\
      \uparrow &  &  &  &  &  &    \\
      P_0\{m-i\} & & &  && &     \\
  \end{array}
\]
of projective $A_m$-modules. The maps go from $P_j$ to
$P_{j\pm 1}$ for various $j$ and are given by right multiplications by
$(j|j\pm 1),$ which maps $P_j\subset A_m$ to $P_{j\pm 1}\subset A_m.$
The commutativity of the diagram follows from the relations
$(j|j-1|j)=(j|j+1|j).$
Moreover, every row and every column is a complex, since
$(j|j+1|j+2)=(j|j-1|j-2)=0,$ so that, after inserting minus signs in
appropriate places, the diagram above becomes a bicomplex. The
associated total complex is a finite length projective resolution
of the module $\simp_i$ (use that every column is
acyclic except in the uppermost component).
 $\square$

\subsection{The functors $\U_i$.} \label{functors-ui}

For $1\le i\le m$ let $P_i\otimes\hsm \iP$ be the $A_m$-bimodule obtained by
tensoring $P_i$ and $\iP$ over $\Z.$ Let $\U_i:\amod\to \amod$
be the functor of tensoring with this bimodule,
  \begin{equation} \label{uu:equation}
    \U_i(M)= P_i \oo \hsm _iP\oo_{A_m} M.
  \end{equation}
This functor is exact (since $P_i\oo \hsm _iP$ is right
projective) and takes projective modules to projective modules\footnote{
$i=0$ case is excluded, for  although  $\U_0$ can be defined, it does not
satisfy the isomorphisms in Theorem~\ref{TL-relations}.}.

 \begin{theorem} \label{TL-relations} There are functor isomorphisms
 \begin{eqnarray} \label{three-U}
  \U_i\U_{i+ 1} \U_i & \cong & \U_i \{ 1\}
  \hspace{0.15in}\mbox{ for }1\le i\le m-1, \\
 \label{iss-1}
  \U_i\U_{i- 1} \U_i & \cong & \U_i \{ 1\}
  \hspace{0.15in}\mbox{ for }2\le i\le m, \\
  \U_i ^ 2 & \cong & \U_i \oplus \U_i\{ 1\}\hspace{0.1in}
  \mbox{ for }1\le i\le m,   \\
  \U_i \U_j & = & 0 \hspace{0.3in}\mbox{ for }|i-j|> 1. \label{iss-3}
 \end{eqnarray}
 \end{theorem}
\proof
 We will consider only (\ref{three-U}); the other isomorphisms
are proved by similar arguments.
The functor on the left hand side of (\ref{three-U})
is given by tensoring with the bimodule
$P_i \oo \hsm _i P\oo_{A_m} P_{i+1} \oo \hsm _{i+1}P \oo_{A_m}
P_i \oo \hsm _iP.$ The graded abelian groups $_i P\oo_{A_m} P_{i+1}$ and
$_{i+1}P \oo_{A_m} P_i $ are free cyclic, one in degree $0$, the other in
degree $1.$ Thus $\U_i\U_{i+ 1} \U_i$ is isomorphic to the composition of
tensoring with $P_i \oo \hsm _iP$ and shifting the grading by $\{ 1\}.$ This is exactly
the functor on the right hand side of (\ref{three-U}). $\square$

\begin{remark} The functors $\U_i$ provide a functor realization of the Temperley-Lieb
algebra (compare to Section~\ref{tl-algebra} and \cite[Section
4.1.2]{bernstein-frenkel-khovanov98}).
\end{remark}

\subsection{Categories of complexes.} \label{cats-compl}
For an abelian category $\ab$ denote by $P(\ab)$ the homotopy category of
bounded complexes of projective objects in $\ab.$ An object $M$ of $P(\ab)$
is a bounded complex of projective objects:
  \begin{equation*}
  M=(M^i,\partial^i), \hsm \hsm \partial^i:M^i\to M^{i+1}, \hsm \hsm
  \partial^{i+1}\partial^i=0.
  \end{equation*}
A morphism $f$ from $M$ to $N$ is a collection of maps $f^i: M^i \to N^i$
that intertwine differentials in $M$ and $N,$
 $f^{i+1}\partial_M^i= \partial_N^i f^i$ for all $i.$
Two morphisms $f,g$  are equal (homotopic)
if there are maps $h^i\in \Hom_{\ab}(M^i,N^{i-1})$
such that $f^i-g^i = h^{i+1} \partial_M^i + \partial_N^{i-1} h^i$ for all $i.$

As usual, $[k]$ denotes the self-equivalence of
$P({\ab})$ which shifts a complex $k$ degrees to the left:
$M[k]^i = M^{k+i}, \hsm \partial_{M[k]} = (-1)^{k} \partial_M.$
Given a map of complexes $f: M \to N,$ the cone of $f$
is the complex $C(f)=M[1]\oplus N$ with the
differential $\partial_{C(f)} (x,y)= (-\partial_Mx, f(x) + \partial_Ny).$
We refer the reader
to \cite{gm} for more information about homotopy and derived categories.

If $\ab$ has finite homological dimension, then
$P(\ab)$ is equivalent to the bounded derived category $D^b(\ab)$ of
$\ab.$   This applies to ${\ab}= \Amod$ by Proposition
\ref{prop:finite-hom-dim}. We will use $P(\Amod)$ instead of
$D^b(\Amod)$ throughout, and denote it by ${\mathcal C}_m = P(\Amod).$

The shift automorphism $\{k\}$ of $\Amod$ extends in the obvious way
to an automorphism of ${\mathcal C}_m$ which we also denote by $\{k\}$.
Note that $\{k\}$ and $[k]$ are quite different. In fact, if we forget
about the module structure and differential, objects of ${\mathcal
C}_m$ are bigraded abelian groups, so that it is natural to have
two shift-like functors. Correspondingly, for any two such objects
$M,N$ there is a bigraded morphism group (of finite total rank)
\[
\bigoplus_{r_1,r_2 \in \Z} \Hom(M,N[r_1]\{-r_2\}).
\]

Let $R=(R^j,\partial^j)$ be a bounded complex of $A_m$-bimodules which
are left and right projective, i.e. each $R^j$ is a projective
left $A_m$-module and a projective right $A_m$-module.
One associates to  $R$ the  functor $\cR: \Ccat\to \Ccat$ of tensoring with $R:$
  \begin{equation*}
    \cR(M)= R\oo_{A_m} M.
  \end{equation*}

  \vspace{0.1in}

\subsection{The functors $\cR_i$} \label{comp-cat}
Define homomorphisms $\beta_i,\gamma_i, 1\le i \le m,$ of
$A_m$-bimodules
  \begin{eqnarray*}
  \beta_i & : & P_i \oo \hspace{0.03in} _iP\lra  A_m, \\
  \gamma_i & : & A_m \lra  P_i \oo \hspace{0.03in} _iP \{ -1\},
  \end{eqnarray*}
by
  \begin{eqnarray}
  \beta_i ((i)\oo (i)) & = & (i) \label{beta}, \\
  \gamma_i (1)& = &  (i-1|i)\oo (i|i-1) + (i+1|i)\oo (i|i+1)\notag \\
            &   & + (i)\oo (i|i-1|i) + (i|i-1|i)\oo (i).
  \label{gamma}
 \end{eqnarray}
(when $i=m$ we omit the term $(i+1|i)\oo (i|i+1)$ from the sum for
 $\gamma_i(1).$)

These homomorphisms induce natural transformations, also denoted
$\beta_i,\gamma_i,$ between the functor $\U_i$, suitably shifted, and
the identity functor:
  \begin{eqnarray*}
  \beta_i & : & \U_i \lra \mbox{Id}, \\
  \gamma_i & : & \mbox{Id} \lra \U_i \{ -1\}.
  \end{eqnarray*}

Let ${R}_i$ be the complex of bimodules
  \begin{equation*}
  R_i = \{ 0\lra P_i \oo \hspace{0.03in} _iP \stackrel{\beta_i}{\lra} A_m
  \lra 0 \}
  \end{equation*}
with $A_m$ in degree $0.$
Denote by ${R}'_i$ the complex of bimodules
  \begin{equation*}
  R_i'= \{ 0\lra A_m \stackrel{\gamma_i}{\lra}
  P_i \oo \hspace{0.03in} _iP \{ -1\} \lra 0 \}
  \end{equation*}
with $A_m$ in degree $0.$
Next, we introduce functors $\cR_i$ and $\cR_i'$ of tensoring with the complexes
${R}_i$ and ${R}'_i$ respectively:
  \begin{equation*}
  \cR_i(M)= {R}_i \oo_{A_m}M, \hspace{0.2in}
  \cR'_i(M)= {R}'_i \oo_{A_m}M, \hspace{0.2in} M\in \mathrm{Ob}(\Ccat).
  \end{equation*}

The functor  $\cR_i$ can be viewed as the cone of $\beta_i,$ and
$\cR_i'$ as the cone of $-\gamma_i,$ shifted by $[-1].$

  \begin{prop}
     \label{inverse}  The functors $\cR_i$ and $\cR_i'$ are mutually inverse
      equivalences of categories, i.e., there are functor isomorphisms
     \begin{equation*}
      \cR_i \cR_i' \cong \mathrm{Id} \cong \cR_i' \cR_i.
     \end{equation*}
  \end{prop}

\proof The functor $\cR_i\cR_i'$ is given by
tensoring with the complex ${R}_i\oo_{A_m} {R}_i'$ of $A_m$-bimodules. We
write down this complex explicitly below. First consider a  commutative
square of $A_m$-bimodules and $A_m$-bimodule homomorphisms
   \begin{equation}
   \begin{CD}
   \label{CD-simple}
   P_i \oo \hsm _iP @>{\beta_i}>> A_m   \\
   @VV{\tau}V      @VV{\gamma_i}V   \\
   P_i \oo Q\oo \hsm _iP \{ -1\}
   @>{\delta}>> P_i \oo \hsm _iP\{ -1\}.
   \end{CD}
   \end{equation}
Here $Q$ is the graded abelian group $_iP\oo_{A_m}P_i.$ It is a free group of
rank $2$ with a basis
$u_1= (i)\oo (i) , u_2=(i|i-1|i)\oo (i).$
The maps $\tau$ and $\delta$ are given by
  \begin{eqnarray*}
    \tau (x\oo y) & = &  x\oo u_1 \oo (i|i-1|i)y + x \oo u_2 \oo y, \\
    \delta( x\oo u_1 \oo y) & = &  x \oo y,  \\
    \delta( x\oo u_2 \oo y) & = & x (i|i-1|i) \oo y.
   \end{eqnarray*}
Define a map $\xi: A_m \to P_i \oo Q\oo \hsm _iP\{ -1\}$ as the composition
of $\gamma_i$ and the map $P_i \oo \hsm _iP\{ -1\} \to
  P_i \oo Q\oo \hsm _iP\{ -1\}$ given by $x\oo y \to x\oo u_1 \oo y.$
Notice that $\delta\xi= \gamma_i.$

We make the diagram (\ref{CD-simple}) anticommutative by putting
a minus sign in front of $\delta.$ Denote by $N= (N^j,\partial^j)$ the
complex of bimodules, associated to this anticommutative diagram.
It has $3$ nonzero components
  \begin{eqnarray*}
    N^{-1} & = &  P_i \oo \hsm _iP,   \\
    N^0 & = &    A_m  \oplus  (P_i \oo Q \oo \hsm _iP \{ -1\}), \\
    N^1 & = &    P_i \oo \hsm _iP\{ -1\},
  \end{eqnarray*}
and the differential $\partial^j : N^j \to N^{j+1}$ reads
  \begin{eqnarray*}
    \partial^{-1} & = & \beta_i + \tau, \\
    \partial^0  & = & (\gamma_i,- \delta), \\
    \partial^j  & = & 0 \hspace{0.2in} \mbox{ for }j\not= -1,0.
  \end{eqnarray*}

It is easy to verify that $N$ is isomorphic to the complex
${R}_i\oo_{A_m} {R}'_i.$ Thus, we want to
show that, in the homotopy category of complexes of $A_m$-modules,
tensoring with $N$ is equivalent to the identity functor, i.e., to
the functor of tensoring with $A_m,$ considered as an $A_m$-bimodule.
We will decompose $N$ as a direct sum of a complex isomorphic to $A_m$
and an acyclic complex. Namely, $N$ splits into  a direct sum of $3$
complexes, $N = T_{-1} \oplus T_0 \oplus T_1$ where
  \begin{eqnarray*}
    T_{-1}& = &  \{  0 \to N^{-1} \to  \partial^{-1}(N^{-1}) \to 0 \},  \\
    T_0   & = &   \{ (a,\xi (a))| a \in A_m \subset N^0 \},   \\
    T_ 1  & = &  \{  0 \to  P_i \oo u_1 \oo \hsm _i P \to N^1 \to 0\}.
  \end{eqnarray*}
$T_{-1}$ is a subcomplex of $N$ generated by $N^{-1}.$
Notice that $\partial^{-1}: N^{-1}\to N^{0}$ is injective
(since $\tau,$ a composition of $d^{-1}$ and a projection, is injective)
and, therefore,  $T_{-1}$ is acyclic.
The complex $T_1$ is nonzero only in degrees $0$ and $1.$ Its
 degree $0$ component is $P_i\oo u_1 \oo \hsm _iP,$ considered
as a subbimodule of $N^0,$ and the degree $1$ component is $N^1.$
The differential of $T_1$ induces an isomorphism between degree $0$ and
degree $1$ components, and, thus, $T_1$ is acyclic.
Finally, the complex $T_0$ is
concentrated in horizontal degree $0,$ and consists of the set of
pairs $(a,\xi(a))\in N^0,$ where $a\in A_m.$ Clearly, $T_0$
is isomorphic to the complex which has $A_m$ in horizontal
degree $0$ and is trivial in all other degrees. Therefore, tensoring with
$T_0$ is equivalent to the identity functor, while tensoring
with $T_{\pm 1}$ is equivalent to tensoring with the zero complex.

It now follows that
$\cR_i \cR_i'$ is isomorphic to the identity functor. A similar
computation gives an  isomorphism of $\cR_i'\cR_i$ and the identity
functor. $\square$

We have seen that the functor $\cR_i$ is invertible, with the inverse functor
isomorphic to $\cR_i'.$ In view of this and for future convenience
we  will denote $\cR_i'$ by $\cR_i^{-1}.$

  \begin{theorem}
   \label{theorem-braid}
    There are functor isomorphisms
    \begin{eqnarray}
    \cR_i \cR_{i+1} \cR_i & \cong & \cR_{i+1} \cR_i \cR_{i+1},
    \label{braid-3} \\
    \cR_i \cR_j  & \cong & \cR_j \cR_i \hspace{0.1in}\mbox{ for }|i-j|>1.
    \label{braid-2}
   \end{eqnarray}
  \end{theorem}

\proof
The commutativity relation (\ref{braid-2}) follows at once
from the obvious isomorphism ${R}_i\oo_{A_m} {R}_j\cong {R}_j\oo_{A_m} {R}_i$
of complexes of bimodules. Indeed, both complexes are isomorphic to the complex

  \begin{equation*}
   0 \lra (P_i \oo \hsm _iP)\oplus (P_j \oo \hsm _jP)
   \stackrel{\beta_i+\beta_j}{\lra } A_m \lra 0.
  \end{equation*}
Therefore, we concentrate on proving (\ref{braid-3}), which is equivalent to
  \begin{equation}
   \label{braid-3a}
   \cR_{i+1}^{-1}\cR_i \cR_{i+1} \cong \cR_i\cR_{i+1}\cR_i^{-1}.
  \end{equation}
The functor on the left side is given by the tensor product with the complex
of bimodules $R'_{i+1}\oo_{A_m} R_i \oo_{A_m} R_{i+1},$ isomorphic to  the cone
of the map of complexes
  \begin{equation}
   \label{onemore}
   R'_{i+1} \oo_{A_m} P_i \oo \hsm _iP \oo_{A_m} R_{i+1} \lra R'_{i+1}
   \oo_{A_m} R_{i+1}.
  \end{equation}
From the proof of Proposition~\ref{inverse} we know that
$R'_{i+1}\oo_{A_m} R_{i+1}$ is isomorphic to the direct sum of an acyclic
complex and the complex $\{ \dots \lra 0 \lra A_m\lra 0\lra \dots \}.$
Factoring out by the acyclic subcomplex we get a complex of bimodules
(quasi-isomorphic to (\ref{onemore}))
  \begin{equation}
  \label{slide1}
   R'_{i+1} \oo_{A_m} P_i \oo \hsm _iP \oo_{A_m} R_{i+1} \lra  A_m.
  \end{equation}
Notice that $R'_{i+1}\oo_{A_m} P_i$
is isomorphic to the complex of left modules
$\{ 0 \lra P_i \lra P_{i+1} \lra 0\}$ where $P_i$ sits in the horizontal
degree $0;$ and $_iP\oo_{A_m} R_{i+1}$ is isomorphic to the complex of right
modules $\{ 0 \lra \hsm _{i+1}P\lra \hsm _iP \lra 0\},$ where $_iP$
sits in degree $0.$ Therefore, the map (\ref{slide1}) has the form

  \begin{equation*}
   \begin{CD}
   0 @>>> P_i \oo\hsm  _{i+1}P   @>{\partial^{-1}}>>
   {\begin{array}{c}
   (P_i \oo \hsm _iP) \\
   \oplus   \\
   (P_{i+1}\oo \hsm _{i+1}P)
   \end{array}}
   @>{\partial^0}>> P_{i+1}\oo \hsm _iP @>>> 0   \\
   @VVV            @VVV         @V{e}VV       @VVV  @VVV  \\
   0   @>>> 0 @>>> A_m @>>> 0 @>>> 0
   \end{CD}
  \end{equation*}
for some map $e$ of bimodules. $e$ is completely
described by two integers $a_1$ and $a_2,$ where
\begin{equation*}
e((i)\oo (i))=a_1 , \hspace{0.1in} e((i+1)\oo (i+1))= a_2.
\end{equation*}
 The condition $e \partial^{-1}=0$ implies $a_1+a_2=0.$
Moreover, if $a_1\not= \pm 1,$ we can quickly come to a
contradiction by picking a prime divisor $p$ of $a_1,$ reducing all
modules and functors mod $p$ and concluding that the
functor $\cR_{i+1}^{-1}\cR_i \cR_{i+1}$ (invertible over any
base field) decomposes as a nontrivial direct sum in characteristic $p.$
Thus, $a_1=1$ or $-1.$ Multiplying each element of the bimodule
$A_m$ by $-1,$ if necessary, we can assume that $a_1=1$ and $a_2=-1.$

The right hand side functor of (\ref{braid-3a}) is given by tensoring
with the complex of bimodules $R_i \oo_{A_m} R_{i+1} \oo_{A_m} R'_i.$
Since $R_i \oo_{A_m} P_{i+1}$ is isomorphic to $\{ 0 \to P_i \to P_{i+1}\to
0 \}$ and $_{i+1}P\oo_{A_m} R'_i$ is isomorphic to
$\{ 0 \to \hsm _{i+1}P \to \hsm _iP\to 0\},$
the right hand side of (\ref{braid-3a})
is isomorphic to the functor of tensoring with the  cone
 \begin{equation*}
   \begin{CD}
   0 @>>> P_i \oo \hsm _{i+1}P   @>{\partial^{-1}}>>
   {\begin{array}{c}
   (P_i \oo \hsm _iP) \\
    \oplus \\
   (P_{i+1}\oo \hsm _{i+1}P)
   \end{array}}
   @>{\partial^0}>> P_{i+1}\oo \hsm _iP @>>> 0   \\
   @VVV            @VVV         @V{f}VV       @VVV  @VVV  \\
   0   @>>> 0 @>>> A_m @>>> 0 @>>> 0
   \end{CD}
  \end{equation*}
for some bimodule map $f.$ Using the
 same argument as for $e$ we can assume that $f((i)\oo (i))= 1$ and
$f((i+1)\oo (i+1))=-1.$ Therefore, $e=f$ and the two functors in
(\ref{braid-3a}) are isomorphic. $\square$

\begin{definition} A weak action of a group $G$ on a category $\mathcal Q$
is a choice of functors ${\mathcal F}_g: {\mathcal Q}\to {\mathcal Q}$ for
each $g\in G$ such that
${\mathcal F}_1$ is the identity functor and functors ${\mathcal F}_{fg}$
and ${\mathcal F}_f{\mathcal F}_g$ are isomorphic for all $f,g\in G.$
\end{definition}

To each element $\sigma$ of the braid group $B_{m+1}$ we associate
a complex ${R}_{\sigma}$ of graded $A_m$-bimodules as follows.
To $\sigma_i$ associate the complex $R_i$ and to $\sigma_i^{-1}$ the complex
$R'_i.$ In general,
fix a decomposition of $\sigma$ as a product of generators $\sigma_i$ of
the braid group and their inverses, $\sigma = \tau_1 \dots \tau_k.$
Then define ${R}_{\sigma}$ as the tensor product of corresponding complexes:
  \begin{equation*}
   {R}_{\sigma} \stackrel{\mbox{\scriptsize{def}}}{=} {R}_{\tau_1} \oo_{A_m}
   {R}_{\tau_2} \oo \dots \oo_{A_m} {R}_{\tau_k}.
  \end{equation*}

Finally, define $\cR_{\sigma}: \Ccat\to \Ccat$ as the functor of tensoring
with the bimodule ${R}_{\sigma}$:
  \begin{equation*}
    \cR_{\sigma}(M)\stackrel{\mbox{\scriptsize{def}}}{=}
    {R}_{\sigma}\oo_{A_m} M, \hspace{0.3in}
    M\in \mathrm{Ob}(\Ccat).
  \end{equation*}
Proposition \ref{inverse} and Theorem \ref{theorem-braid}
together are equivalent to

\begin{prop} The functors $\cR_{\sigma}$ for $\sigma\in B_{m+1}$ define a weak
action of the braid group on the category $\Ccat.$
\end{prop}

A weak action of $G$ on a category $\mathcal Q$ is called {\it faithful}
if for any $g\in G, g\not= 1,$ the functor ${\mathcal F}_g$ is not isomorphic
to the identity functor.
We will prove in Section~\ref{normal-complex} that the above braid group action is
faithful.

\subsection{Miscellaneous} \label{miss}

In this subsection we link the Burau representation with the braid group
action in the derived category $D^b(\amod).$ We also relate $A_m$ and
highest weight categories.

\subsubsection{The Grothendieck group of $\amod$ and the Burau representation}

Recall that the Grothendieck group $K(\amod)$ of the category $\amod$ is the
group formed by virtual graded $A_m$-modules $[M]-[N],$ with a suitable
equivalence relation. The functor $\{ 1\} $ makes $K(\amod)$ into a module
over $\Z[q,q^{-1}].$ In fact, it is a free $\Z[q,q^{-1}]$-module of rank
$m+1.$ As a basis one can take the images of projective modules $P_i.$
The functor $\U_i$
is exact and induces a $\Z[q,q^{-1}]$-linear map $[\U_i]$  on $K(\amod).$
The functor $\cR_i$ is the cone of $\beta_i: \U_i \to \mbox{Id}$ and so it
induces a linear map $[\cR_i]= [\mbox{Id}]-[\U_i]$ on
$K(\amod).$ Since the $\cR_i$ satisfy braid group relations, the
operators $[\cR_i]$ define a representation of the braid group on $K(\amod).$
The generators of the braid group act as follows:
  \begin{eqnarray*}
  \sigma_i [P_i] & = & - q [P_i],  \\
  \sigma_i [P_{i+1}] & = &  [P_{i+1}]-[P_i],  \\
  \sigma_i [P_{i-1}] & = &  [P_{i-1}]-q[P_i],   \\
  \sigma_i [P_j]     & = &  [P_j] \mbox{ for }|i-j|>1.
  \end{eqnarray*}

The braid group action on the Grothendieck group preserves the
$\Z[q,q^{-1}]$-submodule of $K(\amod)$ spanned by
$[P_1], [P_2], \dots , [P_m].$ Denote this submodule by $K'.$ Comparing with
the Burau representation (see \cite{birman}) we obtain

\begin{prop} The representation on $K(\amod)$ is equivalent to the Burau
representation, and its restriction to $K'$ is equivalent to the reduced
Burau representation.
\end{prop}

\subsubsection{Ring $A_m$ and highest weight categories}

Let $\sl_{m+1}$ be the Lie algebra of traceless $(m+1)\times (m+1)$-matrices
with complex coefficients. Let $\mf{h}$ be the subalgebra of traceless
diagonal matrices and $\mf{p}$ be the subalgebra of matrices with zeros in
all off-diagonal entries of the first column. Let $Z$ be the center of the
universal enveloping algebra $U(\sl_{m+1})$ and denote by $Z_0$ the maximal
ideal of $Z$ which is the kernel of the augmentation
homomorphism $Z\to \C.$

Denote by $\Oo'_m$ the category whose objects are finitely-generated
$U(\sl_{m+1})$-modules that are $U(\mf{h})$-diagonalizable,
$U(\mf{p})$-locally finite and annihilated by some power of the maximal
central ideal $Z_0.$ Morphisms are $U(\sl_{m+1})$-module maps
(a  $U(\sl_{m+1})$-module $M$ is said to be $U(\mf{p})$-locally finite if
for any $x\in M$ the vector space $U(\mf{p})x$ is finite-dimensional).

By a base change we can make the ring $A_m$ into an algebra over the
field of complex numbers, $A_m^{\C}\stackrel{\mbox{\scriptsize{def}}}{=}
A_m \oo_{\Z} \C.$ Let $A_m^{\C}$-mod denote the category of
finitely-generated left $A_m^{\C}$-modules (no grading this time).
We learned the following result and its proof from Maxim Vybornov:

\begin{prop}
\label{highest-weight}
The categories $\Oo'_m$ and $A_m^{\C}\mathrm{-mod}$ are equivalent.
\end{prop}

\proof
The category $\Oo'_m$ consists of modules from a regular
block of the Bernstein-Gelfand-Gelfand category $\Oo$ for $\sl_{m+1}$
which are locally $U(\mf{p})$-finite, where $\mf{p}$ is the parabolic
subalgebra described above. Let $Q_0, \dots, Q_m$ be indecomposable projective
modules in $\Oo'_m,$ one for each isomorphism class of indecomposable
projectives. Form the projective module $Q= Q_0 \oplus \dots \oplus Q_m.$
The category $\Oo'_m$ is equivalent to the category of modules over
the algebra $E=\mbox{End}_{\Oo'_m}(Q).$ We claim that $E$ is
isomorphic to $A_m^{\C}.$ Indeed, from \cite{bgs} we know that $E$ has a
structure
of a $\Z_+$-graded algebra, $E= \oplusop{i\ge 0} E_0,$ such that  $E$  is
multiplicatively generated by $E_1$ over $E_0$ with quadratic defining
relations. Moreover, $E$ is Koszul and the category of representations of
its Koszul dual algebra $E^!$ is equivalent to a certain singular block of
the category $\Oo$ for $\sl_{m+1}$ (this is just a special case of the
Beilinson-Ginzburg-Soergel parabolic-singular duality \cite{bgs}).
The algebra $E^!,$ describing this singular block, was written down by
Irving (see \cite[Section 6.5]{ir}).
Recall that $A_m$ is the quotient ring of
the path ring of the quiver associated to the graph $\Gamma_m$ by certain
relations among paths of length $2.$ The algebra $E^!$ is the quotient
algebra of the path algebra of the quiver with the same graph $\Gamma_m$ but
different set of quadratic relations:
  \begin{equation*}
  (i|i-1|i) = (i|i+1|i), \hspace{0.3in}
  (m|m-1|m) = 0.
  \end{equation*}
An easy computation establishes that the Koszul dual to $E^!$ is
isomorphic to $A_m^{\C}.$ Therefore, algebras $E$ and $A_m^{\C}$ are
isomorphic. $\square$

\emph{Remark:} The grading which makes $A_m$ into a quadratic and Koszul
algebra is the grading by the length of paths. This is different from
the grading on $A_m$ defined in Section 1b and used throughout this paper.

\vsp

An alternative proof of Proposition~\ref{highest-weight}, which also
engages Koszul duality, goes as follows.
The category $\Oo'_m$ is equivalent to the category of perverse sheaves on
the projective space $\mathbb{P}^m,$ which are locally constant on each
stratum $\mathbb{P}^i\setminus \mathbb{P}^{i-1}$ of the stratification of
$\mathbb{P}^m$ by a chain of projective spaces
$\mathbb{P}^0 \subset \mathbb{P}^1 \subset \dots \subset \mathbb{P}^m$
(this is a special case of \cite[Theorem 3.5.3]{bgs}). Since the closure of
each strata is smooth, the simple perverse sheaf $\mc{L}_i$ is, up to a shift
in the derived category, the continuation by $0$ of the constant sheaf
on $\mathbb{P}^i.$ The ext algebra
$\oplusop{i,j} \mbox{Ext}(\mc{L}_i,\mc{L}_j)$
can be easily written down, in particular,
  \begin{equation*}
  \mbox{Ext}^{k}(\mc{L}_i,\mc{L}_j)
   = H^{k - |i-j|}(\mathbb{P}^i \cap \mathbb{P}^j, \C).
  \end{equation*}
This algebra is then seen to be Koszul dual to $A_m^{\C}$,
which implies Proposition~\ref{highest-weight}.

We thus matched the category of $A_m^{\C}$-modules with the category
$\Oo'_m.$ The functors $\U_i$ in the category $ A_m^{\C}\mbox{-mod}$ have a simple
interpretation in the framework of highest weight categories.  This
interpretation is easy to guess, since the $\U_i$'s are exact functors while
translation functors, and, more generally, projective functors
(see \cite{bg} and references therein) constitute the main examples of exact functors
in $\Oo.$

Precisely, under the equivalence of categories
$A_m^{\C}\mbox{-mod}\cong \Oo'_m,$
the functor $\U_i$ (formula (\ref{uu:equation})
becomes the functor of translation across the $i$-th wall.
$\U_i$ is the composition of tensoring with
$_iP$ over $A_m$ and then tensoring with $P_i$ over $\C.$ Translation across
the wall functor, on the other hand, is the composition of translation on
and off the wall functors. For our particular parabolic subcategory of
a regular block, translation on the wall functor takes it to the
category equivalent to the category of $\C$-vector spaces. Translation on
and off the $i$-th wall functors are two-sided adjoints, and the adjointness
morphisms give rise to natural transformations between $\U_i$ and the
identity functor. These transformations should coincide with natural
transformations $\beta_i$ and $\gamma_i$
(see formulas (\ref{beta}) and (\ref{gamma})).

\subsubsection{The nature of bimodule maps $\beta_i, \gamma_i$}

The bimodule homomorphism $\beta_i: P_i \oo \hsm _iP \to A_m$ is simply a
direct summand of the multiplication map $A_m \oo A_m \to A_m.$
The homomorphism $\gamma_i,$ which goes the other way, is only slightly more
mysterious and can be interpreted as follows. For simplicity of the discussion
change the base ring from $\Z$ to a field $\F$ and ignore the grading of
$A_m.$ One checks that for $i>0$ the projective module $P_i$ is also
injective, and for any $A_m$-module $M$ there is a nondegenerate pairing
  \begin{equation*}
   \mbox{Hom}(P_i,M)\times \mbox{Hom}(M,P_i) \to \mbox{Hom}(P_i,P_i) \to \F,
  \end{equation*}
which defines an isomorphism $\mbox{Hom}(P_i,M)
\cong \mbox{Hom}(M,P_i)^{\ast},$ functorial in $M.$
Armed with the latter, we can rewrite
the evaluation morphism
   \begin{equation*}
   M\oo \mbox{Hom}(M,P_i)\to P_i
   \end{equation*}
as a functorial map $M \to P_i \oo \mbox{Hom}(P_i,M),$ which, in turn,
determines a bimodule map $A_m \to P_i \oo \hsm _iP.$ This map is
precisely $\gamma_i.$

\subsubsection{Temperley-Lieb algebra at roots of unity}
\label{tl-algebra}

The Temperley-Lieb algebra $TL_n(q)$ has generators $U_1, \dots, U_n$
and defining relations
\begin{eqnarray*}
 U_i^2 & = & (q+q^{-1}) U_i, \\
 U_i U_{i\pm 1}U_i & = & U_i, \\
 U_i U_j & = & U_j U_i \hspace{0.2in} |i-j|>1.
\end{eqnarray*}
The Temperley-Lieb algebra is a quotient of the Hecke algebra of the
symmetric group. Just like
the Hecke algebra, the Temperley-Lieb algebra is semisimple
for a generic $q\in \C$ and acquires a radical at certain roots of unity.
Let $q= e^{\frac{i\pi}{r}}$ and $r> 2$. A theorem of Paul Martin (see \cite{martin},
Chapter 7, for the exact statement) says that $TL_n(q)$ is isomorphic to the direct sum
of matrix algebras and algebras Morita equivalent to $A_m.$ Alternative proofs of this
result can be found in Westbury \cite{westbury}, Goodman and Wenzl \cite{gw}.


\section{Geometric intersection numbers\label{sec:intersection-numbers}}

This part of the paper collects the necessary results about curves on surfaces.
The first section, in which geometric intersection numbers are defined, is a
mild generalization of \cite[Expos{\'e} III]{fathi-laudenbach-poenaru}. As an
easy application, we give a proof of the theorem of Birman and Hilden which was
mentioned in Section \ref{subsec:results}. After that we introduce the bigraded
version of geometric intersection numbers. Finally, we consider a way of
decomposing curves on a disc into certain standard pieces. The basic idea of
this is again taken from \cite[Expos{\'e} IV]{fathi-laudenbach-poenaru}.

%
%

\subsection{Geometric intersection numbers\label{subsec:curves}}

Let $S$ be a compact oriented surface, possibly with boundary, and $\Delta \subset S \setminus
\partial S$ a finite set of marked points. $\Diff(S,\partial S;\Delta)$ is the group of
diffeomorphisms $f: S \longrightarrow S$ with $f|\partial S = \id$ and $f(\Delta) = \Delta$.
By a {\em curve in $(S,\Delta)$} we mean a subset $c \subset S$ which is either a simple
closed curve in $S^o = S \setminus (\partial S \cup \Delta)$ and essential (not contractible
in $S^o$), or else the image of an embedding $\gamma: [0;1] \longrightarrow S$ which is
transverse to $\partial S$ and such that $\gamma^{-1}(\partial S \cup \Delta) = \{0;1\}$.
Thus, our curves are smooth and unoriented. Two curves $c_0,c_1$ are called isotopic if one
can be deformed into the other by an isotopy in $\Diff(S,\partial S;\Delta)$; we write $c_0
\isotopic c_1$ for this. Note that endpoints on $\partial S$ may not move during an isotopy.

Let $c_0, c_1$ be two curves in $(S,\Delta)$. We say that they have {\em minimal intersection}
if they intersect transversally, satisfy $c_0 \cap c_1 \cap \partial S = \emptyset$, and the
following condition is satisfied: take any two points $z_- \neq z_+$ in $c_0 \cap c_1$ which
do not both lie in $\Delta$, and two arcs $\alpha_0 \subset c_0$, $\alpha_1 \subset c_1$ with
endpoints $z_-,z_+$, such that $\alpha_0 \cap \alpha_1 = \{z_-,z_+\}$. Let $K$ be the connected
component of $S \setminus (c_0 \cup c_1)$ which is bounded by $\alpha_0 \cup \alpha_1$. Then
if $K$ is topologically an open disc, it must contain at least one point of $\Delta$.
\includefigure{minimal-intersection}{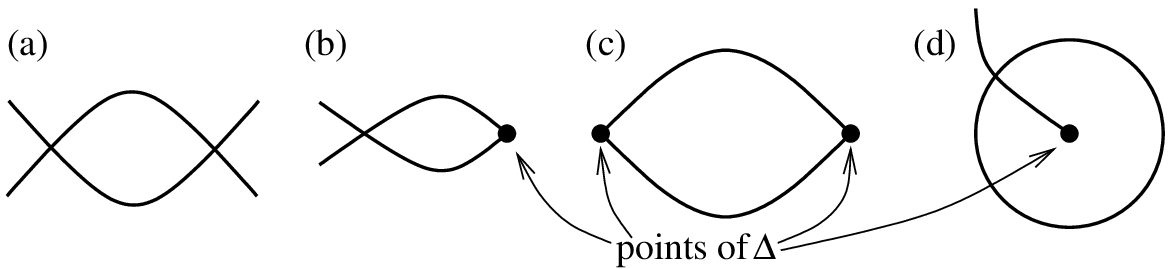}{hb}%

Among the examples in Figure\ \ref{fig:minimal-intersection}, (a) and (b) do not have minimal
intersection while (c) and (d) do.

\begin{lemma} \label{th:no-maps}
Assume that $\Delta = \emptyset$. Let $c_0,c_1$ be two curves in $(S,\emptyset)$ which
intersect transversally and satisfy $c_0 \cap c_1 \cap \partial S = \emptyset$. They have
minimal intersection iff the following property holds: there is no continuous map $v: [0;1]^2
\longrightarrow S$ with $v(s,0) \in c_0$, $v(s,1) \in c_1$, $v(0,t) = z_-$ and $v(1,t) = z_+$
for all $(s,t)$, where $z_- \neq z_+$ are points of $c_0 \cap c_1$.
\end{lemma}

The proof is a very slightly modified version of the equivalence $(2^o) \Leftrightarrow (3^o)$
of \cite[Proposition III.10]{fathi-laudenbach-poenaru}. Note that one of the implications is
obvious. There is a similar result for $\Delta \neq \emptyset$ but we will not need it.

Given two curves $c_0, c_1$ in $(S,\Delta)$ with $c_0 \cap c_1 \cap \partial S = \emptyset$,
one can always find a $c_1' \isotopic c_1$ which has minimal intersection with $c_0$, by
successively killing the unnecessary intersection points. We define the geometric intersection
number $I(c_0,c_1) \in \half\Z$ to be $I(c_0,c_1) = 2$ if $c_0,c_1$ are simple closed curves
with $c_0 \isotopic c_1$, and by the formula
\[
I(c_0,c_1) = |(c_0 \cap c_1') \setminus \Delta | + \half | c_0 \cap c_1' \cap \Delta|
\]
in all other cases. In words, one counts the intersection points of $c_0$ and $c_1'$, common
endpoints in $\Delta$ having weight $1/2$. The exceptional case in the definition can be
motivated by Floer cohomology, or by looking at the behaviour of $I$ on double branched covers
(see Section \ref{subsec:branched-cover}). The fact that $I(c_0,c_1)$ is independent of the
choice of $c_1'$ follows from the next two Lemmas. Once this has been established, it is
obvious that it is an invariant of the isotopy classes of $c_0$ and $c_1$.

\begin{lemma} \label{th:minimal-intersection}
Let $c_0$ be a curve in $(S,\Delta)$. Let $c_1',c_1''$ be two other curves such that both
$(c_0,c_1')$ and $(c_0,c_1'')$ have minimal intersection, and with $c_1' \isotopic c_1''$.
Assume moreover that $c_0 \not\isotopic c_1',c_1''$. Then there is an isotopy rel $c_0$ which
carries $c_1'$ to $c_1''$ (more formally, there is a smooth path $(f_t)_{0 \leq t \leq 1}$ in
$\Diff(S,\partial S;\Delta)$ such that $f_t(c_0) = c_0$ for all $t$, $f_0 = \id$, and
$f_1(c_1') = c_1''$).
\end{lemma}

For closed curves this is \cite[Proposition III.12]{fathi-laudenbach-poenaru}; the proof of
the general case is similar. The assumption that $c_0 \not\isotopic c_1', c_1''$ cannot be
removed, but there is a slightly weaker statement which fills the gap:

\begin{lemma} \label{th:isotopic-curves}
Let $c_0,c_1$ be two curves in $(S,\Delta)$ which are isotopic and have minimal intersection
(by definition, this implies that the $c_i$ are either closed, or else that their endpoints
all lie in $\Delta$). Then there is an isotopy rel $c_0$ which carries $c_1$ to one of the
curves $c_1',c_1''$ shown (in the two cases) in Figure\ \ref{fig:isotopic}.
\end{lemma}
\includefigure{isotopic}{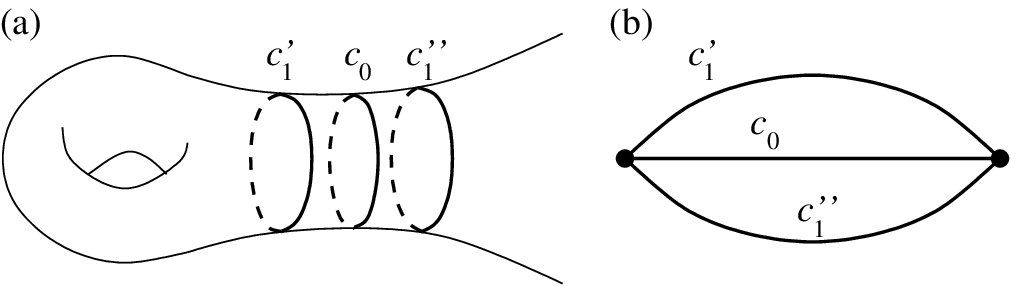}{ht}%

The geometric intersection number can be extended to curves $c_0, c_1$ which intersect on
$\partial S$, as follows: take a nonvanishing vector field on $\partial S$ which is positively
oriented, and extend it to a smooth vector field $Z$ on $S$ which vanishes on $\Delta$. Let
$(f_t)$ be the flow of $Z$. Set $c_0^+ = f_t(c_0)$ for some small $t>0$, and define
$I(c_0,c_1) := I(c_0^+,c_1)$. This depends on the orientation of $S$, and is also no
longer symmetric.

%
%

\subsection{Curves on a disc\label{subsec:disc}}

Let $D$ be a closed disc, with some fixed orientation, and $\Delta \subset D \setminus
\partial D$ a set of $m+1 \geq 2$ marked points. $\Conf_{m+1}(D \setminus \partial D)$ denotes
the configuration space of unordered $(m+1)$-tuples of points
in $D\setminus \partial D.$ As in Section
\ref{subsec:results} we write $\Diffeo = \Diff(D,\partial D;\Delta)$. There is a well-known
canonical isomorphism
\begin{equation} \label{eq:canonical-isomorphism}
\pi_1(\Conf_{m+1}(D \setminus \partial D),\Delta) \iso \pi_0(\Diffeo)
\end{equation}
defined as follows: take the evaluation map $\Diff(D,\partial D) \longrightarrow \Conf_{m+1}(D
\setminus \partial D)$, $f \longmapsto f(\Delta)$. This is a Serre fibration, with the fibre
at $\Delta$ being exactly $\Diffeo$, and \eqref{eq:canonical-isomorphism} is the boundary map
in the resulting sequence of homotopy groups.

From now on, by a curve we mean a curve in $(D,\Delta)$, unless otherwise specified. A {\em
basic set of curves} is a collection of curves $b_0,\dots,b_m$ which looks like that in Figure\
\ref{fig:basic-curves} (more rigorously, this means that the $b_i$ can be obtained by applying
some element of $\Diffeo$ to the collection drawn in Figure\ \ref{fig:basic-curves}).

\begin{lemma} \label{th:curves-and-diffeos}
Let $b_0,\dots,b_m$ be a basic set of curves. If $f \in \Diffeo$ satisfies $f(b_i) \isotopic
b_i$ for all $i$, then $[f] \in \pi_0(\Diffeo)$ is the identity class.
\end{lemma}

The proof is not difficult (it rests on the fact that one can isotop such an $f$ to a $g$
which satisfies $g(b_i) = b_i$ for all $i$). A basic set of curves determines a preferred
isomorphism between $B_{m+1}$, which we think of as an abstract group defined by the standard
presentation, and $\pi_0(\Diffeo)$ resp.\ $\pi_1(\Conf_{m+1}(D \setminus \partial D),\Delta)$.
Concretely, the $i$-th generator $\sigma_i$ ($1 \leq i \leq m$) of $B_{m+1}$
goes to the isotopy class
$[t_i] \in \pi_0(\Diffeo)$ of the half-twist along $b_i$. This half-twist is a diffeomorphism
which is trivial outside a small neighbourhood of $b_i$, which reverses $b_i$ itself, and such
that the image of $b_{i-1}$ is as shown in Figure\ \ref{fig:half-twist}(a). The corresponding
element in the fundamental group of the configuration space is the path which rotates the two
endpoints of $b_i$ around each other by $180^{\mathrm{o}}$, see Figure\ \ref{fig:half-twist}(b).
\includefigure{half-twist}{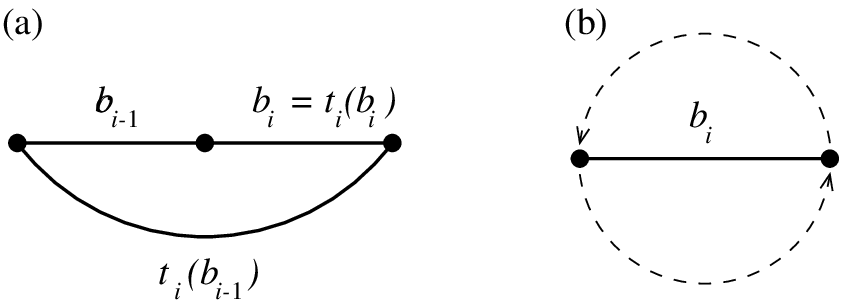}{ht}%

Fix a basic set of curves $b_0,\dots,b_m$. A curve $c$ will be called admissible if it is
equal to $f(b_i)$ for some $f \in \Diffeo$ and $0 \leq i \leq m$. The endpoints of an
admissible curve must lie in $\Delta \cup (b_0 \cap \partial D)$; conversely any curve with
such endpoints is admissible. There is an obvious action of $\pi_0(\Diffeo)$ on the set of
isotopy classes of admissible curves. Lemma \ref{th:curves-and-diffeos} is a faithfulness
result for this action: it shows that only the identity element acts trivially.

\begin{lemma} \label{th:like-ck}
Let $c$ be an admissible curve. Assume that there is a $k \in \{0,\dots,m\}$ such that
$I(b_i,c) = I(b_i,b_k)$ for all $i = 0,\dots,m$. Then
\[
c \isotopic
\begin{cases}
 b_0 \text{ or } \tau_0^{-1}(b_0) & \text{ if } k = 0,\\
 b_k \text{ or } \tau_k^{\pm 1}(b_k) & \text{ if } 1 \leq k < m, \\
 b_m & \text{ if } k = m.
\end{cases}
\]
Here $\tau_0,\dots,\tau_{m-1}$ are the positive Dehn twists along the closed curves
$l_0,\dots,l_{m-1}$ shown in Figure\ \ref{fig:loops}.
\end{lemma}
\includefigure{loops}{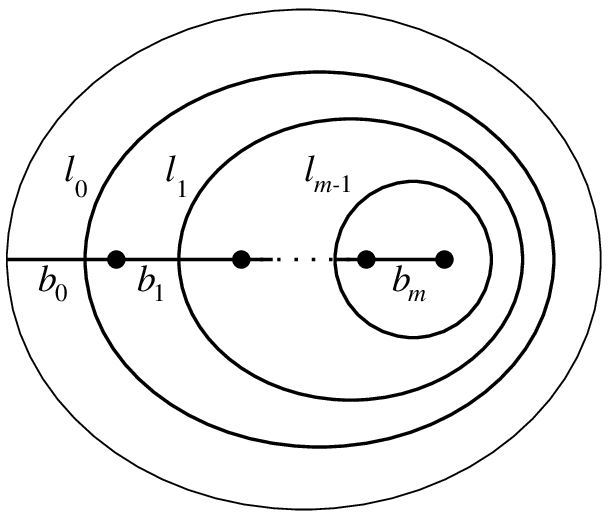}{ht}

\proof We will only explain the proof for $k = m$; the other cases are similar. We may assume
that $m \geq 2$, since the statement is fairly obvious for $m = 1$. Let $c$ be a curve
satisfying the conditions in the Lemma. Since $I(b_{m-1},c) = I(b_{m-1},b_m) = 1/2$ is not
integral, $c$ and $b_{m-1}$ have precisely one endpoint in common. This common endpoint could
either be $b_{m-2} \cap b_{m-1}$ or $b_{m-1} \cap b_m$, but the first case is impossible
because $I(b_{m-2},c) = I(b_{m-2},b_m) = 0$. We have now shown that $c$ and $b_m$ have one
endpoint in common. It follows that the other endpoints of these two curves must also be the
same, since $I(b_m,c) = I(b_m,b_m) = 1$ is integral. Next, we can assume that $c$ has minimal
intersection with all the $b_i$ (this is because the ``shortening moves'' which kill unnecessary
points of $b_i \cap c$ do not create new intersection points with any other $b_j$). Again by
using geometric intersection numbers, it follows that only the endpoints of $c$ lie on $b_0
\cup \dots \cup b_m$. This means that we can consider $c$ as a curve on the surface obtained
by cutting $D$ open along $b_0 \cup \dots \cup b_m$. The rest is straightforward. \qed

\begin{lemma} \label{th:detect-identity}
If $f \in \Diffeo$ satisfies $I(b_j,f^2(b_k)) = I(b_j,f(b_k)) = I(b_j,b_k)$ for all $j,k =
0,\dots m$, then $[f] \in \pi_0(\Diffeo)$ is the identity class.
\end{lemma}

\proof Applying Lemma \ref{th:like-ck} to $c = f(b_j)$ shows that there are numbers $\nu_0 \in
\{-1;0\}$ and $\nu_1,\dots,\nu_{n-1} \in \{-1;0;1\}$, such that
\[
f(b_j) \isotopic
 \begin{cases}
 \tau_j^{\nu_j}(b_j) & \text{if } j = 0,\dots,m-1, \\
 b_m & \text{if } j = m.
 \end{cases}
\]
Since the $\tau_j$ commute, and $\tau_j(b_k) = b_k$ for all $k \neq j$, the map $g =
\tau_0^{\nu_0} \tau_1^{\nu_1} \dots \tau_{m-1}^{\nu_{m-1}} \in \Diffeo$ satisfies $f(b_j)
\isotopic g(b_j)$ for all $j$. By Lemma \ref{th:curves-and-diffeos} it follows that $f$ and
$g$ lie in the same isotopy class. Applying the same argument to $f^2$ shows that it can be
written, up to isotopy, as a product $\tau_0^{\mu_0} \dots \tau_{m-1}^{\mu_{m-1}}$ with
numbers $\mu_j$ satisfying the same properties as $\nu_j$. In $\pi_0(\Diffeo)$ we have
therefore
\[
[\tau_0^{2\nu_0} \dots \tau_{m-1}^{2\nu_{m-1}}] = [f^2] = [\tau_0^{\mu_0} \dots
\tau_{m-1}^{\mu_{m-1}}].
\]
Since the $\tau_j$ generate a free abelian subgroup of rank $m$ in $\pi_0(\Diffeo)$, one has
$2\nu_j = \mu_j$ for all $j$. Now $|\mu_j| \leq 1$ which implies that $\nu_j = 0$, hence that
$[f] = [\id]$. \qed

%
%

\subsection{The double branched cover\label{subsec:branched-cover}}

In this section we assume that $D$ is embedded in $\C$.  The double cover of $D$ branched
along $\Delta$ is, by definition,
\[
p_S: S = \{ (x_0,x_1) \in \C \times D \suchthat x_0^2 + h(x_1) = 0 \} \longrightarrow D, \quad
p_S(x_0,x_1) = x_1,
\]
where $h \in \C[z]$ is a polynomial which has simple zeros exactly at the points of $\Delta$.

\begin{lemma} \label{th:preimage}
Let $c$ be a curve in $(D,\Delta)$ such that $c \cap \Delta \neq \emptyset$ (in future, we
will say that $c$ meets $\Delta$ if this is the case). Then $p_S^{-1}(c)$ is a curve in
$(S,\emptyset)$.
\end{lemma}

\proof Let $z$ be a point of $c \cap \Delta$, and $\gamma: [0;1) \longrightarrow c$ a local
parametrization of $c$ near that point. One can write $h(\gamma(t)) = -t\psi(t)$ for some
$\psi \in \smooth([0;1),\C^*)$. Choose a square root $\sqrt{\psi}$. Then a local smooth
parametrization of $p_S^{-1}(c)$ near $\{(0,z)\} = p_S^{-1}(z)$ is given by $(-1;1)
\longrightarrow S$, $t \longmapsto (t\sqrt{\psi(t^2)},\gamma(t^2))$. The smoothness of
$p_S^{-1}(c)$ at all other points is obvious. As for the topology, there are (because of the
assumption $c \cap \Delta \neq \emptyset$) only two possibilities: either $c$ is a path
joining two points of $\Delta$, and then $p_S^{-1}(c)$ is a simple closed curve which is not
homologous to zero; or $c$ connects a point of $\Delta$ with a point of $\partial D$, and then
$p_S^{-1}(c)$ is a path connecting two points of $\partial S$. \qed

\begin{lemma} \label{th:lifting-curves}
Let $c_0,c_1$ be curves on $(D,\Delta)$, both of which meet $\Delta$. Assume that they have
minimal intersection and are not isotopic. Then $p_S^{-1}(c_0),p_S^{-1}(c_1)$ also have
minimal intersection.
\end{lemma}

\proof It is clear that $p_S^{-1}(c_0)$ and $p_S^{-1}(c_1)$ intersect transversally, and do
not intersect on $\partial S$. Let $U \subset D$ be a thickening of $c_0 \cup c_1$, such that
$D \setminus U$ is a surface with corners. More precisely, the thickening (shown in Figure\
\ref{fig:corners}) should be done in such a way that a point of $(c_0 \cap c_1) \setminus
\Delta$ gives rise to four corners of $D \setminus U$, and a point of $c_0 \cap c_1 \cap
\Delta$ to two corners. Take some connected component $K$ of $D \setminus U$ which satisfies
$K \cap \partial D = \emptyset$. Then at least one of the following three assertions is true:
(1) $K$ is not contractible; (2) it has at least three corners; (3) it has two corners and
contains at least one point of $\Delta$. The point here is that $K$ cannot be a $2$-gon (a
disc with two corners) and disjoint from $\Delta$, since our assumptions exclude all the four
situations shown in Figure\ \ref{fig:minimal-intersection}: (a),(b) because of the minimal
intersection, (c) because $c_0 \not\isotopic c_1$, and (d) because $c_0,c_1$ are not closed
curves.
\includefigure{corners}{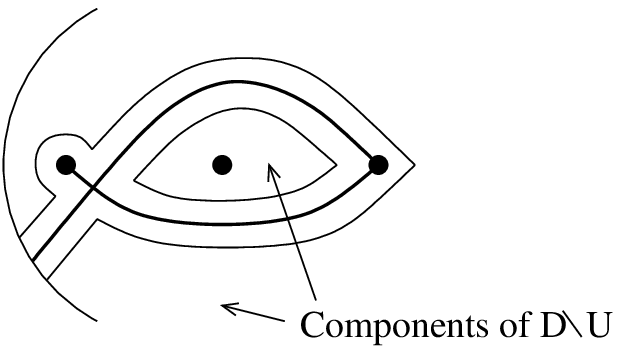}{hb}

If (1) holds then the preimage $p_S^{-1}(K)$ may consist of one or two connected components,
but none of them can be contractible. If (2) holds the connected components of $p_S^{-1}(K)$
may be contractible, but each of them has at least three corners. If (3) holds then
$p_S^{-1}(K)$ is connected with at least four corners, since its boundary is a double cover of
the boundary of $K$. We have now shown that no connected component of $S \setminus
p_S^{-1}(U)$ can be a $2$-gon and disjoint from $\partial S$. Essentially by definition, this
means that $p_S^{-1}(c_0)$ and $p_S^{-1}(c_1)$ have minimal intersection. \qed

\begin{lemma} \label{th:lifting-isotopies}
Let $c_0, c_1$ be two curves on $(D,\Delta)$, both of which meet $\Delta$, such that $c_0 \cap
c_1 \cap \partial D = \emptyset$. Assume that $c_0 \not\isotopic c_1$. Then $p_S^{-1}(c_0)
\not\isotopic p_S^{-1}(c_1)$.
\end{lemma}

\proof {\bf Case 1:} {\em $c_0$ or $c_1$ has an endpoint on $\partial D$.} Then $c_0 \cap
\partial D \neq c_1 \cap \partial D$, because we have assumed that the intersection of these
two sets is zero. It follows that $p_S^{-1}(c_0) \cap \partial S \neq p_S^{-1}(c_1) \cap
\partial S$, so these two curves cannot be isotopic. {\bf Case 2:} {\em $c_0,c_1$ are
contained in $D \setminus \partial D$}. This means that both of them are paths joining two
points of $\Delta$. Since the statement is one about isotopy classes, we may assume that
$c_0,c_1$ have minimal intersection. Then the same holds for $p_S^{-1}(c_0),p_S^{-1}(c_1)$ by
Lemma \ref{th:lifting-curves}. Assume that, contrary to what we want to show, $p_S^{-1}(c_0)
\isotopic p_S^{-1}(c_1)$. Then Lemma \ref{th:isotopic-curves} implies that $p_S^{-1}(c_0) \cap
p_S^{-1}(c_1) = \emptyset$, so that $c_0 \cap c_1 = \emptyset$, and in particular $c_0 \cap
\Delta \neq c_1 \cap \Delta$. But then $p_S^{-1}(c_0)$ and $p_S^{-1}(c_1)$ are not even
homologous, which is a contradiction. \qed

\begin{lemma} \label{th:comparing-gin}
Let $c_0,c_1$ be curves in $(D,\Delta)$, both of which meet $\Delta$. Then
\[
I(p_S^{-1}(c_0),p_S^{-1}(c_1)) = 2\,I(c_0,c_1).
\]
\end{lemma}

\proof It is enough to prove the statement for $c_0 \cap c_1 \cap \partial D = \emptyset$,
since the other situation reduces to this by definition. {\bf Case 1:} $c_0 \isotopic c_1$.
Because we are assuming that $c_0 \cap c_1 \cap \partial D = \emptyset$, this means that
$c_0,c_1$ are paths connecting two points of $\Delta$; then $I(c_0,c_1) = 1$. On the other
hand $p_S^{-1}(c_0)$ and $p_S^{-1}(c_1)$ are simple closed curves and isotopic, so that
$I(p_S^{-1}(c_0),p_S^{-1}(c_1)) = 2$. {\bf Case 2:} $c_0 \not\isotopic c_1$. Then
$p_S^{-1}(c_0) \not\isotopic p_S^{-1}(c_1)$ by Lemma \ref{th:lifting-isotopies}. Choose a
$c_1' \isotopic c_1$ which has minimal intersection with $c_0$. By Lemma
\ref{th:lifting-curves} one has $I(p_S^{-1}(c_0),p_S^{-1}(c_1)) = |p_S^{-1}(c_0) \cap
p_S^{-1}(c_1')| = 2| (c_0 \cap c_1') \setminus \Delta| + |c_0 \cap c_1' \cap \Delta| =
2I(c_0,c_1)$. \qed

Any $f \in \Diffeo$ can be lifted in a unique way to a homeomorphism $f_S$ of $S$ with
$f_S|\partial S = \id$. If $f$ is holomorphic in some neighbourhood of each point of $\Delta$,
then $f_S$ is again smooth. This defines a lifting homomorphism $\pi_0(\Diffeo)
\longrightarrow \pi_0(\Diff(S,\partial S))$.

\begin{prop}[Birman and Hilden] \label{th:birman-hilden}
The lifting homomorphism is injective.
\end{prop}

\proof Choose a set of basic curves $b_0,\dots,b_{m-1}$ in $(D,\Delta)$. Consider a map $f \in
\Diffeo$ which is holomorphic near each point of $\Delta$, and assume that $[f]$ lies in the
kernel of the lifting homomorphism. Then $p_S^{-1}(f(c)) = f_S(p_S^{-1}(c)) \isotopic
p_S^{-1}(c)$ for any curve $c$ which meets $\Delta$. Therefore, by Lemma
\ref{th:comparing-gin},
\[
I(f(b_j),b_k) = \half I(p_S^{-1}(f(b_j)),p_S^{-1}(b_k)) = \half I(p_S^{-1}(b_j),p_S^{-1}(b_k))
= I(b_j,b_k)
\]
for $j,k = 0,\dots,m$. The same clearly holds for $f^2$. Lemma \ref{th:detect-identity} shows
that $[f]$ must be the identity class. \qed

We have chosen this proof, which is not the most direct one, because it serves as a model for
later arguments. One can in fact avoid geometric intersection numbers entirely, by
generalizing Lemma \ref{th:lifting-isotopies} to the case when $c_0 \cap c_1 \cap \partial D$
is not necessarily empty, and then using Lemma \ref{th:curves-and-diffeos}. Birman and
Hilden's original paper \cite{birman-hilden73} contains a stronger result (they identify the
image of the lifting homomorphism) and their approach is entirely different.

%
%

\subsection{Bigraded intersection numbers\label{subsec:bigraded-curves}}

Let $P = PT(D \setminus \Delta)$ be the (real) projectivization of the tangent bundle of $D
\setminus \Delta$. We want to introduce a particular covering of $P$ with covering group
$\Z^2$. Take an oriented trivialization of $TD$; this allows one to identify $P = \RP{1} \times (D
\setminus \Delta)$. For every point $z \in \Delta$ choose a small loop $\lambda_z: S^1
\longrightarrow D \setminus \Delta$ winding positively once around it. The classes $[point
\times \lambda_z]$ together with $[\RP{1} \times point]$ form a basis of $H_1(P;\Z)$. Let $C
\in H^1(P;\Z^2)$ be the cohomology class which satisfies $C([point \times \lambda_z]) =
(-2,1)$ and $C([\RP{1} \times point]) = (1,0)$. Define $\tP$ to be the covering classified by
$C$, and denote the $\Z^2$-action on it by $\Pchi$. This is independent of the choices we
have made, up to isomorphism.

Any $f \in \Diffeo$ induces a diffeomorphism of $P$ which, since it preserves $C$, can be
lifted to an equivariant diffeomorphism of $\tP$. There is a unique such lift, denoted by
$\tf$, which acts trivially on the fibre of $\tP$ over any point $T_z\partial D \in P$, $z \in
\partial D$. For any curve $c$ there is a canonical section $s_c: c \setminus \Delta
\longrightarrow P$ given by $s_c(z) = T_zc$. A {\em bigrading} of $c$ is a lift $\tilde{c}$ of
$s_c$ to $\tP$. Pairs $(c,\tc)$ consisting of a curve and a bigrading are called {\em bigraded
curves}; we will often write $\tc$ instead of $(c,\tc)$. $\Pchi$ defines a $\Z^2$-action
on the set of bigraded curves, and the lifts $\tilde{f}$ yield a $\Diffeo$-action on the same
set. There is an obvious notion of isotopy for bigraded curves, and one obtains induced
actions of $\Z^2$ and $\pi_0(\Diffeo)$ on the set of isotopy classes.

An alternative formulation is as follows. Assume that $D$ is embedded in $\C$,
and let $h \in \C[z]$ be as in the previous section.
The embedding $D\subset \C$ gives a preferred trivialization of $TD,$ which we use to
identify $P$ with $\C^*/\R \times (D \setminus \Delta).$ The map
\[
\delta_P: P \longrightarrow (\C^*/\Rg) \times (\C^*/\Rg), \quad (\zeta,z) \longmapsto
(h(z)^{-2}\zeta^2,-h(z))
\]
represents $C$. Hence one can take $\tP$ to be the pullback of the universal covering $\R^2
\longrightarrow (\C^*/\Rg)^2$, $(\xi_1,\xi_2) \mapsto (\exp(2\pi i \xi_1), \exp(2\pi i
\xi_2))$. Then a bigrading of a curve $c$ is just a continuous lift $\tc: c \setminus \Delta
\longrightarrow \R^2$ of the map $c \setminus \Delta \longrightarrow (\C^*/\Rg)^2$, $z \mapsto
\delta_P(T_zc)$; and the $\Z^2$-action is by adding constant functions to $\tc$. To put it
even more concretely, take an embedding $\gamma: (0;1) \longrightarrow D$ which parametrizes
an open subset of $c \setminus \Delta$. Then the map which one has to lift can be written as
\begin{equation} \label{eq:explicit-deltap}
\delta_P(T_{\gamma(t)}c) = (h(\gamma(t))^{-2}\gamma'(t)^2, -h(\gamma(t))).
\end{equation}
\begin{lemma}
A curve $c$ admits a bigrading iff it is not a simple closed curve.
\end{lemma}

\proof If $c$ is not a simple closed curve, $c \setminus \Delta$ is contractible, so that the
pullback $s_c^*\tP \longrightarrow c \setminus \Delta$ is the trivial covering. Conversely, if
$c$ is a simple closed curve and bounds a region of $D$ containing $k > 0$ points of $\Delta$,
then $s_c^*C([c]) = \pm (2-2k,k) \neq 0$, which means that $s_c^*\tP$ does not admit a
section. \qed

\begin{lemma} \label{th:self-isotopy}
The $\Z^2$-action on the set of isotopy classes of bigraded curves is free. Equivalently, a
bigraded curve $\tc$ is never isotopic to $\chi(r_1,r_2)\tc$ for any
$(r_1,r_2) \neq 0$.
\end{lemma}

\proof Let $c$ be a curve which connects two points $z_0,z_1 \in \Delta$. Let $(f_t)_{0 \leq t
\leq 1}$ be an isotopy in $\Diffeo$ with $f_0 = \id$ and $f_1(c) = c$. Without any real loss
of generality we may assume that $f_1|c = \id$. {\bf Claim:} {\em the closed path $\kappa_z:
[0;1] \longrightarrow D \setminus \Delta$, $\kappa_z(t) = f_t(z)$ is freely nullhomotopic for
each $z \in c \setminus \Delta$}. The free homotopy class of $\kappa_z$ is evidently
independent of $z$. On the other hand, for $z$ close to $z_0$ it must be a multiple of the
free homotopy class of $\lambda_{z_0}$, and the corresponding statement holds for $z$ close to
$z_1$; which proves our claim.

As a consequence one can deform the isotopy $(f_t)_{0 \leq t \leq 1}$, rel endpoints, into
another isotopy $(g_t)$ such that $g_t(c)$ agrees with $c$ in a neighbourhood of $z_0,z_1$ for
all $t$. By considering the preferred lifts $\tilde{g}_t$ one sees immediately that
$\tilde{f}_1(\tc) = \tilde{g}_1(\tc) = \tc$ for any bigrading $\tc$. The proof for the
remaining types of curves, which satisfy $c \cap \partial \Delta \neq \emptyset$, is much
easier, and we leave it to the reader. \qed

\begin{lemma} \label{th:self-shift}
Let $c$ be a curve which joins two points of $\Delta$, $t_c \in \Diffeo$ the half twist along
it, and $\tilde{t}_c$ its preferred lift to $\tP$. Then $\tilde{t}_c(\tc) =
\Pchi(-1,1)\tc$ for any bigrading $\tc$ of $c$.
\end{lemma}

\proof Since $t_c(c) = c$, one has $\tilde{t}_c(\tilde{c}) = \Pchi(r_1,r_2)\tc$ for some
$(r_1,r_2) \in \Z^2$. Take an embedded smooth path $\beta: [0;1] \longrightarrow D \setminus
\Delta$ from a point $\beta(0) \in \partial D$ to the unique point $\beta(1) \in c$ which is a
fixed point of $t_c$; Figure\ \ref{fig:self-shift} illustrates the situation. Consider the
closed path $\pi:[0;2] \longrightarrow P$ given by $\pi(t) = \R\beta'(t) \subset
T_{\beta(t)}D$ for $t \leq 1$, and by $\pi(t) = Dt_c(\R\beta'(2-t))$ for $t \geq 1$. It is not
difficult to see that (with respect to the basis of $H_1(P;\Z)$ which we have used before)
$[\pi] = -[\RP{1} \times point] - [point \times \lambda_z]$, where $z$ is one of the endpoints
of $c$. Therefore $(r_1,r_2) = -C([\pi]) = C([\RP{1} \times point]) + C([point \times
\lambda_z]) = (-1,1)$. \qed
\includefigure{self-shift}{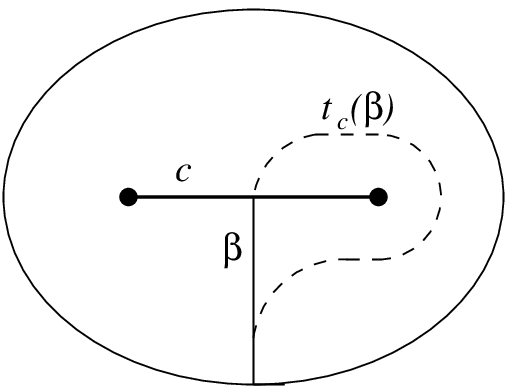}{ht}%

Let $(c_0,\tc_0)$ and $(c_1,\tc_1)$ be two bigraded curves, and $z \in D \setminus \partial D$
a point where $c_0$ and $c_1$ intersect transversally; $z$ may lie in $\Delta$ or not. Fix a
small circle $l \subset D \setminus \Delta$ around $z$. Let $\alpha: [0;1] \longrightarrow l$
be an embedded arc which moves clockwise around $l$, such that $\alpha^{-1}(c_0) = \{0\}$ and
$\alpha^{-1}(c_1) = \{1\}$. If $z \in \Delta$ then $\alpha$ is unique up to a change of
parametrization; otherwise there are two possibilities, which are distinguished by their
endpoints (see Figure\ \ref{fig:local-index}). Take a smooth path $\pi: [0;1] \longrightarrow P$
with $\pi(t) \in P_{\alpha(t)}$ for all $t$, going from $\pi(0) = T_{\alpha(0)}c_0$ to $\pi(1)
= T_{\alpha(1)}c_1$, such that $\pi(t) \neq T_{\alpha(t)}l$ for all $t$. One can picture $\pi$
as a family of tangent lines along $\alpha$ which are all transverse to $l$ (see again Figure\
\ref{fig:local-index}). Lift $\pi$ to a path $\tilde{\pi}: [0;1] \longrightarrow \tP$ with
$\tilde{\pi}(0) = \tc_0(\alpha(0))$; necessarily $\tc_1(\alpha(1)) =
\Pchi(\mu_1,\mu_2)\tilde{\pi}(1)$ for some $\mu_1,\mu_2$. We define the local index of
$\tc_0,\tc_1$ at $z$ to be $\mubigr(\tc_0,\tc_1;z) = (\mu_1,\mu_2) \in \Z^2$; it is not
difficult to see that this is independent of all choices.
\includefigure{local-index}{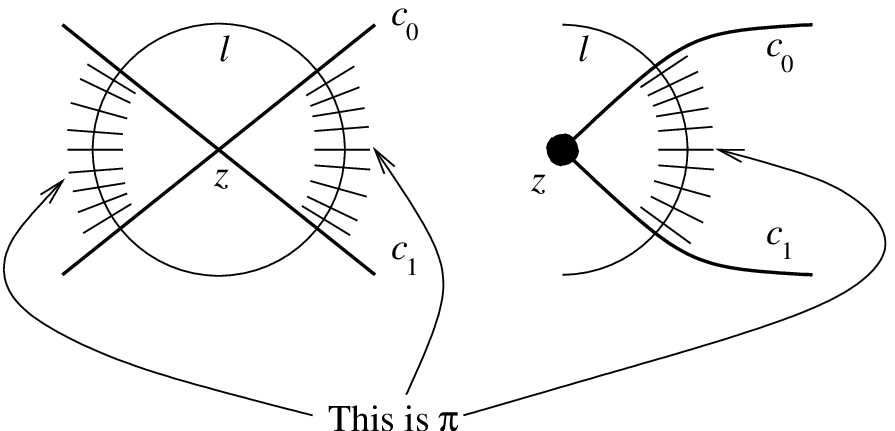}{ht}%

Now let $\tc_0$, $\tc_1$ be two bigraded curves such that $c_0 \cap c_1 \cap \partial D =
\emptyset$. Choose a $c_1' \isotopic c_1$ which has minimal intersection with $c_0$. There is
a bigrading $\tc_1'$ of $c_1'$, which is unique by Lemma \ref{th:self-isotopy}, such that
$\tc_1' \isotopic \tc_1$. The bigraded intersection number $\Ibigr(\tc_0,\tc_1) \in
\Z[q_1,q_1^{-1},q_2,q_2^{-1}]$ is defined by
\[
\Ibigr(\tc_0,\tc_1) = (1 + q_1^{-1}q_2)\Big(\sum_{z \in (c_0 \cap c_1') \setminus \Delta}
q_1^{\mu_1(z)}q_2^{\mu_2(z)}\Big) + \Big(\sum_{z \in c_0 \cap c_1' \cap \Delta}
q_1^{\mu_1(z)}q_2^{\mu_2(z)}\Big),
\]
where $(\mu_1(z),\mu_2(z)) = \mubigr(\tc_0,\tc_1';z)$. The proof that this is independent of
the choice of $c_1'$, and is an invariant of the isotopy classes of $\tc_0, \tc_1$, is
basically the same as for ordinary geometric intersection numbers. The only non-obvious case
is when $c_0$ and $c_1$ are isotopic, because then there are two essentially different choices
for $c_1'$. However, an explicit computation shows that both choices lead to the same result.
Bigraded intersection numbers can be extended to curves with $c_0 \cap c_1 \cap \partial D
\neq \emptyset$ just like ordinary ones, by taking a flow $(f_t)$ which moves $\partial D$ in
the positive sense, lifting the induced flow on $P$ to one $(\tilde{f}_t)$ on $\tP$
such that $\tilde{f}_0=\id,$
and setting $\Ibigr(\tc_0,\tc_1) = \Ibigr(\tilde{f}_t(\tc_0),\tc_1)$ for
small $t>0$. We list some elementary properties of $\Ibigr$:
\begin{Blist}
\item
$I(c_0,c_1)$ is obtained from $\Ibigr(\tc_0,\tc_1)$ by setting $q_1 = q_2 = 1$ and dividing by
two;
\item \label{item:bigraded-map}
$\Ibigr(\tf(\tc_0),\tf(\tc_1)) = \Ibigr(\tc_0,\tc_1)$ for any $f \in \Diffeo$;
\item \label{item:bigraded-shift}
$\Ibigr(\tc_0,\Pchi(r_1,r_2)\tc_1) = \Ibigr(\Pchi(-r_1,-r_2)\tc_0,\tc_1) =
q_1^{r_1}q_2^{r_2}\Ibigr(\tc_0,\tc_1)$.
\item \label{item:symmetry}
If $c_0 \cap c_1 \cap \partial D = \emptyset$ and $\Ibigr(\tc_0,\tc_1) = \sum_{r_1,r_2}
a_{r_1,r_2} q_1^{r_1}q_2^{r_2}$ then $\Ibigr(\tc_1,\tc_0) = \sum_{r_1,r_2} a_{r_1,r_2}
q_1^{-r_1}q_2^{1-r_2}$.
\end{Blist}

The only item which requires some explanation is the last one. Assume that $c_0,c_1$ have
minimal intersection. The local index has the following symmetry property:
\[
\mubigr(\tc_1,\tc_0;z) =
 \begin{cases}
 (1,0) - \mubigr(\tc_0,\tc_1;z) & \text{if } z \notin \Delta,\\
 (0,1) - \mubigr(\tc_0,\tc_1;z) & \text{if } z \in \Delta.
 \end{cases}
\]
One sees immediately that \ref{item:symmetry} holds for the contributions coming from points
of $c_0 \cap c_1 \cap \Delta$. Now take a point $z \in (c_0 \cap c_1) \setminus \Delta$, and
assume for simplicity that $\mubigr(\tc_0,\tc_1;z) = (0,0)$. Then the contribution of $z$ to
$\Ibigr(\tc_0,\tc_1)$ is $1 + q_1^{-1}q_2$. Since $\mubigr(\tc_1,\tc_0;z) = (1,0)$, the same
point contributes $q_1(1 + q_1^{-1}q_2) = q_1 + q_2$ to $\Ibigr(\tc_1,\tc_0)$, which is in
agreement with \ref{item:symmetry}. We point out that the term $(1+q_1^{-1}q_2)$ in the
definition of bigraded intersection numbers was essential in this computation.

%
%

\subsection{Normal form}
\label{normalform}

Fix a set of basic curves $b_0, \dots, b_m$. In addition choose curves $d_0,\dots,d_m$ as in
Figure\ \ref{fig:bc} which divide $D$ into regions $D_0, \dots, D_{m+1}$. We say that an
admissible curve $c$ is in normal form if it has minimal intersection with all the $d_k$. This
can always be an achieved by an isotopy. The corresponding uniqueness theorem is a slight
generalization of Lemma \ref{th:minimal-intersection}:
\includefigure{bc}{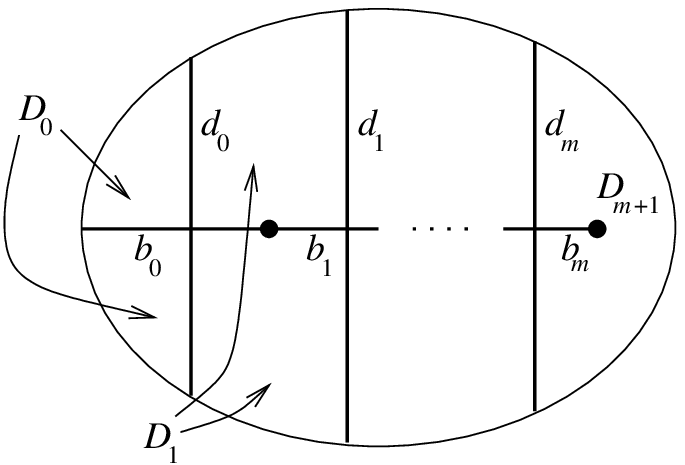}{htb}%

\begin{lemma} \label{th:normal-form}
Let $c_0,c_1$ be two isotopic admissible curves, both of which are in normal form. Then there
is an isotopy rel $d_0 \cup d_1 \cup \dots \cup d_m$ which carries $c_0$ to $c_1$.
\end{lemma}

\begin{figure}
\begin{center}
\epsfig{file=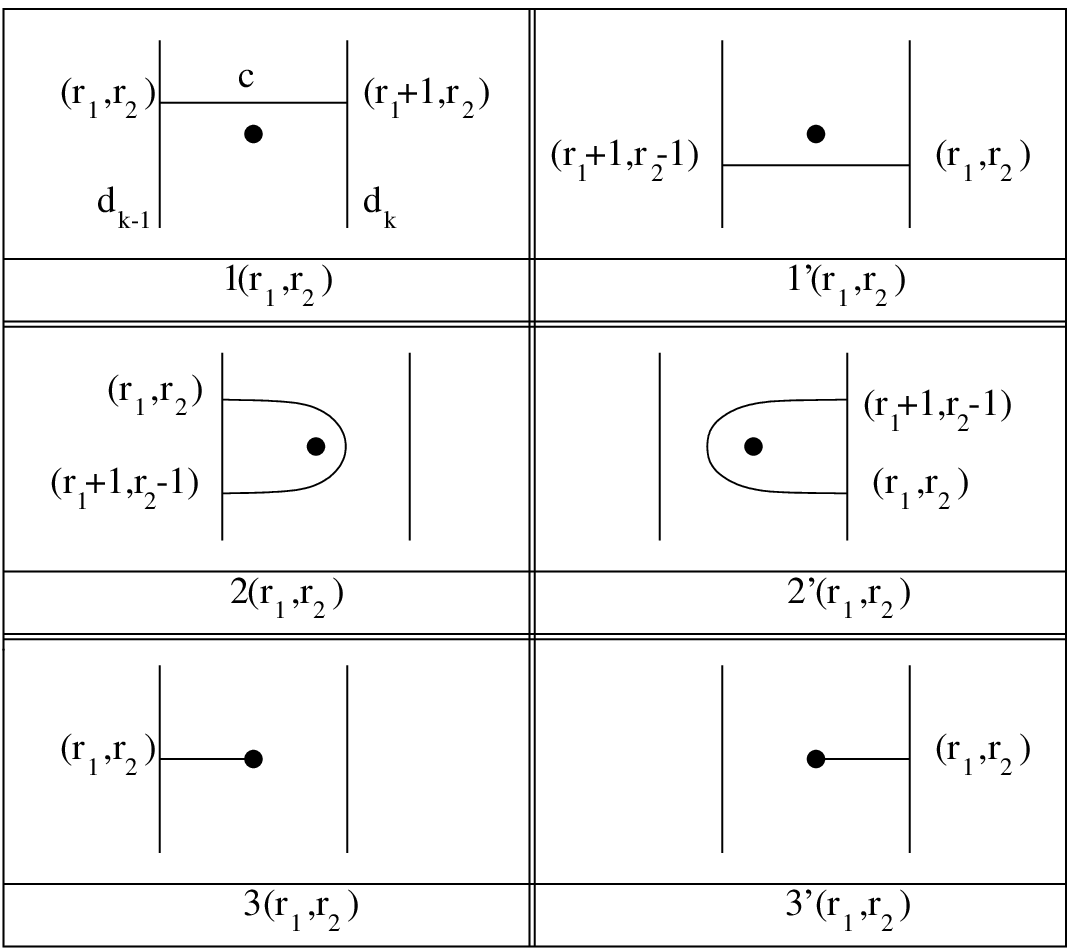} \vspace{-1em}
\caption{Connected components of $d \cap D_k$, for $d$ in
normal form and $1 \leq k \leq m$} \label{fig:typesa} \vspace{2em}

\epsfig{file=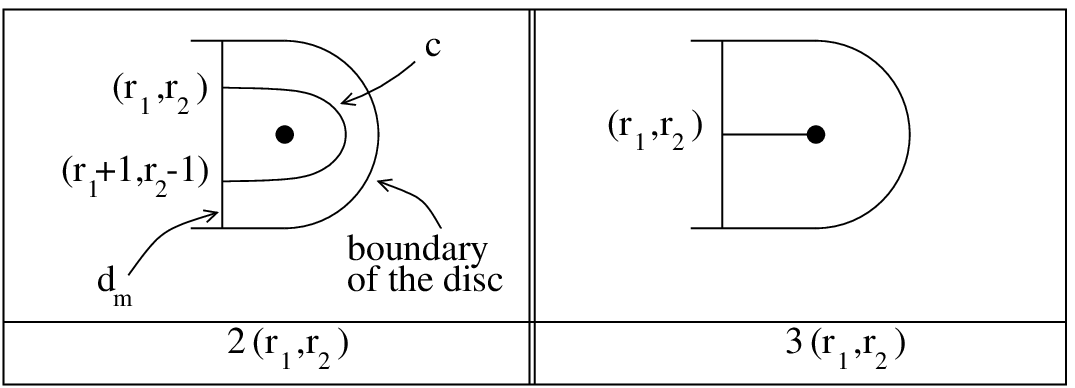} \vspace{-1em}
\caption{The same for $k = m+1$} \label{fig:typesb} \vspace{2em}

\epsfig{file=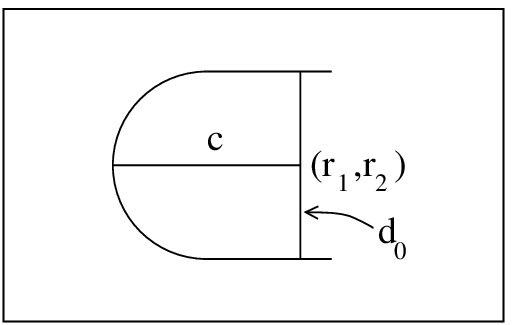} \vspace{-1em}
\caption{The same for $k = 0$} \label{fig:typesc}
\end{center}
\end{figure}

Let $c$ be a curve in normal form. Each connected component of $c \cap D_k$ belongs to one of
finitely many classes or types. For $1 \leq k \leq m$ there are six such types, denoted by
$1$, $1'$, $2$, $2'$, $3$, $3'$ (see Figure\ \ref{fig:typesa}; the numbers $(r_1,r_2)$ here and
in the subsequent figures are for later use, and the reader should ignore them at present).
For $k = m+1$ there are two types, which are analogues of the types $2$ and $3$ considered
before; we use the same notation for them (Figure\ \ref{fig:typesb}). For $k = 0$ there is a
single type, for which we will not need a proper name (Figure\ \ref{fig:typesc}). Conversely, it
is not difficult to prove that an admissible curve $c$ which intersects all the $d_k$
transversally, and such that each connected component of $c \cap D_k$ belongs to one of the
types listed in Figures \ \ref{fig:typesa}--\ref{fig:typesc}, is already in normal form.

\begin{remark}
For a completely accurate formulation, one should consider the group $\Diffeo_k$ of
diffeomorphisms $f: D_k \longrightarrow D_k$ which satisfy $f|\partial D \cap D_k = \id$,
$f(d_{k-1}) = d_{k-1}$, $f(d_k) = d_k$, and $f(\Delta \cap D_k) = \Delta \cap D_k$. Our list
of types classifies each connected component of $c \cap D_k$ up to an isotopy in $\Diffeo_k$.
\end{remark}

Lemma \ref{th:normal-form} implies that the number of connected components of each type, as well
as their relative position and the way in which they join each other,
is an invariant of the isotopy class of $c$.

For the rest of this section, $c$ is an admissible curve in normal form.
The points of $\cross(c) = c\cap (d_0\cup d_1\cup \dots \cup d_m)$
are called {\em crossings}, and those in the subset $c\cap d_k$ are called
{\em $k$-crossings} of $c.$ The connected components of $c\cap D_k,$
$0\leq k \leq m,$ are called \emph{segments} of $c.$
A segment is \emph{essential} if its endpoints are both crossings (as opposed
to points of $\Delta \cup \partial D$). Thus, the essential segments are precisely
those of type $1,1',2$ or $2'.$ Note that the basic curves $b_0, \dots, b_m$ have no
essential segments.

$c$ can be reconstructed up to isotopy by listing its crossings, plus the types of
essential segments bounded by consecutive crossings aw one travels along $c$ from
one end to the other. Conversely, Lemma~\ref{th:normal-form} shows that this
combinatorial data is an invariant of the isotopy class of $c.$

\clearpage

We now discuss how the normal form changes under the half-twist $t_k$ along $b_k.$
The curve $t_k(c),$ in general, is not in normal form.
$t_k(c),$ however, has minimal intersection with each $d_i$ for $i\not= k.$
To bring $t_k(c)$ into normal form we only need to
simplify its intersections with $d_k.$ More precisely, there is a bijection between
connected components of  the intersection of $c$ and $t_k(c)$ with $D_k \cup D_{k+1}.$
Take a connected component $y$ of $c\cap (D_k \cup D_{k+1})$ and deform (inside
$D_k \cup D_{k+1}$) the twisted component $t_k(y),$ keeping its intersection with
$d_{k-1}\cup d_{k+1}\cup \Delta$ fixed,
so that it has minimal intersection with $d_k.$
After simplifying each connected component of  $c\cap (D_k \cup D_{k+1})$
in this way we end up with normal form for $t_k(c).$
We collect our observations into the following

\begin{proposition} \label{bijections}
(a) The normal form of $t_k(c)$ coincides with $c$
outside of $D_k\cup D_{k+1}.$ The curve $t_k(c)$ can be brought into normal form by
an isotopy inside $D_k\cup D_{k+1}.$

(b) Now assume that $t_k(c)$ is in normal form.
There is a natural bijection between $i$-nodes of $c$ and $t_k(c)$
for each $i\not= k.$ There is a natural bijection between connected components of
intersections of $c$ and $t_k(c)$ with $D_k \cup D_{k+1}.$
\end{proposition}

Define a $k$\emph{-string of} $c$ as a connected component of $c\cap (D_k\cup D_{k+1}).$
Denote by $\mathrm{st}(c,k)$ the set of $k$-strings of $c.$

Define a $k$\emph{-string} as a curve in $D_k\cup D_{k+1}$ which is a connected component of
$c\cap (D_k\cup D_{k+1})$ for some admissible curve $c$ (recall that we assume $c$ to be in
normal form).

We say that two $k$-strings are isotopic (or belong to the same isotopy class)
if there is a deformation of one into the other
via diffeomorphisms $f$ of $D'=D_k\cup D_{k+1}$ which satisfy $f(d_{k-1})= d_{k-1},
f(d_{k+1})= d_{k+1}$ and $f(\Delta \cap D')= \Delta \cap D'.$

\begin{figure}
\begin{center}
\epsfig{file=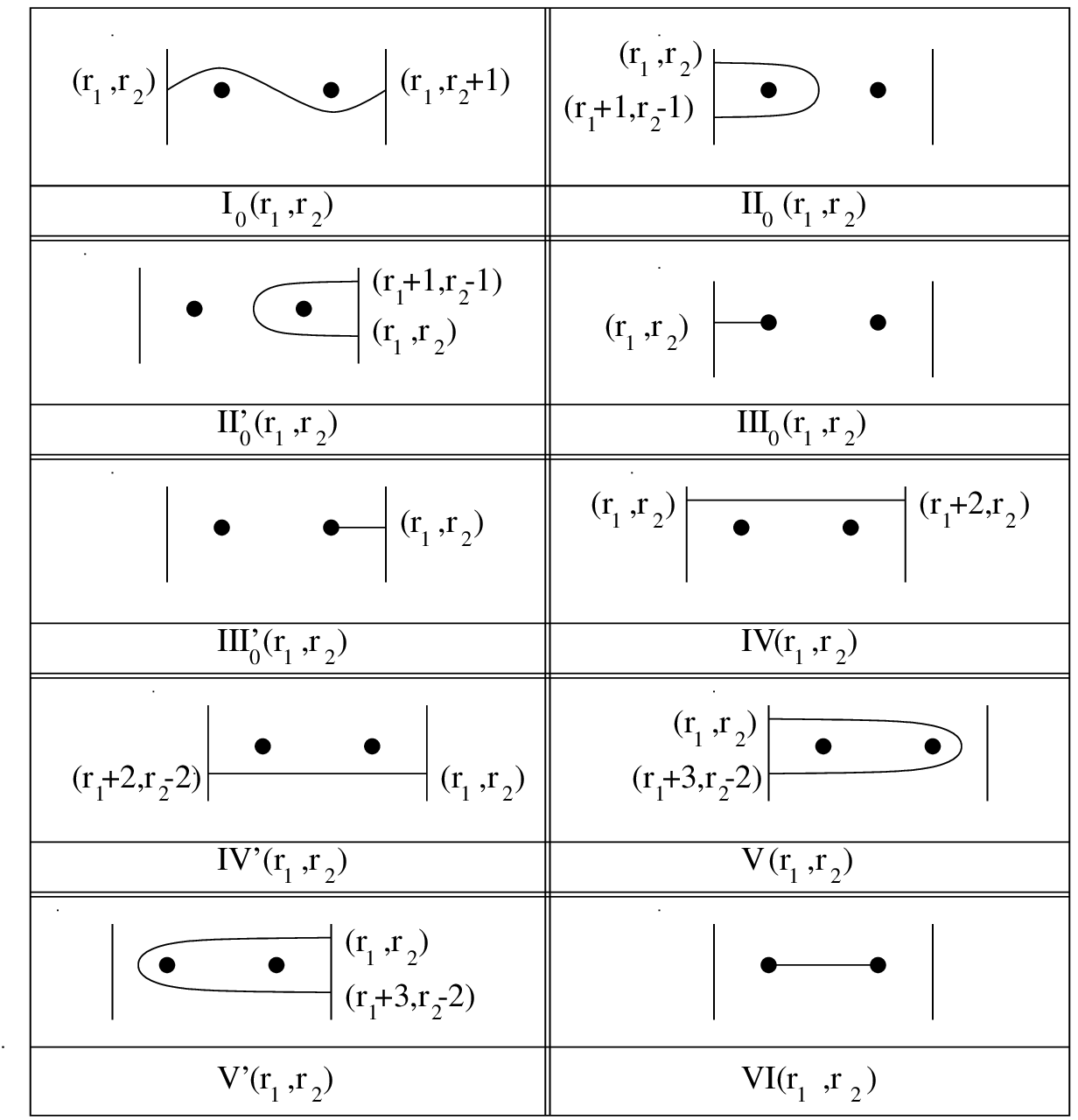} \vspace{-1em} \caption{Isotopy classes of $k$-strings,
for $1 \leq k < m$}
\label{fig:infinite-typesa} \vspace{2em}
\epsfig{file=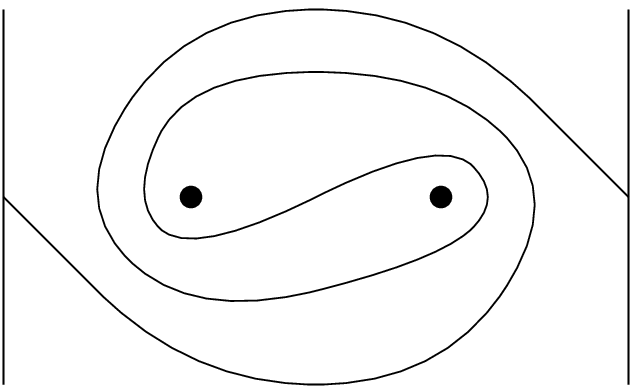} \caption{Type $I_3$} \label{fig:itype-example}
\end{center}
\end{figure}

\begin{figure}
\begin{center}
\epsfig{file=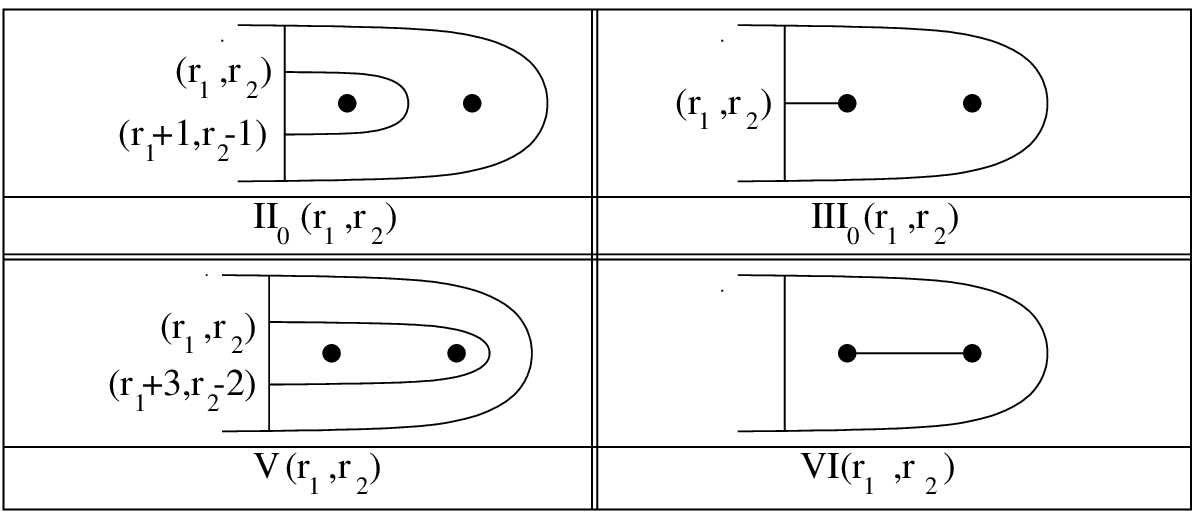} \vspace{-1em} \caption{Isotopy classes of $m$-strings}
\label{fig:infinite-typesb} \vspace{2em}
\epsfig{file=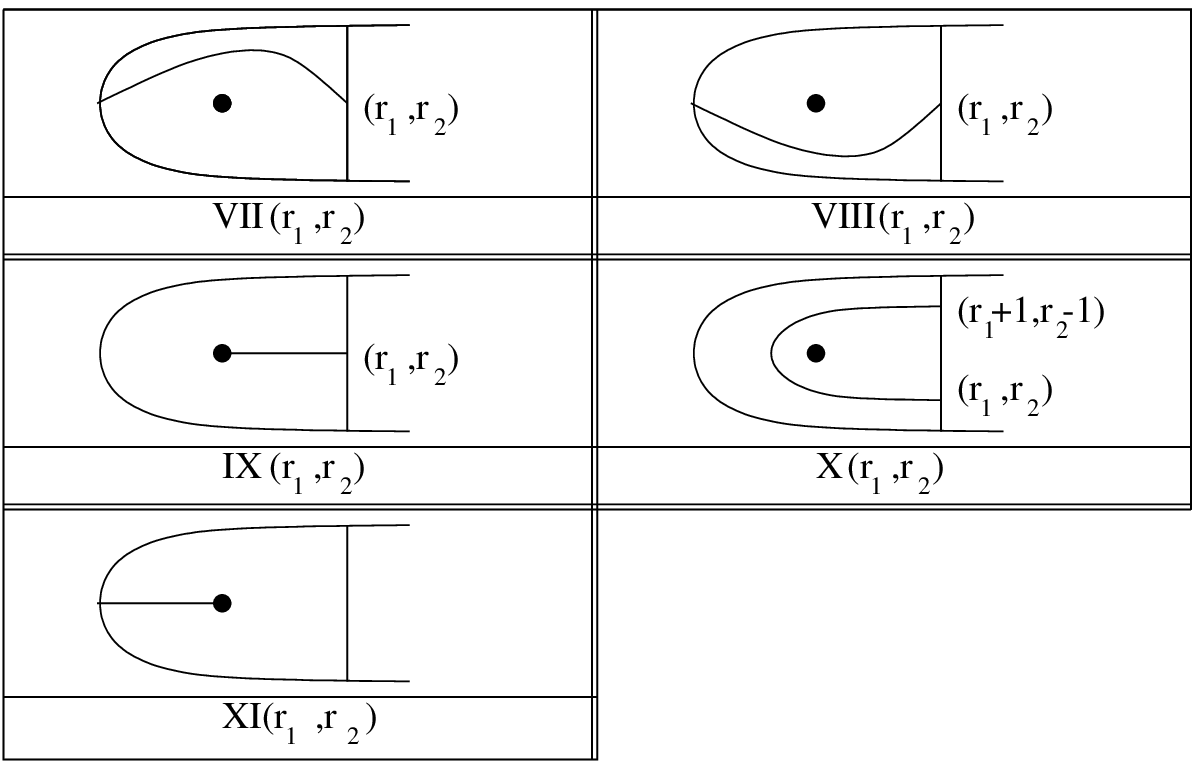} \vspace{-1em} \\ \caption{Isotopy classes of $0$-strings}
\label{fig:infinite-typesc}
\end{center}
\end{figure}

Isotopy classes of $k$-strings can be divided into types as follows:
there are five
infinite families $I_u$, $II_u$, $II'_u$, $III_u$, $III'_u$ $(u \in \Z)$ and five exceptional
types $IV$, $IV'$, $V$, $V'$ and $VI$. The exceptional types, and the members $u = 0$ of the
families, are drawn in Figure\ \ref{fig:infinite-typesa}. The rule for generating the other
members is that the $(u+1)$-st is obtained from the $u$-th by applying $t_k$. As an example,
Figure\ \ref{fig:itype-example} shows type $I_3$.
For $k = m$ there is a similar list, which consists of two families and two exceptional types
(Figure\ \ref{fig:infinite-typesb}).
Finally, for $k = 0$ there are just five exceptional types (Figure\
\ref{fig:infinite-typesc}).

According to our definition, $k$-strings are already assumed to be in ``normal form,'' i.e. to
have minimal intersection with $d_k.$
Define \emph{crossings} and \emph{essential segments} of $k$-strings in the same way as for
admissible curves in normal form. Denote by $\nd(g)$ the set of crossings of a $k$-string $g.$

\begin{lemma} \label{th:almost-normal-i}
If $k>0$, the geometric intersection number
$I(b_k,c)$ can be computed as follows: every $k$-string of $c$
which is of type $I_u$, $II_u$, $II'_u$ or $VI$ contributes $1$, those of type
$III_u$, $III'_u$ contribute $1/2$, and the rest zero. Similarly, for $k = 0$ the types $VIII$
and $X$ contribute $1$, the types $IX$ and $XI$ contribute $1/2$, and $VII$ zero.
\end{lemma}

\proof[Outline of proof.] We will discuss only the case $k>0$. The first step is the following
result, which is quite easy to prove: there is an
isotopy rel $d_0 \cup \dots \cup d_{k-1} \cup d_{k+1} \cup \dots \cup d_m$ which brings $c$
into minimal intersection with $b_k$. This yields a lower bound for the geometric intersection
number under discussion. For instance, if $c \cap (D_k \cup D_{k+1})$ has a $k$-string
 of type $I_0$, this gives rise to one intersection point with $b_k$ which cannot be
removed by our isotopy. Conversely, one can construct an explicit isotopy whereby this lower
bound is attained. \qed

We will now adapt the discussion to bigraded curves. Choose bigradings $\tb_k, \td_k$ of $b_k,
d_k$ for $0 \leq k \leq m$, such that
\[
 \Ibigr(\td_k,\tb_k) = 1 + q_1^{-1}q_2, \quad \Ibigr(\tb_k,\tb_{k+1}) = 1.
\]
These conditions determine the bigradings uniquely up to an overall shift
$\Pchi(r_1,r_2)$. Let $\tc$ be a bigrading of an admissible curve $c$
(as before, $c$ is assumed to be in normal form). Let $a\subset c$ be a connected
component of $c \cap D_k$ for some $k$, and $\ta$ the part of $\tc$
which lies over $a$. Evidently, $\ta$ is determined by $a$ together with the local index
$\mubigr(\tb_{k-1},\ta;z)$ or $\mubigr(\td_k,\ta;z)$ at any point $z \in (d_{k-1} \cup d_k)
\cap a$ (if there is more than one such point, the local indices determine each other). The
resulting classification of the pieces $(a,\ta)$ into various types, together with the local
indices, is indicated by the pairs of integers in Figures \ \ref{fig:typesa}--\ref{fig:typesc}.
For example, if we have a connected component $(a,\ta)$ of type $1(r_1,r_2)$ whose endpoints
are $z_0 \in d_{k-1} \cap a$ and $z_1 \in d_k \cap a$, the local indices are
$\mubigr(\td_{k-1},\ta;z_0) = (r_1,r_2)$ and $\mubigr(\td_k,\ta;z_1) = (r_1+1,r_2)$.

In Section~\ref{subsec:bigraded-curves} we described a canonical lift
$\tilde{f}$ of a diffeomorphism $f\in \Diffeo$
to a diffeomorphism of $\tP.$ Denote by $\tilde{t}_k$ the canonical lift of the twist
$t_k$ along the curve $b_k.$

\begin{prop} \label{t-action}
Diffeomorphisms $\tilde{t}_i$ for $1\le i \le m$ induce a braid group action
on the set of isotopy classes of admissible bigraded curves. In particular, if $\tc$ is an
admissible bigraded curve, the following pairs of bigraded curves are isotopic
  \begin{eqnarray*}
    \tilde{t}_i \tilde{t}_{i+1}\tilde{t}_i (\tc) & \simeq &
    \tilde{t}_{i+1}\tilde{t}_i \tilde{t}_{i+1}(\tc), \\
    \tilde{t}_i \tilde{t}_j(\tc)  & \simeq & \tilde{t}_j \tilde{t}_i(\tc) \hspace{0.3in} |i-j|> 1.
  \end{eqnarray*}
\end{prop}

\vspace{0.2in}

A crossing of $c$ will also be called a \emph{crossing} of $\tc.$ We denote by $\nd(\tc)$
the set of crossings of $\tc,$ it is the same set as $\nd(c).$ Crossings of $\tc$
come with more information than crossings of $c,$ namely, each crossing comes with
a local index, which is an element of $\Z^2.$

The local index $(r_1,r_2)$ of a crossing $x$ of $\tc$ will be denoted
$(x_1,x_2),$ to emphasize that the index is a function of the crossing.
Let $x_0$ denote the index of the vertical curve which contains the crossing:
$x\in d_{x_0}\cap c.$

\emph{Essential segments} of $\tc$ are defined as essential segments of $c$
together with bigradings; the bigrading can be expressed by assigning local
indices to the ends of the segment.

Define a $k$\emph{-string} of $\tc$ as a connected component of $\tc\cap (D_k\cup D_{k+1}),$
together with the bigrading induced from that of $\tc.$
Denote by $\mathrm{st}(\tc,k)$ the set of $k$-strings of $\tc.$

Define a \emph{bigraded} $k$\emph{-string} as a bigraded curve in $D_k\cup D_{k+1}$ which
is a connected component of  $\tc\cap (D_k\cup D_{k+1})$ for some bigraded curve $\tc.$

Figures\ \ref{fig:infinite-typesa}--\ref{fig:infinite-typesc} depict the isotopy classes
of bigraded $k$-strings.
The types $VI$ and $XI$ do not fit immediately into our notational scheme, since they
do not intersect $d_{k-1} \cup d_{k+1}$. Instead, given a bigraded $k$-string $\tg$
with the underlying $k$-string $g$ of type $VI$ or $XI$ ($\tg$ is then a bigraded
admissible curve),  we say that $\tg$ has type $VI(r_1,r_2),$ respectively $XI(r_1,r_2),$ if
$\tg = \chi(r_1,r_2) \tb_k .$

The next result is the bigraded analogue of Lemma \ref{th:almost-normal-i}.

\begin{lemma} \label{th:boring}
Let $(c,\tc)$ be a bigraded curve. Then $\Ibigr(\tb_k,\tc)$ can be
computed by adding up contributions from each bigraded $k$-string of $\tc.$ For $k>0$ the
contributions are listed in the following table:
\[
\begin{array}{c|c|c|c|c|c|c|c|c|c}
 I_0(0,0) & II_0(0,0) & II'_0(0,0)    & III_0(0,0) & III'_0(0,0) & IV & IV' & V & V'
  & VI(0,0) \\
 \hline
 q_1+q_2  & q_1+q_2   & q_1q_2^{-1}+1 & q_2        & 1           & 0  & 0   & 0 & 0
  & 1+ q_2
 \end{array}
\]
and the remaining ones can be computed as follows: to determine the contribution of a
component of type, say, $I_u(r_1,r_2)$ one takes the contribution of $I_0(0,0)$ and multiplies
it by $q_1^{r_1}q_2^{r_2}(q_1^{-1}q_2)^u$. A parallel result holds for $k = 0,$
where the relevant contributions are
\[
\begin{array}{c|c|c|c|c}
 VII(0,0)        & VIII(0,0)       & IX(0,0) & X(0,0) & XI(0,0) \\
 \hline
 0               & q_1q_2^{-1} + 1 & 1       & q_1q_2^{-1} + 1  & 1
\end{array}
\]


\end{lemma}

The proof consists of a series of tedious elementary verifications, which we omit. The only
point worth while discussing is the dependence of the parameter $u$. By combining Lemma
\ref{th:self-shift} with properties \ref{item:bigraded-map} and \ref{item:bigraded-shift} of
$\Ibigr$, one computes that
\[
\Ibigr(\tb_k,\tilde{t}_k(\tc)) = \Ibigr(\tilde{t}^{-1}_k(\tb_k),\tc) =
\Ibigr(\Pchi(1,-1)\tb_k,\tc) = (q_1^{-1}q_2) \Ibigr(\tb_k,\tc).
\]
The same holds for the local contributions coming from a single $k$-string, which
explains the occurrence of the factor $(q_1^{-1}q_2)^u$.

The bigraded intersection number of $\tb_k$ and a bigraded $k$-string $\tg$
(denoted by $\Ibigr(\tb_k,\tg)$) is defined as the contribution, described in
Lemma~\ref{th:boring}, of $\tg$ to the bigraded
intersection number $\Ibigr(\tb_k, \tc).$ Note that $\Ibigr(\tb_k,\tg)$
depends only on the isotopy class of $\tg.$

%
%
%
%
%
%

\section{Admissible curves and complexes of projective modules}
\label{normal-complex}

In this section we prove Theorem~\ref{th:gin} and Corollary~\ref{th:faithful} of the
introduction. In Section~\ref{complex-to-a-curve} we associate a complex of $A_m$-modules
to a bigraded curve. In Section ~\ref{twist-twist}
we prove that this construction relates the braid group actions on bigraded curves
and in the category of complexes. Section~\ref{subsec:bigrdim} interprets bigraded
intersection numbers as dimensions of homomorphisms spaces between complexes associated
to bigraded curves.

Throughout this section $\tc$ denotes a bigraded admissible curve in normal form.

%
%
%
%

\subsection{The complex associated to an admissible curve}
\label{complex-to-a-curve}

\vsp

We associate to $\tc$ an object $L(\tc)$ of the category $\Ccat.$ It is
a complex of projective modules, containing one suitably shifted copy of
$P_i$ for every $i$-crossing of $\tc.$ Start by defining $L(\tc)$ as a
bigraded $A_m$-module:
  \begin{equation}\label{define-L}
     L(\tc) = \bigoplusop{x\in \mathrm{nd}(\tc)} P(x),
  \end{equation}
where  $P(x) = P_{x_0}[-x_1]\{ x_2\}$ (see Section~\ref{normalform}, after
Proposition~\ref{t-action}, for definition of $x_0,x_1,x_2$). For every
$x,y\in \mathrm{nd}(\tc)$ define
  \begin{equation*}
     \partial_{yx}: P(x)\to P(y)
  \end{equation*}
by the following rules

\begin{itemize}
\item If $x$ and $y$ are the endpoints of an essential
      segment and $y_1=x_1+1$ then
      \begin{enumerate}
       \item If $x_0 = y_0$ (then also $x_2= y_2 + 1$) then
       \begin{equation*}
          \partial_{yx}: \hspace{0.1in} P(x) \lra
          P(y)\cong P(x)[-1]\{ 1\}
       \end{equation*}
       is the multiplication on the right by $(x_0| x_0-1 | x_0)\in A_m.$
       \item
       If $x_0= y_0\pm 1$ then $\partial_{yx}$ is the right multiplication
       by $(x_0| y_0)\in A_m;$
      \end{enumerate}
\item otherwise $\partial_{xy}=0.$
\end{itemize}
Now define the differential as
\begin{equation*}
   \partial = \sum_{x,y} \partial_{xy}.
\end{equation*}

\begin{lemma} $(L(\tc),\partial)$ is a complex of graded projective
 $A_m$-modules with a grading-preserving differential.
   \end{lemma}

\proof The equation $\partial^2=0$ follows from relations
$\partial_{zy}\partial_{yx}=0$ for any triple $z,y,x$ of crossings, which are
implied by defining relations in the ring $A_m$:
if there is an arrow from $x$ to $y$ and from $y$ to $z$
then $\partial_{zy}\partial_{yx}: P(x)\to P(z)$ is the right multiplication
by a certain product of generators $(i|i\pm 1)$ of $A_m$, and this
product is equal to zero in all cases (the whole point is that $(i|i+1)$
can never be followed by $(i+1|i)$ since that would mean that we're dealing
with a closed curve around one point of $\Delta,$ and that has been excluded).
$\partial$ is grading-preserving
since each $\partial_{yx}$ is. The latter property is easily checked on
a case by case basis, for each type of essential segments. $\square$

\vsp

Here is a less formal way to describe this complex. Assign the module $P_i$
to each intersection of $c$ with $d_i.$ Notice
(see Figures\ \ref{fig:typesa}--\ref{fig:typesc}) that $x_1= y_1 \pm 1$
if $x$ and $y$ are two ends of an essential segment.
Consequently, every essential segment can be oriented in a canonical way, from
the endpoint $x$ to the endpoint $y$ with $y_1 = x_1+1.$ If the segment has
type $1,$ the orientation is from left to right; if type $1',$ the orientation
is from right to left; if type $2$ or $2',$ the orientation is clockwise.
In general, the orientation is clockwise around the marked point in the
center of the region $D_k$ containing the segment.

There is a natural choice of a homomorphism from $P_{x_0}$ to $P_{y_0}.$
Depending on the type of the segment, the homomorphism is the right
multiplication by $(i|i-1|i),  (i|i \pm 1)$ or $(i\pm 1 | i).$
Namely, if $x_0=y_0=i$ (i.e., both crossings belong to the same vertical
line $d_i$), the homomorphism is the right multiplication by $(i|i-1|i).$
If $x_0=y_0\pm 1=i$ (i.e., the two crossings lie on neighbouring vertical
curves $d_i$ and $d_{i\pm 1}$), the homomorphism is the
right multiplication by $(i|i\pm 1).$

\begin{figure}[h]
\begin{center}
\epsfig{file=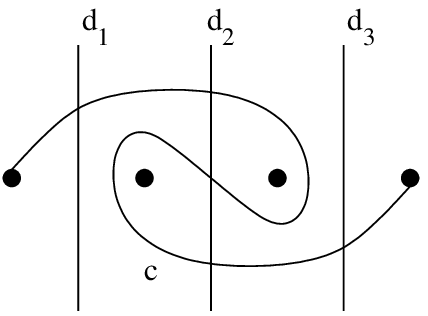}
\vspace{0.1in}
\caption{A curve $c$} \label{fig:curve-c}
\end{center}
\end{figure}

\begin{figure}[h]
\begin{center}
\epsfig{file=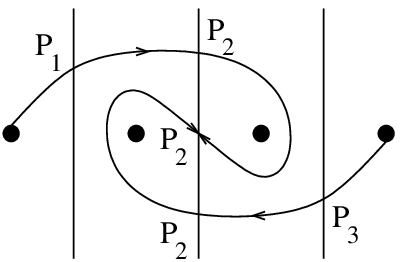}
\vspace{0.1in}
\caption{Orient essential segments and assign projective modules
$P_i$ to intersections of $c$ with $d_i.$}
\label{fig:assign}
\vspace{0.15in}
\end{center}
\end{figure}

\begin{figure}[!h]
\begin{center}
\epsfig{file=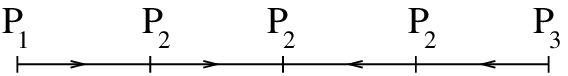}
\vspace{0.1in}
\caption{Stretch out $c$ and remove two inessential segments at the ends
of $c.$}
\label{fig:straighten}
\vspace{0.15in}
\end{center}
\end{figure}

This construction gives us a chain of projective
modules and maps between them, as an example in
Figures\ \ref{fig:curve-c}--\ref{fig:straighten} indicates.

\begin{figure}[!h]
\begin{center}
\epsfig{file=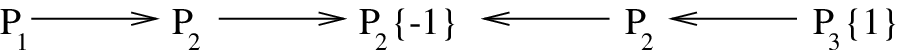}
\vspace{0.1in}
\caption{Replace endpoints of essential segments with the modules
associated to these endpoints  and shift gradings to make homomorphisms
grading-preserving.}
\label{fig:fold}
\end{center}
\end{figure}

The composition of any two consecutive homomorphisms is zero.
The homomorphisms are not grading-preserving, in general. We can fix that
by starting with one of these modules and appropriately shifting the grading
of its neighbours in the chain, then the grading of the neighbours of the
neighbours, etc.

\begin{figure}[h]
\begin{center}
\epsfig{file=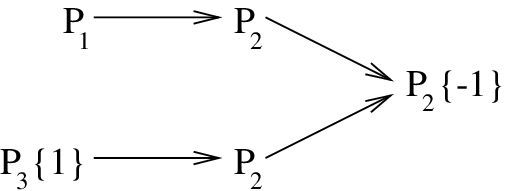} \vspace{-1em}
\caption{Fold the diagram to make all arrows go from left to right.}
\label{fig:shifty-shifts}
\end{center}
\end{figure}

We get different results by starting with different modules, but, in fact,
the difference could only be in the overall shift of the grading.

Fold the resulting diagram made of projective modules and arrows
(homomorphisms) between them so that all arrows go from left to right.
We call this diagram the \emph{folded diagram} of $\tc.$
By summing up modules in different columns of the folded diagram
we get a complex of projective $A_m$-modules. For the example depicted in
Figures \ \ref{fig:curve-c}--\ref{fig:shifty-shifts} the complex is
\begin{equation*}
\dots \lra 0 \lra P_1\oplus (P_3\{ 1\}) \stackrel{\partial}{\lra}
    P_2 \oplus P_2 \stackrel{\partial}{\lra}
    P_2\{ -1\} \lra 0 \lra \dots
\end{equation*}

This object is
defined only up to an overall shift by $[j]\{ k\}.$ This is where
the bigrading of $\tc$ kicks in.
We pick a crossing $x$ of $\tc$ and make the overall shift in the bigrading
so that the projective module $P_{x_0},$ corresponding to $x,$ gets the
bigrading $[-x_1]\{ x_2\}.$ The resulting complex, $L(\tc),$
defined earlier more formally, does not depend on the choice of a crossing $x.$

\vsp

\begin{lemma} \label{lemma:shifts-relation}
An $(r_1,r_2)$ shift of a bigraded curve translates into the $[-r_1]\{ r_2\}$
shift in the category $\Ccat$:
\begin{equation*}
L(\chi(r_1,r_2)\tc) \cong L(\tc)[-r_1]\{ r_2\}.
\end{equation*}
\end{lemma}

\proof Immediately follows  from formula (\ref{define-L}) and
description of the differential in  $L(\tc).$    $\square$

Just like admissible curves, $k$-strings have crossings and essential segments
(see Section~\ref{normalform}). In particular, to a bigraded $k$-string $\tg$
we will associate an object $L(\tg)$ of the category $\Ccat.$ As a bigraded
abelian group,
\begin{equation*}
L(\tc) = \bigoplusop{x\in \mathrm{nd}(g)} P(x),
\end{equation*}
where the sum is over all crossings of $g,$ i.e. intersections of $g$ with
$d_{k-1}\cup d_k \cup d_{k+1}$ and the differential is read off the
essential segments of $g$ in the same way as for admissible curves.

\begin{figure}[h]
\begin{center}
\epsfig{file=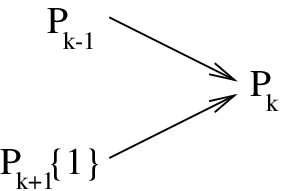}
\vspace{0.05in}
\label{fig:foldstring}
\end{center}
\end{figure}

\emph{Example: } The folded diagram of the bigraded $k$-string $I_0(0,0)$
(see Figure~\ref{fig:infinite-typesa}) is depicted in
Figure~\ref{fig:foldstring}.

The complex $L(I_0(0,0))$ is
$0 \to P_{k-1}\oplus (P_{k+1}\{ 1\}) \to P_k \to 0. $

\vsp

\emph{Remark:} If $\tg$ is a bigraded $k$-string of $\tc,$ then $L(\tg)$
is an abelian subgroup of $L(\tc),$ but not, in general, a subcomplex or a
quotient complex.

\vsp

\emph{Remark:} Lemma~\ref{lemma:shifts-relation} holds for bigraded
$k$-strings as well:
\begin{equation}\label{eqn:shifts-strings}
L(\chi(r_1,r_2)\tg) \cong L(\tg)[-r_1]\{ r_2\}.
\end{equation}

%
%
%
%
%
%

\subsection{Twisting by an elementary braid}
\label{twist-twist}

In the previous section we associated an element of $\Ccat$ to
an admissible bigraded curve. Here we will prove that this map
intertwines the braid group action on bigraded curves with the
braid group action in the category $\Ccat:$

\begin{theorem}
\label{main-theorem}
For a braid $\sigma\in B_m$ and an admissible bigraded curve $\tc$
the objects $\cR_{\sigma}L(\tc)$ and $L(\sigma \tc)$ of the category
$\Ccat$ are isomorphic.
\end{theorem}

It suffices to prove the theorem
when $\sigma$ is an elementary braid, $\sigma \in \{ \sigma_1, \dots,
\sigma_m \}.$ In our notations (see Proposition~\ref{t-action}),
the braid $\sigma_k$ acts on bigraded curves as the Dehn twist
$\tilde{t}_k$ and on an element of $\Ccat$ by tensoring it with the
complex of bimodules $R_k.$ Thus, we want to prove

\begin{prop}
    \label{curve-curve}
    For an admissible bigraded curve $\tc$ and $1\le k \le m,$
    the objects ${R}_k\oo_{A_m} L(\tc)$ and $L(\tilde{t}_k (\tc))$ of
    $\Ccat$ are isomorphic.
\end{prop}

\emph{Proof of the proposition.} Denote $L(\tc)$ by $L.$ The complex
\begin{equation*}
\cR_k L=R_k \oo_{A_m} L
\end{equation*}
is the cone of the map $\U_k(L)\to L.$ We start by
carefully analyzing $\U_k(L)=P_k\oo\hsm _kP\oo_{A_m} L.$ Note that
\begin{equation*}
\U_k (P_j) = \left\{ \begin{array}{ll}
                     P_k \oplus P_k \{ 1\}  & \mbox{ if $j=k$}, \\
                     P_k                    & \mbox{ if $j=k+1$}, \\
                     P_k \{ 1\}             & \mbox{ if $j=k-1$}, \\
                     0                      & \mbox{otherwise}.
                     \end{array} \right.
\end{equation*}
In particular, if a crossing $x$ of $\tc$ does not belong to any $k$-string of
$\tc$ then $\U_k(P(x))=0.$

Pick a bigraded $k$-string $\tg$ of $\tc.$ We have an inclusion of
abelian groups $L(\tg) \subset L(\tc).$

\begin{lemma} \label{applied}
The functor $\U_k,$ applied to this inclusion of groups, produces
an inclusion of complexes
  \begin{equation*}
    \U_k L(\tg) \subset \U_k L(\tc).
  \end{equation*}
\end{lemma}

\emph{Proof of the lemma.} To prove that $\U_k L(\tg)$ is a subcomplex of
$\U_k L(\tc)$ we check that
\begin{equation*}
\U_k(\partial_{xy}): \U_k P(x) \lra \U_k P(y)
\end{equation*}
is zero whenever crossings $x$ and $y$ do not lie in the same $k$-string of $c.$
Indeed, for such $x$ and $y$ there are two possibilities

\begin{enumerate}
\item Either $x$ or $y$ do not lie in a $k$-string
of $c.$ Then $\U_k P(x)=0,$ respectively $\U_k P(y)=0$ and
$\U_k(\partial_{xy})=0.$
\item
$x$ and $y$ belong to different $k$-strings. Then either $\partial_{xy}=0$
or $x_0 = y_0 = k\pm 1$ and the homomorphism
$\partial_{xy}: P(x) \to P(y)$ is the right multiplication by
$(x_0| k | x_0).$ Applying $\U_k$ to this homomorphism produces the trivial
homomorphism.
\end{enumerate}
The lemma follows. $\square$

From this Lemma we deduce

\begin{prop} There is a natural direct sum decomposition of complexes
\begin{equation*}
\U_k L \cong \bigoplusop{\tg\in \mathrm{st}(\tc,k)} \U_k L(\tg).
\end{equation*}
\end{prop}

Earlier we established
a bijection between $k$-strings of $\tc$ and $\tilde{t}_k(\tc)$
(proposition~\ref{bijections})

From the same proposition
we have a natural bijection between $j$-crossings
 of $\tc$ and $\tilde{t}_k (\tc),$ for each $j$ such that  $|j-k|>1.$
If $x$ is such a
crossing of $\tc$ and $x'$ is the associated crossing of $\tilde{t}_k (\tc)$ then
$x_1=x'_1 $ and $x_2 = x'_2.$ Therefore, the complexes $L$ and
$L(\tilde{t}_k(\tc))$ have a large common piece, namely, the abelian
group $\oplusop{x}P(x),$ where the sum is over all crossings $x$ of $\tc$
such that $|k-x_0|>1.$ Denote this direct sum by $W.$

We want to contruct a homotopy equivalence
\begin{equation*}
\cR_k(L) \cong L(\tilde{t}_k (\tc)).
\end{equation*}
We will obtain it as a composition of homotopy equivalences
\begin{equation*}
\cR_k(L) = T_0 \cong T_1 \cong \dots \cong T_{s-1}\cong
T_s = L(\tilde{t}_k(\tc)).
\end{equation*}
Every intermediate complex $T_i$ will have a structure similar to that of
$\cR_k(L)$ and $L(\tilde{t}_k(\tc)).$ As an $A_m$-module, $T_i$ will be
isomorphic to the direct sum of  $W$ and $A_m$-modules underlying certain
complexes $T_i(\tg),$ homotopy equivalent to $\cR_k(\tg),$ over all
bigraded $k$-strings $\tg$ in $\tc:$

\begin{equation*}
 T_i = W\oplus (\bigoplusop{\tg\in \mathrm{st}(\tc,k)} T_i(\tg)).
\end{equation*}
Elementary homotopies
$T_i \cong T_{i+1}$ will be of two different kinds:

(i) homotopies ``localized'' inside $T_i(\tg),$ for some
bigraded $k$-string $\tg.$ Namely, the homotopy will be the identity on
$W$ and on
$T_i(\tilde{h}),$ for all $k$-strings $\tilde{h}$ of $\tc$ different from
$\tg.$ The complexes $T_i(\tg)$ and $T_{i+1}(\tg)$ will differ as follows:
$T_{i+1}(\tg)$ will be the quotient of $T_i(\tg)$ by an acyclic subcomplex,
or $T_{i+1}(\tg)$ will be realized as a subcomplex of $T_i(\tg)$ such that
the quotient complex is acyclic.
The ``localized'' homotopy $T_i(\tg)\cong T_{i+1}(\tg)$ will respect the
position of these complexes as subcomplexes of $T_i $ and $T_{i+1}$ and
will naturally extend to homotopy equivalence $T_i \cong T_{i+1}.$

(ii) Isomorphisms of complexes $T_i \cong T_{i+1}.$ These will come from
natural isomorphisms of $A_m$-modules that underlie complexes $T_i$ and
$T_{i+1}$ and a direct sum decomposition of the module $T_i= (T_i)^+
\oplus (T_i)^-.$ The isomorphism will take $(v^+,v^-)\in T_i$ to
$(v^+,-v^-).$

\vsp

From now on we proceed on a case by case basis, working with different
types of bigraded $k$-strings
(Figures~\ref{fig:infinite-typesa}-\ref{fig:infinite-typesc}).
We start with some easy cases.

\vsp

{\bf Case 1.}
If $\tg$ is of type $VI,$
then $\tc\simeq \tc_k,$ up a shift in the bigrading.
The half twist $\tilde{t}_k(\tc)$ shifts the bigrading by $(-1,1).$
On the other hand, the complex $R_k \oo_{A_m} P_k$ is chain homotopic
to $P_k[1]\{ 1\}.$ The topological
shift matches the algebraic (see Lemma~\ref{lemma:shifts-relation}
and formula (\ref{eqn:shifts-strings})). This gives an isomorphism
$L(\tilde{t}_k (\tc_k))\cong \cR_k (L(\tc_k)).$

\vsp

%
%
%
%

{\bf Case 2.}
If $\tg$ has one of the types $IV, IV', V$ or
$V'$ then the half-twist $\tilde{t}_k$ preserves $\tg.$
At the same time, we have

\begin{lemma} If $\tg$ is of type $IV,IV',V$ or $V',$ then
the complex $\U_k L(\tg)$ is acyclic.
\end{lemma}

\emph{Proof of the lemma.} Explicit computation. For instance, if $\tg$ has type
$IV,$ then $L(\tg)$ is the complex
\begin{equation}
\label{not-quite-exact}
0 \lra P_{k-1}\{r_2\} \lra  P_k\{r_2\} \lra P_{k+1}\{r_2\} \lra 0.
\end{equation}
The functor $\U_k$ is the tensor product with the bimodule
$P_k \oo \hsm _kP.$ First we tensor over $A_m$ with the right module
$\hsm _kP,$ which gets us a complex of abelian groups. Then
we tensor the result with $P_k$ over $\Z.$ The tensor product of the
complex (\ref{not-quite-exact})
with $\hsm _kP$ over $A_m$ is an acyclic complex of abelian
groups (exercise for the reader). The lemma follows for $\tg$ of type $IV.$
Similar computations establish the lemma for the other three types of $\tg.$
$\square$

Since $\tilde{t}_k(\tg) \simeq \tg,$ we have a natural isomorphism of complexes
  $L(\tilde{t}_k(\tg)) \cong L(\tg).$
Moreover, $\U_k (L(\tg)),$ which is
a direct summand  of the complex
$\U_k (L),$ is acyclic, and we see that the complex $\cR_k L$ is homotopic
to its subcomplex obtained by throwing away $\U_k (L(\tg))$ for all $\tg$ of
types IV,IV',V and V'.

\vsp

%
%
%
%

{\bf Case 3.} Types $III$ and $III'.$

\vspace{0.1in}

Assume $\tg$ has type $III_u(r_1,r_2).$ The homotopies constructed below commute
with shifts in the bigrading, and,  without loss of generality,
we set $r_1=r_2=0.$ We will also set $r_1=r_2=0$ in our treatment of
types $II,II'$ and $I.$

\begin{figure}[h]
\begin{center}
\epsfig{file=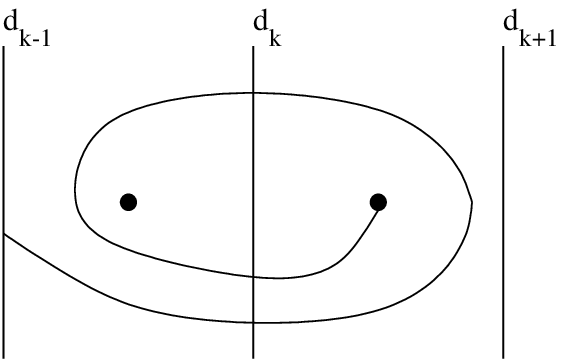} \vspace{-1em}
\label{fig:typeIII3}
\end{center}
\end{figure}

Consider first the case $u\ge 0.$ A $k$-string of type $III_3$ is depicted
on Figure~\ref{fig:typeIII3}.

The complex $L(\tg)$ is
\begin{equation*}
0 \lra P_k\{ u\} \lra P_k\{ u-1\} \lra \dots \lra P_k\{ 1\} \lra P_{k-1}\lra 0.
\end{equation*}

Let us understand the relation between $L(\tg)$ and $L.$
The complex $L$ contains $L(\tg)$ as an abelian subgroup.
The part of the folded diagram of $\tc$ that contains the folded diagram
of $\tg$ has the form

\begin{equation*}
  P_k\{ u\} \to P_k\{ u-1\} \to \cdots \to P_k\{ 1\} \to P_{k-1}
  \leftrightarrow \nabla ,
\end{equation*}
where $\nabla$ denotes the complement to $L(\tg)$ in $L.$
More precisely, $\nabla$ denotes the direct sum
$\oplusop{x} P(x)\subset L$ over all crossings $x$ of $\tc$ that do not belong
to $\tg.$ Notation $P_{k-1}\leftrightarrow \nabla$ means that
$\nabla$ relates to $L(\tg)$ through the module $P_{k-1}.$ Namely,
in the folded diagram of $L$ there
might be an arrow $P(x)\to P_{k-1}$ or an arrow $P_{k-1}\to P(x),$
where $x$ is a crossing of $\tc$ which does not belong to $\tg.$

We will use this notation throughout the rest of the proof of
Proposition~\ref{curve-curve}. $\nabla$ will denote either such a complement
to $L(\tg)$ in $L,$ or the complement to $\cR_k(L(\tg))$ in $\cR_k L,$
or any similar complement, and the meaning will always be clear from the
context. Whenever we modify $\cR_k L(\tg)$
to a homotopy equivalent complex $L',$ we will do so relative to $\nabla,$
so that our homotopy equivalence extends to a homotopy equivalence
between $\cR_k L$ and the complex naturally built out of $L'$ and
$\nabla.$ When we are in such situation, we say that the homotopy
\emph{respects} $\nabla.$

The complex $\cR_k(L(\tg)),$ written as a bicomplex, and considered
as a part of $\cR_k L,$  has the form
\[
\begin{array}{ccccccccc}
   P_k\{ u+1\} \oplus P_k\{ u\} & \to & \dots & \to &
    P_k\{ 2\} \oplus P_k\{ 1\} & \to & P_k \{ 1\}  & &   \\
   \downarrow                   &      &  \dots &     &
       \downarrow              &      &   \downarrow  & &  \\
   P_k\{ u\}                    & \to & \dots & \to &
    P_k \{ 1\}                 & \to & P_{k-1} & \leftrightarrow & \nabla .
\end{array}
\]
Note that, in view of Lemma~\ref{applied}, $\nabla$ connects
 to $\cR_k(L(\tg))$ only through the module $P_{k-1}$ in the bottom
right corner of the diagram.

Complex $\U_k(L(\tg))$ (the top row in the above diagram) is isomorphic
to the direct sum of $P_k\{ u+1\}$ and acyclic complexes
\begin{equation*}
0\lra P_k\{ i\}\stackrel{\mathrm{id}}{\lra} P_k\{ i\} \lra 0,
\hspace{0.3in} 1\le i \le u
\end{equation*}
(we leave this as an exercise for the reader, it follows at once after the
differential in $\U_k(L(\tg))$ is written down),
so that  $\cR_k( L(\tg))$ is homotopic to
 \begin{equation*}
 P_k\{ u+1\} \lra  P_k\{ u\} \lra \dots \lra P_k \{ 1\} \lra P_{k-1}
 \longleftrightarrow \nabla,
 \end{equation*}
which is exactly $L(\tilde{t}_k(\tg))$ (note that the homotopy respects
$\nabla$).

\vspace{0.2in}

We next treat the case $III_u$ where $u<0.$ As before, set $r_1=r_2=0,$
this does not reduce the generality since our homotopies will commute with
bigrading shifts.

\begin{figure}[h]
\begin{center}
\epsfig{file=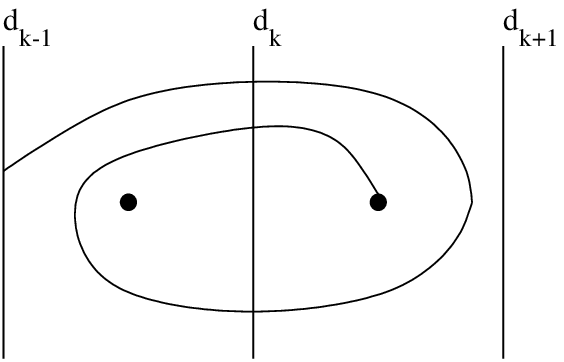} \vspace{-1em}
\caption{$k$-string of type $III_{-3}$}
\label{fig:typeIIIminus3}
\end{center}
\end{figure}

The complex $L(\tg)$ (considered as a part of $L$) has the form (compare to
Figure~\ref{fig:typeIIIminus3})
\begin{equation*}
   \nabla \longleftrightarrow P_{k-1}\lra P_k \lra P_k\{ -1\} \lra
   \cdots \lra P_k\{ 1+u\}.
\end{equation*}

The complex $\cR_k(L(\tg)),$ in relation to $\cR_k(L),$ has the form
\[
\begin{array}{ccccccccc}
   P_k\{1\} & \to &
   {\left(\begin{array}{c}P_k  \\   \oplus   \\   P_k\{1\}  \end{array}
    \right)}
   & \to &
   {\left(\begin{array}{c}P_k\{-1\}\\  \oplus   \\   P_k  \end{array}
    \right)}
   & \to &
   \dots & \to &
   {\left(\begin{array}{c}P_k\{1+u\}\\   \oplus   \\   P_k\{2+u\}  \end{array}
    \right)}  \\
  \downarrow & & \downarrow & & \downarrow & & & &   \downarrow \\
  P_{k-1}  & \to & P_k & \to & P_k\{ -1\} & \to &  \dots & \to & P_k\{ 1+u\}
  \\
  \updownarrow & &  &  &  & & & & \\
  \nabla       & &  &  &  & & & &
  \end{array}
\]

The complex $\U_k(L(\tg))$ (top row in the diagram above) is isomorphic to
the direct sum of $P_k\{ 1+u\}$ and acyclic complexes
  \begin{equation*}
  0\lra P_k \{ i\}\stackrel{\mathrm{id}}{\lra} P_k\{ i\} \lra 0,
  \hspace{0.2in}\mbox{ for }\hspace{0.2in}  1\ge i \ge 2+u.
  \end{equation*}
Therefore, $\cR_k(L(\tg))$ is homotopic to the folding of
 \[
 \begin{array}{ccccccccccc}
   & &  & & & & & & & & P_k\{1+u\}  \\
   & &  & & & & & & & & \mbox{\scriptsize{id}}\downarrow \\
   \nabla & \leftrightarrow & P_{k-1} & \to & P_k & \to & P_k\{ -1\} & \to &
   \dots & \to & P_k\{ 1+u\} \\
 \end{array}
 \]
 (an explicit computation shows that the vertical arrow is the identity map).
This complex has an acyclic subcomplex
 \begin{equation*}
   0 \to P_k\{ 1+u\} \stackrel{\mathrm{id}}{\lra} P_k\{ 1+u\} \to 0,
 \end{equation*}
and the quotient is isomorphic to $L(\tilde{t}_k(\tg)).$ This sequence of
two homotopy equivalences respects $\nabla.$ Our consideration of case
$III$ is now complete.
Type $III'$ can be worked out entirely parallel to $III,$ and we omit
the computation.

\vspace{0.2in}

%
%
%
%

{\bf Case 4.} Types $II$ and $II'.$ Let $\tg$ has type $II_u(r_1,r_2).$
Without loss of generality we can specialize to $r_1=r_2=0.$

\begin{figure}[h]
\begin{center}
\epsfig{file=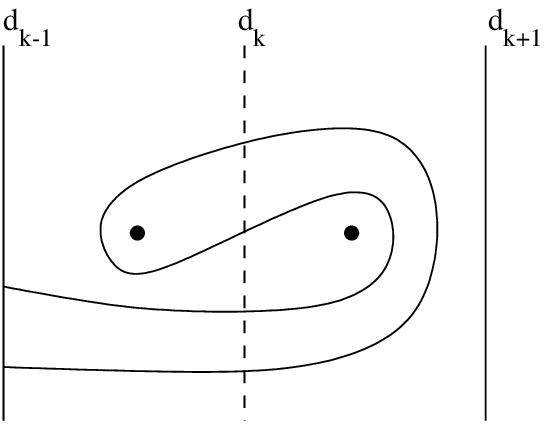} \vspace{-1em}
\label{fig:typeII2}
\end{center}
\end{figure}

We start with the case $u>0.$ Figure~\ref{fig:typeII2} depicts
an example of a $k$-string of that type.

The complex $L(\tg)$ has the form (the first line of arrows continues into
the third, and the second into the fourth)

\[
\begin{array}{cccccclcc}
 P_k\{ u\} & \lra      & P_k\{ u-1\}   & \lra & \dots & & & &  \\
           & \searrow  &               &      &       & & & &  \\
           &           & P_k\{ u-1\}   & \lra & \dots & & & &  \\
           &           & \hspace{0.1in}&      &       & & & &  \\
 \dots     & \lra      &     P_k\{ 1\} & \lra &
 P_{k-1}   & \longleftrightarrow & \nabla         &         &  \\
           &           &               &      &       & & & &  \\
 \dots     & \lra      &     P_k\{ 1\} & \lra &
 P_k       & \lra                & P_{k-1}\{ -1\} &
   \longleftrightarrow & \nabla.
\end{array}
\]

We next write down the complex $\U_k(L(\tg)):$

\[
\begin{array}{cccccclcc}
 P_k\{ u\}\oplus P_k\{ u+1\} & \to      & P_k\{ u-1\}\oplus P_k\{ u\}
           & \to & \dots & & & &  \\
           & \searrow  &               &      &       & & & &  \\
           &           & P_k\{ u-1\}\oplus P_k\{u\}
           & \to & \dots & & & &  \\
           &           & \hspace{0.1in}&      &       & & & &  \\
 \dots     & \to      &     P_k\{ 1\}\oplus P_k\{2\} & \to &
 P_{k}\{1\}   & \to & 0        &         &  \\
           &           &               &      &       & & & &  \\
 \dots     & \to      &     P_k\{ 1\}\oplus P_k\{2\} & \to &
 P_k\oplus P_k\{1\}       & \to                & P_k &
   \to & 0
\end{array}
\]

It is easily checked that this complex is isomorphic to direct sum of
acyclic complexes $0\to P_k\{u\}\stackrel{\mathrm{id}}{\lra}P_k\{u\}\to 0,$
(two complexes for each $i,1\le i\le u$ and one for $i=0$) and the complex
\begin{equation}
\begin{array}{clccc}
\label{asubcomplex}
 P_k\{ u\}\oplus P_k\{ u+1\} & 
           \stackrel{\scriptstyle{\psi_1}}{\lra}
           & P_k\{ u\}  & \lra & 0  \\
           & \searrow \scriptstyle{\psi_1}  &            &      &    \\
           &                           & P_k\{u\}   & \lra & 0,
\end{array}
\end{equation}
where $\psi_1$ is the projection onto the first summand, taken with the minus
sign.

Denote by $L'$ the subcomplex of $\cR_k(L(\tg))$ which, as an abelian group,
in the direct sum of $L(\tg)$ and (\ref{asubcomplex}). Note that $L'$ is
homotopy equivalent to $\cR_k(L(\tg)),$ and that this homotopy equivalence
respects $\nabla.$ Complex $L'$ consists of the "central part"
\begin{equation}
\begin{array}{ccccccc}
\label{complex-prime}
   &  P_k\{ u\} &  \stackrel{\psi_1}{\leftarrow} &
    P_k\{ u\}\oplus P_k\{ u+1\} & \stackrel{\psi_1}{\to} &
      P_k\{ u\} &      \\
   &  \downarrow &            & \scriptstyle{\psi_2}\downarrow        &     &
      \downarrow &     \\
   \leftarrow &  P_k\{u-1\} & \leftarrow & P_k\{ u\}  & \to  &
   P_k\{ u-1\} & \to
\end{array}
\end{equation}
and two chains of $P_k\{i\}$'s going to the left and right of the center.
$\psi_2$ is given by
\begin{equation}
\label{psi-2}
\psi_2(a,b)= a + b(k|k-1|k), \hspace{0.2in} a\in P_k\{u\}, \hspace{0.1in}
                                             b\in P_k\{ u+1\}.
\end{equation}

The differential in this complex is injective on $P_k\{u\}$ situated
in the center of the top row of (\ref{complex-prime}).
Thus, $P_k\{u\}$  generates an acyclic subcomplex
$0 \to P_k\{ u\} \stackrel{\mathrm{id}}{\lra}P_k\{u\}\to 0$ of $L'.$
The quotient of $L'$ by this acyclic subcomplex is isomorphic to
(the diagram below only depicts the central part of the complex, don't forget
to add two tails of $P_k$'s)

\[
\dots \longleftarrow  P_k\{u-1\} \stackrel{-(k|k-1|k)}{\longleftarrow}
 P_k\{ u\} \stackrel{(k|k-1|k)}{\longrightarrow} P_k\{ u-1\}
\longrightarrow \dots
\]

Denote this quotient by $L''.$ The $A_m$-modules underlying complexes
$L''$ and $L(\tilde{t}_k(\tg))$ are canonically isomorphic. The
two differentials differ only in one place: in $L''$ there is
a minus sign in the left arrow. We can get rid of this minus sign with the
isomorphism $L''\cong L(\tilde{t}_k(\tg))$ that takes $a\in P(x)$ to $\pm a$
depending on whether $P(x)$ is to the left or right of the ``minus''
arrow.

\vspace{0.3in}

We next treat the case $u=0.$ The complex $L(\tg)$ is
\begin{equation*}
\nabla \leftrightarrow P_{k-1} \lra P_{k-1}\{ -1\} \leftrightarrow \nabla,
\end{equation*}
and the complex $\cR_k(L(\tg)),$ for $\tg$ of type
$II_0,$ is the folding of
\begin{equation}
\begin{array}{ccccccc}
\label{2by2-start}
 &  &  P_k\{ 1\} & \stackrel{0}{\lra} &  P_k &  &  \\
 &  &  \downarrow &  &  \downarrow & & \\
\nabla & \leftrightarrow & P_{k-1} & \lra & P_{k-1}\{ -1\} & \leftrightarrow &
\nabla .
\end{array}
\end{equation}

We want to find a homotopy equivalence with the complex
$L(\tilde{t}_k(\tg)),$ which is
\begin{equation}
\begin{array}{ccccccc}
\label{2by2-finish}
 &  &  P_k\{ 1\} & \stackrel{-(k|k-1|k)}{\lra} &  P_k &  &  \\
 &  &  \downarrow &  &  \downarrow & & \\
\nabla & \leftrightarrow & P_{k-1} &  & P_{k-1}\{ -1\} & \leftrightarrow &
\nabla .
\end{array}
\end{equation}
In fact, these two complexes are isomorphic, via the isomorphism which
takes an element
\begin{equation*}
{\left( \begin{array}{cc} x & y \\  z & w \end{array}\right)},
\hspace{0.2in} x\in P_k\{1\}, y\in P_k, z\in P_{k-1}, w\in P_{k-1}\{-1\}
\end{equation*}
of (\ref{2by2-finish}) to
\begin{equation*}
{\left( \begin{array}{cc} x & -y - z(k-1|k)\\  z & -w
\end{array}\right)}
\end{equation*}
of the complex (\ref{2by2-start}).

Since our map takes $w\in P_{k-1}\{-1\}$ to $-w,$ when extending it
to the remaining part of the complex we will send $v$ to $-v$ for any
$v$ that belongs to $\nabla$ on the right side of the diagram
(\ref{2by2-finish}). This is to insure that our map commutes with
differentials, i.e., it is a map of complexes.

\vspace{0.3in}

Case $u<0.$ Let $\tg$ be of type $II_u$ for $u<0.$
The complex $L(\tg)$ is the folding of

\[
\begin{array}{cclcccccc}
            &            &  \nabla & \longleftrightarrow & P_{k-1}\{-1\} &
 \lra & P_k\{ -1\} & \lra & \dots  \\
   &           &               &      &       & & & &  \\
 \nabla & \longleftrightarrow & P_{k-1} & \lra & P_k & \lra & P_k\{-1\}  &
 \lra  & \dots   \\
   &           & \hspace{0.05in}&     &       &  & & &  \\
 &     &    &    &   \dots & \lra  & P_k\{ u+1\} & \lra         & P_k\{u\} \\
 &     &    &    &         &       &             & \nearrow     &  \\
 &     &    &    &   \dots & \lra  & P_k\{u+1\}  &              &
\end{array}
\]

The situation is very similar to the case $u>0.$ The complex $\U_k(L(\tg))$
is isomorphic to the direct sum of acyclic complexes and of the complex
which is the folding of
\begin{equation}\label{letsfold}
 P_k\{ u+1\} \stackrel{\psi_1}{\lra}
    P_k\{ u\}\oplus P_k\{ u+1\} \stackrel{\psi_1}{\longleftarrow} P_k\{ u+1\},
\end{equation}
where $\psi_1=(0,-\mbox{id}).$

The complex $L',$ which is the direct sum of $L(\tg)$ and
(\ref{letsfold}) inside $\cR_k(L(\tg)),$
consists of the central part, which is the folding of
\begin{equation}
\begin{array}{ccccccc}
\label{central-part}
   &  P_k\{ u+1\} &  \stackrel{\psi_1}{\lra} &
    P_k\{ u\}\oplus P_k\{ u+1\} & \stackrel{\psi_1}{\longleftarrow} &
      P_k\{ u+1\} &      \\
   &  \mbox{\scriptsize{id}}\downarrow &  & \scriptstyle{\psi_2}
   \downarrow  &  &  \mbox{\scriptsize{id}}\downarrow &  \\
   \to &  P_k\{u+1\} & \to & P_k\{ u\}  & \leftarrow  &
   P_k\{ u+1\} & \leftarrow
\end{array}
\end{equation}
and of two long tails of $P_k\{i\}$'s. The map $\psi_2$ is given by the
formula (\ref{psi-2}).

The module $P_k\{u\}$ in the center of the top row of (\ref{central-part})
generates an acyclic complex inside $L'.$ The quotient of $L'$ by this
acyclic subcomplex is isomorphic to $L(\tilde{t}_k(\tg))$ (things work here
in the same way as in the case $u<0,$ treated earlier).

%
%
%
%

{\bf Case 5.} Type $I.$ This case is very similar to $II.$ Write down
the complex $L(\tg),$ compute $\U_k(L(\tg)),$ throw away long tails of
acyclic subcomplexes in $\U_k(L(\tg)),$ we'll be left with
(\ref{asubcomplex}) (case $u>0$) or (\ref{letsfold}) (case $u\le 0$).
Form the direct sum $L'$ of $L(\tg)$ and what's left of $\U_k(L(\tg)).$
Complex $L'$ is a subcomplex of $\cR_k(L(\tg)).$ Quotient out
$L'$ by an acyclic subcomplex in the same way as in case $II.$ The
result is isomorphic to $L(\tilde{t}_k(\tg)).$

\vsp

This completes our case-by-case analysis. For each $k$-string $\tg$ of
$\tc$ we constructed a homotopy equivalence
$\cR_k(L(\tg))\cong L(\tilde{t}_k(\tg))$ which extends naturally to
$\nabla.$ Putting all these equivalences together we obtain
a homotopy equivalence $\cR_k(L)\cong L(\tilde{t}_k(\tilde{c})).$
The proof of proposition~\ref{curve-curve} is complete. $\square$

\vsp

Proposition~\ref{curve-curve} implies (taking $\tc= \tb_k$)

\begin{corollary}
\label{main-corollary}
For a braid $\sigma\in B_m$ objects $\cR_{\sigma}P_k$
  and $L(\sigma \tilde{b}_k)$ of the category $\Ccat$ are isomorphic.
\end{corollary}

%
%
%
%
%
%

\subsection{Bigraded intersection numbers and dimensions of homomorphism
spaces} \label{subsec:bigrdim}

The main result of this section:

\begin{prop} \label{th:dimequals} For $\sigma, \tau\in B_{m+1}, s_1,s_2\in \Z$
and $0\le k,j\le m$ the abelian group
  \begin{equation*} \mathrm{Hom}_{\Ccat}
  (\cR_{\tau}P_k, \cR_{\sigma}P_j[s_1]\{ -s_2\})
  \end{equation*}
 is free. The Poincar{\'e} polynomial
  \begin{equation*}
  \sum_{s_1,s_2}\mathrm{rk}(\mathrm{Hom}_{\Ccat}
  (\cR_{\tau}P_k, \cR_{\sigma}P_j[s_1]\{ -s_2\})) q_1^{s_1} q_2^{s_2}
  \end{equation*}
 is equal to the bigraded intersection number $\Ibigr(f_{\tau}(\tilde{b}_k),
 f_{\sigma}(\tilde{b}_j)).$
\end{prop}

\proof The braid group acts in the category $\Ccat.$ In particular,
\begin{equation*}
  \mbox{Hom}_{\Ccat}(\cR_{\tau}P, \cR_{\tau}Q)\cong
  \mbox{Hom}_{\Ccat}(P,Q).
\end{equation*}
The bigraded intersection number is also preserved by the braid group action:
\begin{equation*}
  \Ibigr(f_{\tau}(\tc_1), f_{\tau}(\tc_2))=\Ibigr(\tc_1,\tc_2).
\end{equation*}
It suffices, therefore, to prove the proposition in the case $\tau =1.$
Moreover, Corollary~\ref{main-corollary} tells us that we can substitute
$L(\sigma\tb_j)$ for $\cR_{\sigma}P_j$ in this proposition.

 \begin{lemma}
  \label{contributions}
 For a bigraded admissible curve $\tc$ there is a natural isomorphism
  \begin{equation*}
    \mathrm{Hom}_{\Ccat}(P_k, L(\tc)[s_1]\{ -s_2\}) \cong
    \bigoplusop{\tg\in \mathrm{st}(\tc,k)}
    \mathrm{Hom}_{\Ccat}(P_k,L(\tg)[s_1]\{ -s_2\})
  \end{equation*}
 for all $s_1,s_2\in \Z$ and $0\le k \le m.$
 \end{lemma}
In other words, homomorphisms from $P_k$ to the complex associated to an
admissible curve come from homomorphisms of $P_k$ to complexes associated to
$k$-strings of the curve.

\emph{Proof of the lemma.} The complex $\bigoplusop{\tg\in \mathrm{st}(\tc,k)}
L(\tg)$ is obtained from $L(\tc)$ by throwing away all modules $P_j$ for
$|j-k|>1,$ and by changing differentials
  \begin{equation*}
    \partial_{xy}: P(x)\to P(y)
   \end{equation*}
to zero whenever $x$ and $y$ belong to boundaries of two different
$k$-strings of $\tc.$ The lemma follows since
  \begin{equation*}
    \mbox{Hom}_{\Ccat}(P_k, P_j[s_1]\{ s_2\})=0
  \end{equation*}
for all $s_1,s_2\in \Z$ and all $j$ such that $|j-k|>1;$ and
since maps
 \begin{equation*}
   \mbox{Hom}_{\Ccat}(P_k,P(x))\lra \mbox{Hom}_{\Ccat}(P_k,P(y))
 \end{equation*}
induced by $\partial_{xy}$ as above are zero maps.$\square$

We know that $\Ibigr(\tilde{b}_k, \tc)$ is obtained by summing up
contributions from all $k$-strings of $\tc$ (see Lemma~\ref{th:boring}).
Similarly, Lemma~\ref{contributions} tells us that the space
of homomorphisms from $P_k$ to a (possibly shifted) $L(\tc)$ is the
direct sum of contributions corresponding to $k$-strings of $\tc.$
Therefore, it suffices to prove

\begin{lemma}
\label{free-abelian}
For any bigraded $k$-string $\tg$ the abelian group
   \begin{equation*} \mathrm{Hom}_{\Ccat}(P_k, L(\tg)[s_1]\{ -s_2\})
   \end{equation*}
  is free. The Poincar{\'e} polynomial
   \begin{equation*}
   p(\tg)\stackrel{\scriptstyle{\mathrm{def}}}{=}\sum_{s_1,s_2}
   \mathrm{rk}(\mathrm{Hom}_{\Ccat}
   (P_k, L(\tg)[s_1]\{ -s_2\})) q_1^{s_1} q_2^{s_2}
  \end{equation*}
 is equal to the bigraded intersection number $\Ibigr(\tilde{b}_k,\tg).$
\end{lemma}

\emph{Proof of the lemma.}
Let us examine how shifting a bigraded $k$-string by $(r_1,r_2)$
changes its Poincar{\'e} polynomial. If $\tg'= \chi(r_1,r_2)\tg$ then
$L(\tg')= L(\tg) [-r_1]\{ r_2\}$ and $p(\tg')= q_1^{r_1}q_2^{r_2}p(\tg).$
This shift matches the corresponding shift in the bigraded intersection
number:
  \begin{equation*}
    \Ibigr(\tilde{b}_k,\tg')= q_1^{r_1} q_2^{r_2}\Ibigr(\tilde{b}_k,\tg).
  \end{equation*}
Therefore, it suffices to prove the lemma in the case when the bigraded
$k$-string $\tg$ has its parameters $(r_1,r_2)$ (see
Figures\ \ref{fig:infinite-typesa}--\ref{fig:infinite-typesc}) set to $0.$

We start by treating the case $0< k < m.$ Figure~\ref{fig:infinite-typesa}
depicts 10 possible types of bigraded $k$-strings,
and we have already reduced to the $r_1=r_2=0$ case.
Notice that types $I,II, II',III$ and $III'$ have an extra
integral parameter $u.$ Denote by $Y$ one of the above $5$ types.
The $k$-string $Y_u$ is given by twisting $Y_0$ by
$(\tilde{t}_k)^u.$ Since
  \[
  \begin{array}{l}
   \Ibigr(\tilde{b}_k,  Y_u)= \\
    \Ibigr(\tilde{b}_k , (\tilde{t}_k)^u Y_0 )=
    \Ibigr((\tilde{t}_k)^{-u} \tilde{b}_k, Y_0) =
    \Ibigr(\chi(u,-u) \tilde{b}_k, Y_0)=  \\
    (q_1^{-1}q_2)^u \Ibigr(\tilde{b}_k, Y_0),
  \end{array}
  \]
(the third equality uses Lemma~\ref{th:self-shift}) and, at the same time,
  \[
  \begin{array}{l}
      p(Y_u)
         = \\
      \sumop{s_1,s_2} \mathrm{rk}(\mathrm{Hom}_{\Ccat}
         (P_k, L((\tilde{t}_k)^uY_0)[s_1]\{ -s_2\})) q_1^{s_1} q_2^{s_2}
         = \\
      \sumop{s_1,s_2} \mathrm{rk}(\mathrm{Hom}_{\Ccat}
         (P_k, (\cR_k)^uL(Y_0)[s_1]\{ -s_2\})) q_1^{s_1} q_2^{s_2}
         = \\
      \sumop{s_1,s_2} \mathrm{rk}(\mathrm{Hom}_{\Ccat}
         ((\cR_k)^{-u}P_k, L(Y_0)[s_1]\{ -s_2\})) q_1^{s_1} q_2^{s_2}
         = \\
      \sumop{s_1,s_2} \mathrm{rk}(\mathrm{Hom}_{\Ccat}
         (P_k[-u]\{ -u \}, L(Y_0)[s_1]\{ -s_2\})) q_1^{s_1} q_2^{s_2}
         = \\
      \sumop{s_1,s_2} \mathrm{rk}(\mathrm{Hom}_{\Ccat}
         (P_k, L(Y_0)[s_1+u]\{ u-s_2\})) q_1^{s_1} q_2^{s_2}
         =  \\
      (q_1^{-1}q_2)^u p(Y_0),
  \end{array}
  \]
(the second equality follows from Theorem~\ref{main-theorem})
it suffices to treat the case $u=0.$

\begin{lemma}
\label{more-tables}
For $k>0$ the Poincar{\'e} polynomials of bigraded $k$-strings
with parameters $r_1=r_2=0$ and (for types $I,II,II',III,III'$) with $u=0$
are given by the following table:

\[
\begin{array}{c|c|c|c|c|c|c|c|c|c}
 I_0(0,0) & II_0(0,0) & II'_0(0,0)    & III_0(0,0) & III'_0(0,0) &
  IV & IV' & V & V' & VI(0,0) \\
 \hline
 q_1+q_2  & q_1+q_2   & q_1q_2^{-1}+1 & q_2        & 1
    & 0  & 0   & 0 & 0 & 1+ q_2
 \end{array}
\]
The Poincar{\'e} polynomials of $0$-strings with $r_1=r_2=0$ are
\[
\begin{array}{c|c|c|c|c}
 VII(0,0)        & VIII(0,0)       & IX(0,0) & X(0,0) & XI(0,0) \\
 \hline
 0               & q_1q_2^{-1} + 1 & 1       & q_1q_2^{-1} + 1  & 1
\end{array}
\]
\end{lemma}

\emph{Proof of the lemma.} Direct computation. For $0<k<m$ and a bigraded
$k$-string $\tg$ as above write down the complex $L(\tg)$ according to the
prescription in Section~\ref{complex-to-a-curve}. Then compute the hom complex
$\Hom(P_k, L(\tg))$ and the Poincar{\'e} polynomial $P(\tg).$ To simplify the
computation we can look at all homomorphisms from $P_k$ to $L(\tg),$ not just
the ones that preserve the grading, and later sort them out into graded pieces.

Example: $\tg$ has type $I_0(0,0).$ The complex $L(\tg)$ is (see also the end
of Section~\ref{complex-to-a-curve})
\begin{equation*}
  0 \lra P_{k-1}\oplus P_{k+1}\{ 1\} \lra P_k \lra 0.
\end{equation*}
The complex $\Hom(P_k, L(\tg))$ is quasi-isomorphic to
\begin{equation*}
  0 \lra \Z \{ 1\} \stackrel{0}{\lra} \Z \lra 0,
\end{equation*}
and the Poincar{\'e} polynomial $p(\tg) = q_1 + q_2.$ The remaining nine
cases are equally easy to compute, details are left to the reader.
We note that for $\tg$ of the type $IV,IV',V$ and $V'$ the complex
$\Hom(P_k, L(\tg))$ is acyclic and the Poincar{\'e} polynomial is $0.$

Cases $k=0$ and $k=m$ can be treated in the same fashion.
For $k=m,$ the four types of $m$-strings
(see figure~\ref{fig:infinite-typesb})
naturally correspond to the four out of ten types of $k$-strings with
$0<k<m, $ and the computation for the  $0< k < m$  case extends word by
word to $k=m.$

If $k=0,$ there are $5$ cases to consider. A $0$-string
has either one or two crossings, and
the Poincar{\'e} polynomial is easy to compute.
The answer is written in the second table of the lemma.
$\square$

The tables in Lemmas~\ref{th:boring} and \ref{more-tables}
are identical. Lemma~\ref{free-abelian} and
Proposition~\ref{th:dimequals} follow. $\square$

\vspace{0.10in}

Proposition~\ref{th:dimequals}, specialized to $q_1=q_2=1,$ implies
Theorem~\ref{th:gin} of the introduction. Proposition~\ref{th:dimequals},
specialized to $\tau=1,$  and
Lemma~\ref{th:detect-identity} imply that objects $\cR_{\sigma}P_j$
and $P_j$ are isomorphic if and only if $\sigma$ is the trivial
braid (a similar argument was used to prove
Proposition~\ref{th:birman-hilden}).
Therefore, the braid group acts faithfully in the category
$\Ccat$ (that proves Corollary~\ref{th:faithful}).

%
%
%
%
%
%

\section{Floer cohomology}

The definition of Floer cohomology used here is essentially the original one
\cite{floer88c}; far more general versions have been developed in the meantime,
but we do not need them. In fact, we will even impose additional restrictions
whenever this simplifies the exposition. The only exception is that we allow
Lagrangian submanifolds with boundary. After the expository part, a result is
proved which computes Floer cohomology in certain very special situations. This
is based on a simple $\Z/2$-symmetry argument, and will be the main technical
tool later on. To conclude, we discuss briefly the question of grading.

\subsection{Some basic symplectic geometry\label{subsec:symplectic}}

A symplectic manifold with contact type boundary is a compact manifold with boundary
$M$, with
a symplectic form $\o \in \Omega^2(M)$ and a contact one-form $\alpha \in \Omega^1(\partial
M)$, subject to two conditions. One is that $\o|\partial M = d\alpha$, and the other is that
the Reeb vector field $R$ on $\partial M$ should satisfy $\o(N,R)>0$ for any normal vector
field $N$ pointing outwards. A Lagrangian submanifold with Legendrian boundary in
$(M,\o,\alpha)$ is a submanifold $L \subset M$ which intersects $\partial M$ transversally,
such that $\o|L = 0$ and $\alpha|\partial L = 0$ (note that $\partial L = \emptyset$ is
allowed). For brevity, we call these objects simply {\em Lagrangian submanifolds}.

\begin{lemma} \label{th:normal-field}
Let $L \subset M$ be a Lagrangian submanifold. There is a unique $Y \in \smooth(TL|\partial
L)$ which satisfies $i_{Y}\o|T_x(\partial M) = \alpha$ for any $x \in \partial L$.
\end{lemma}

The proof is straightforward, and we omit it. By definition, a Lagrangian
isotopy is a smooth family $(L_t)_{0 \leq t \leq 1}$ of Lagrangian
submanifolds. Infinitesimally, such an isotopy is described at each time $t$ by
a section $X_t$ of the normal bundle to $L_t$, or dually by a one-form
$\beta_t$ on $L_t$. $\beta_t$ is always closed. Moreover, if $Y_t \in
\smooth(TL_t|\partial L_t)$ is the vector field from Lemma
\ref{th:normal-field} for $L = L_t$ then
\begin{equation} \label{eq:normal-relation}
\beta_t|\partial L_t = d(i_{Y_t}\beta_t).
\end{equation}
It follows that the pair $(\beta_t,i_{Y_t}\beta_t)$ defines a class in $H^1(L_t,\partial
L_t;\R)$. The Lagrangian isotopy is called {\em exact} if this class is zero for all $t$. This
means that there is a $K_t \in \smooth(L_t,\R)$ with $K_t|\partial L_t = i_{Y_t}\beta_t$ and
$dK_t = \beta_t$.

A special class of Lagrangian isotopies are those where $\partial L_t$ remains constant for
all $t$, which we call {\em isotopies rel $\partial M$}.

Another important class consists of  those isotopies where $\partial L_t$ moves along the Reeb flow
with positive speed. By this we mean that there is a function $b \in \smooth([0;1], \R^{>0})$
such that $b(t)R,$ when projected to the normal bundle of $\partial L_t,$ is equal to
$X_t|\partial L_t.$ Using
\eqref{eq:normal-relation} one can see that this is the case iff  $i_{Y_t}\beta_t
\in \smooth(\partial L_t,\R)$ is a positive constant function for each $t.$
Another equivalent condition:
if one takes an embedding $\partial L_0 \times [0;1] \longrightarrow \partial
M$ such that the image of $\partial L_0 \times \{t\}$ is $\partial L_t$, then the pullback of
$\alpha$ under it is $\psi(t)\,dt$ for some positive function $\psi$. We call isotopies of
this kind {\em positive}. There are always many positive isotopies, even exact ones, starting
at a given Lagrangian submanifold $L$. For instance, one can carry $L$ along the flow of some
Hamiltonian function $H \in \smooth(M,\R)$ with $H(x) = 1$ and $dH_x = \o(\cdot,R_x)$ at all
points $x \in \partial M$.

Let $(M,\o,\alpha)$ be a symplectic manifold with contact type boundary. One can always extend
$\alpha$ to a one-form $\theta$ on some neighbourhood $U \subset M$ of $\partial M$, such that
$d\theta = \o|U$ holds. The vector field $Z$ dual to $\theta$, $i_Z\o = \theta$, is Liouville
and points outwards along $\partial M$. The embedding $\kappa: (-r_0;0] \times \partial M
\longrightarrow M$, for some $r_0>0$, defined by the flow of $Z$ satisfies $\kappa^*\o =
d(e^r\alpha)$ and hence $\kappa^*\theta = \kappa^*(i_Z\o) = e^r\alpha$, where $r$ is the
variable in $(-r_0;0]$. Fix some $\theta$ and the corresponding $\kappa$. A Lagrangian
submanifold $L \subset M$ is called {\em $\kappa$-compatible} if $\kappa^{-1}(L) \cap
([-r_1;0] \times \partial M) = [-r_1;0] \times \partial L$ for some $0 < r_1 < r_0$.
Equivalently,  $\theta|L \in \Omega^1(L)$ should vanish near $\partial L$.

\begin{lemma} \label{th:compatible}
\begin{theoremlist} \item \label{item:compatible-one}
Any Lagrangian submanifold can be deformed, by an exact isotopy rel $\partial M$, to a
$\kappa$-compatible one.
\item \label{item:compatible-two}
Let $(L_t)_{0 \leq t \leq 1}$ be a Lagrangian isotopy, and assume that $L_0,L_1$ are
compatible with $\kappa$. Then there is another isotopy $(L_t')$ with the same endpoints, and
with $\partial L_t' = \partial L_t$ for all $t$, such that all $L_t'$ are compatible with
$\kappa$. If $(L_t)$ is exact then $(L_t')$ may be chosen to have the same property.
\end{theoremlist}
\end{lemma}

\proof \ref{item:compatible-one} Let $L \subset M$ be a Lagrangian submanifold.
Choose a $0 < r_1 < r_0$ such that $\partial L$ is a deformation retract of $K
= L \cap \kappa([-r_1;0] \times \partial M)$. $\theta|K$ is a closed one-form
which vanishes on $\partial L \subset K$; hence it can be written as the
differential of a function on $K$ which vanishes on $\partial L.$ Extend this
to a function $h \in \smooth(M,\R)$ with $h|\partial M = 0$, and consider the
family of one-forms $\theta_t = \theta - t\,dh$, $0 \leq t \leq 1$. They
satisfy $\theta_t|\partial M = \alpha$ and $d\theta_t = \omega$, and hence
define a family of embeddings $\kappa_t: [-r_2;0] \times \partial M
\longrightarrow M$ for some $r_2 > 0$. Note that $\theta_0 = \theta$ while
$\theta_1|L$ vanishes on a neighbourhood of $\partial L$, namely on $K$.

Let $U$ be a neighbourhood of $\partial M$ which is contained in the image of $\kappa_t$ for
all $t$; we can assume that $\partial M \subset U$ is a deformation retract. $X_t = (\partial
\kappa_t/\partial t) \circ \kappa_t^{-1}|U$ is a symplectic vector field on $U$ which vanishes
on $\partial M$. The corresponding one-form $\o(\cdot,X_t)$ is zero on $\partial M$, and hence
it can be written as the boundary of a function $H_t \in \smooth(U,\R)$ with $H_t|\partial M =
0$. After possibly making $U$ smaller one finds that for each $t$, $X_t$ can be extended to a
Hamiltonian vector field $\widetilde{X}_t$ on the whole of $M$. Let $\phi_t$ be the flow of
$(\widetilde{X}_t)$ seen as a time-dependent vector field. By definition $\phi_t \circ
\kappa_0 = \kappa_t$ on some neighbourhood of the boundary, and, therefore,
\[
\phi_1^*\theta = (\kappa_1\kappa_0^{-1})^*\theta = (\kappa_0^{-1})^*(e^r\alpha) = \theta_1
\quad\text{near }\partial M.
\]
This shows that $\phi_1^*(\theta|\phi_1(L))$ is equal to $\theta_1|L$ near $\partial L$, so
that $\phi_1(L)$ is $\kappa$-compatible. One sees easily that the isotopy $L_t = \phi_t(L)$ is
exact, which finishes the argument. The proof of \ref{item:compatible-two} is just a
parametrized version of \ref{item:compatible-one}. \qed

Assume now that the relative symplectic class $[\o,\alpha] \in H^2(M,\partial M;\R)$ is zero.
Then there is a one-form $\theta$ on $M$ with $\theta|\partial M = \alpha$ and $d\theta = \o$.
Fix such a $\theta$; the corresponding embedding $\kappa$ is defined on all of $\R^{\leq 0}
\times \partial M$. A Lagrangian submanifold $L$ is called {\em $\theta$-exact} if $[\theta|L]
\in H^1(L,\partial L;\R)$ is zero. The next Lemma describes the relationship between the
absolute notion of a $\theta$-exact submanifold and the relative notion of exact isotopy. The
Lemma after that concerns the existence of $\theta$-exact submanifolds in a given Lagrangian
isotopy class. We omit the proofs.

\begin{lemma} \label{th:exact-isotopies}
Let $(L_t)_{0 \leq t \leq 1}$ be a Lagrangian isotopy such that $L_0$ is $\theta$-exact. Then
the isotopy is exact $\Leftrightarrow$ all the $L_t$ are $\theta$-exact.
\end{lemma}

\begin{lemma} \label{th:exist-exact}
\begin{theoremlist}
\item
Let $L \subset M$ be a Lagrangian submanifold such that $H^1(M,\partial M;\R)$ surjects onto
$H^1(L,\partial L;\R)$. Then there is a Lagrangian isotopy $(L_t)$ rel $\partial M$ with $L_0
= L$ and such that $L_1$ is $\theta$-exact.
\item \label{item:deform-to-exact-isotopy}
Let $L_0,L_1 \subset M$ be two $\theta$-exact Lagrangian submanifolds which are isotopic rel
$\partial M$. Assume that $H^1(M,\partial M;\R)$ surjects onto $H^1(L_j,\partial L_j;\R)$.
Then there is an exact isotopy rel $\partial M$ which connects $L_0$ and $L_1$.
\end{theoremlist}
\end{lemma}

There is a similar notion for symplectic automorphisms. Let $\Symp(M,\partial
M,\o)$ be the group of symplectic automorphisms $\phi: M \longrightarrow M$
which are the identity in some neighbourhood (depending on $\phi$) of $\partial
M$. By definition, a Hamiltonian isotopy in this group is one generated by a
time-dependent Hamiltonian function which vanishes near $\partial M$. Call
$\phi \in \Symp(M,\partial M,\o)$ {\em $\theta$-exact} if
$[\phi^*\theta-\theta] \in H^1(M,\partial M;\R)$ is zero. The analogue of Lemma
\ref{th:exact-isotopies} is this: let $(\phi_t)$ be an isotopy in
$\Symp(M,\partial M,\o)$ with $\phi_0$ $\theta$-exact. Then $(\phi_t)$ is a
Hamiltonian isotopy iff all the $\phi_t$ are $\theta$-exact. As for Lemma
\ref{th:exist-exact}, the analogous statement is actually simpler: the subgroup
of $\theta$-exact automorphisms is always a deformation retract of
$\Symp(M,\partial M,\o)$. Finally, note that $\theta$-exact automorphisms
preserve the class of $\theta$-exact Lagrangian submanifolds.

%
%
\subsection{The definition\label{subsec:floer}}

Let $(M,\o,\alpha)$ be a symplectic manifold with contact type boundary, such that
$[\o,\alpha] \in H^2(M,\partial M;\R)$ is zero. Fix a $\theta \in \Omega^1(M)$ with
$\theta|\partial M = \alpha$, $d\theta = \o$. Floer cohomology associates a finite-dimensional
vector space over $\Z/2$, denoted by $HF(L_0,L_1)$, to any pair $(L_0,L_1)$ of $\theta$-exact
Lagrangian submanifolds whose boundaries are either the same or else disjoint. We now give a
brief account of its definition, concentrating on the issues which arise in connection with
the boundary $\partial M$. For expositions emphasizing other aspects see e.g.\ \cite{oh93},
\cite{salamon90}, \cite{salamon-zehnder92}.

We begin with a digression concerning almost complex structures. Let $\xi \longrightarrow
\partial M$ be the symplectic vector bundle $\xi = (\ker\,\alpha, d\alpha)$. Given a
compatible almost complex structure $\bar{j}$ on $\xi$, one can define an almost complex
structure $\overline{J}$ on $\R \times \partial M$, compatible with $d(e^r\alpha)$, by
$\overline{J}(y_1,X + y_2 R) = (-y_2,\bar{j}X + y_1 R)$ for $X \in \xi$ and $y_1,y_2 \in \R$.
The almost complex structures $\overline{J}$ appeared first in \cite[p.\ 529]{hofer93}. The
next Lemma describes their main property; it is taken from \cite[Lemma
3.9.2]{eliashberg-hofer-salamon95} with some modifications.

\begin{lemma} \label{th:convexity}
Let $\Sigma$ be an open subset of $\R \times [0;1]$, and $(\bar{j}_{s,t})$ a smooth family of
compatible almost complex structures on $\xi$ parametrized by $(s,t) \in \Sigma$. Let
$(\overline{J}_{s,t})$ be the corresponding family of almost complex structures on $\R \times
\partial M$. Let $\Lambda_0,\Lambda_1 \subset \partial M$ be two Legendrian submanifolds, and
$u = (u_1,u_2): \Sigma \longrightarrow \R \times \partial M$ a smooth map which satisfies
$\partial u/\partial s + \overline{J}_{s,t}(u)\partial u/\partial t = 0$, such that $u(s,t)
\in \R \times \Lambda_t$ for all $(s,t) \in (\R \times \{0;1\}) \cap \Sigma$. Then $u_1$ has
no local maxima.
\end{lemma}

\proof Let $pr_1,pr_2$ be the projections from $\R \times \partial M$ to the two factors. For
any $(s,t) \in \Sigma$ one has $D pr_1 = \alpha \circ D pr_2 \circ \overline{J}_{s,t}$. Hence
$\partial u_1/\partial s = \alpha(D pr_2 \circ \overline{J}_{s,t} \circ \partial u/\partial s)
= \alpha(\partial u_2/\partial t)$, and similarly $\partial u_1/\partial t = -\alpha(\partial
u_2/\partial s)$. Therefore $\Delta u_1 = -d\alpha(\partial u_2/\partial s,\partial
u_2/\partial t) \leq 0$. By assumption $\partial u_2/\partial s \in T\Lambda_t$ at each point
$(s,t) \in (\R \times \{0;1\}) \cap \Sigma$, and hence $\partial u_1/\partial t = -
\alpha(\partial u_2/\partial s) = 0$. This shows that $u_1$ is subharmonic and satisfies von
Neumann boundary conditions on $(\R \times \{0;1\}) \cap \Sigma$. The maximum principle
implies that such a function cannot have local maxima. \qed

Let $\kappa: \R^{\leq 0} \times \partial M \longrightarrow M$ be the embedding associated to
$\theta$, as in the previous section. We define $\JJ$ to be the space of all smooth families
$\J = (J_t)_{0 \leq t \leq 1}$ of $\o$-compatible almost complex structures on $M$ which have
the following property: there is an $r_0>0$ and a family $(\bar{j}_t)$ of compatible almost
complex structures on $\xi$ such that for each $t$, $\kappa^*J_t|[-r_0;0] \times \partial M$
agrees with $\overline{J}_t$.

Now take two $\theta$-exact Lagrangian submanifolds $L_0,L_1 \subset M$ with $\partial L_0
\cap \partial L_1 = \emptyset$. The first step in defining $HF(L_0,L_1)$ is to deform them in
a suitable way.

\begin{lemma} \label{th:perturb}
There are Lagrangian submanifolds $L'_0,L'_1$ such that $L_j$ and $L'_j$ are joined by an
exact isotopy rel $\partial M$, for $j = 0,1$, and with the following properties: each $L'_j$
is $\kappa$-compatible, and the intersection $L'_0 \cap L'_1$ is transverse.
\end{lemma}

This follows from Lemma \ref{th:compatible}\ref{item:compatible-one} and a well-known
transversality argument. Note that $L_0',L_1'$ are again $\theta$-exact, because of Lemma
\ref{th:exact-isotopies}. For $\J = (J_t) \in \JJ$, let $\moduli(\J)$ be the set of smooth
maps $u: \R \times [0;1] \longrightarrow M$ satisfying Floer's equation
\begin{equation} \label{eq:floer}
\partial u/\partial s + J_t(u)\partial u/\partial t = 0, \quad u(\R \times \{0\}) \subset
L'_0,\;\; u(\R \times \{1\}) \subset L'_1,
\end{equation}
and whose energy $E(u) = \int u^*\o$ is finite. We remind the reader that \eqref{eq:floer} can
be seen as the negative gradient flow equation of the action functional $a$ on the path space
${\mathcal P}(L'_0,L'_1) = \{\gamma \in \smooth([0;1], M) \suchthat \gamma(j) \in L_j',\; j =
0,1\}$. To define $a$ one needs to choose functions $H_j \in \smooth(L_j',\R)$ with $dH_j =
\theta|L_j'$ (such functions exist by the assumption of $\theta$-exactness); then $a(\gamma) =
-\int_{[0;1]} \gamma^*\theta + H_1(\gamma(1)) - H_0(\gamma(0))$.

\begin{lemma} \label{th:bounded}
For any $\J \in \JJ$, there is a compact subset of $M \setminus \partial M$ which contains the
image of all $u \in \moduli(\J)$.
\end{lemma}

This is mainly a consequence of Lemma \ref{th:convexity} and of the definition of $\JJ$. The
other elements which enter the proof are: the convergence of each $u \in \moduli(\J)$ towards
limits $x_{\pm} = \lim_{s \rightarrow \pm\infty} u(s,\cdot)$ in $L_0' \cap L_1'$; the
$\kappa$-compatibility of $L_0',L_1'$; and the fact that $\partial L_0' \cap \partial L_1' =
\partial L_0 \cap \partial L_1 = \emptyset$. We leave the details to the reader. The
importance of Lemma \ref{th:bounded} is that it allows one to ignore the boundary of $M$ in
all arguments concerning solutions of \eqref{eq:floer}.

One can associate to any $u \in \moduli(\J)$ a Fredholm operator $D_{u,\J}: \W^1_u
\longrightarrow \W^0_u$ from $\W^1_u = \{X \in W^{1,p}(u^*TM) \suchthat X(\cdot,0) \in
u^*TL'_0, \; X(\cdot,1) \in u^*TL'_1\}$ to $\W^0_u = L^p(u^*TM)$ ($p>2$ is some fixed real
number) which describes the linearization of \eqref{eq:floer} near $u$. Let
$\moduli_k(x_-,x_+;\J) \subset \moduli(\J)$ be the subset of those $u$ with limits $x_-,x_+
\in L_0' \cap L_1'$ and such that $D_{u,\J}$ has index $k \in \Z$. There is a natural action
of $\R$ on each of these sets, by translation in the $s$-variable. Call $u \in \moduli(\J)$
regular if $D_{u,\J}$ is onto. If this is true for all $u \in \moduli(\J)$ then $\J$ itself is
called regular.

\begin{prop}[Floer-Hofer-Salamon \cite{floer-hofer-salamon94}, Oh \cite{oh96c}] \label{th:fhs}
The regular $\J$ form a $\smooth$-dense subset $\JJ^\reg(L'_0,L'_1) \subset \JJ$.
\end{prop}

\begin{prop}[Floer] \label{th:floer}
For $\J \in \JJ^\reg(L_0',L_1')$, the quotients $\moduli_1(x_-,x_+;\J)/\R$ are finite sets.
Moreover, if we define $\nu(x_-,x_+;\J) \in \Z/2$ to be the number of points mod $2$ in
$\moduli_1(x_-,x_+;\J)/\R$, then $\sum_x \nu(x_-,x;\J) \nu(x,x_+;\J) = 0$ for each $x_-,x_+$.
\end{prop}

Let $CF(L_0',L_1')$ be the $\Z/2$-vector space freely generated by the points
of $L_0' \cap L_1'$, and $d_\J$ the endomorphism of this space which sends the
basis element $\gen{x_+}$ to $\sum_x \nu(x,x_+;\J) \gen{x}$. Theorem
\ref{th:floer} implies that $d_{\J}^2 = 0$, and one sets $HF(L'_0,L'_1;\J) =
\ker\,d_{\J}/\im\,d_{\J}$. Finally, $HF(L_0,L_1)$ is defined as
$HF(L'_0,L'_1;\J)$ for any choice of $L_0',L_1'$ as in Lemma \ref{th:perturb}
and any $\J \in \JJ^\reg(L_0',L_1')$. One can show that this is independent of
the choices up to canonical isomorphism. The proof uses the continuation maps
introduced in \cite[Section 6]{salamon-zehnder92}. We will not explain the
details, but it seems appropriate to mention two simple facts which are used in
setting up the argument: firstly, an exact Lagrangian isotopy rel $\partial M$
which consists only of $\kappa$-compatible submanifolds can be embedded into a
Hamiltonian isotopy in $\Symp(M,\partial M,\o)$. Secondly, there is an analogue
of Lemma \ref{th:bounded} for solutions of the continuation equation, which is
again derived from Lemma \ref{th:convexity}. As well as proving its
well-definedness, the continuation argument also establishes the basic isotopy
invariance property of Floer cohomology:

\begin{prop} \label{th:isotopy-invariance}
$HF(L_0,L_1)$ is invariant (up to isomorphism) under exact isotopies of $L_0$ or $L_1$ rel
$\partial M$.
\end{prop}

We will now extend the definition of Floer cohomology to pairs of Lagrangian submanifolds
whose boundaries coincide. The basic idea is the same as in the case of geometric intersection
numbers. Let $(L_0,L_1)$ be two $\theta$-exact Lagrangian submanifolds with $\partial L_0 =
\partial L_1$, and let $L_0^+$ be a submanifold obtained from $L_0$ by a sufficiently small
Lagrangian isotopy which is exact and positive. Then $\partial L_0^+ \cap \partial L_1 =
\emptyset$, and one sets $HF(L_0,L_1) = HF(L_0^+,L_1)$. The proof that this is independent of
the choice of $L_0^+$ requires a generalization of Proposition \ref{th:isotopy-invariance}:

\begin{lemma} \label{th:moving-boundary}
Let $(L_{0,s})_{0 \leq s \leq 1}$ be an isotopy of $\theta$-exact Lagrangian submanifolds, and
$L_1$ a $\theta$-exact Lagrangian submanifold such that $\partial L_{0,s} \cap \partial L_1 =
\emptyset$ for all $s$. Then $HF(L_{0,s},L_1)$ is independent of $s$ up to isomorphism.
\end{lemma}

\proof[Sketch of proof] It is sufficient to show that the dimension of $HF(L_{0,s},L_1)$ is
locally constant in $s$ near some fixed $\hat{s} \in [0;1]$. Let $L_1'$ be a Lagrangian
submanifold which is $\kappa$-compatible and isotopic to $L_1$ by an exact isotopy rel
$\partial M$. Let $(L'_{0,s})_{0 \leq s \leq 1}$ be an isotopy of Lagrangian submanifolds,
such that each $L'_{0,s}$ is $\kappa$-compatible and isotopic to $L_{0,s}$ by an exact isotopy
rel $\partial M$. We may assume that $L_{0,\hat{s}}'$ intersects $L_1'$ transversally. Take a
$\J \in \JJ^\reg(L'_{0,\hat{s}},L'_1)$. For all $s$ sufficiently close to $\hat{s}$, the
intersection $L_{0,s}' \cap L_1'$ remains transverse, and $\J$ lies still in
$\JJ^\reg(L'_{0,s},L'_1)$. Hence by definition of Floer cohomology
\[
HF(L_{0,s},L_1) = HF(L'_{0,s},L'_1;\J)
\]
for all $s$ close to $\hat{s}$. By considering the corresponding parametrized
moduli spaces one sees that the complex $(CF(L'_{0,s},L'_1),d_\J)$ is
independent of $s$ (near $\hat{s})$ up to isomorphism; compare \cite[Lemma
3.4]{floer88c}. \qed

The same argument can be used to extend Proposition~\ref{th:isotopy-invariance}
to pairs $(L_0,L_1)$ with $\partial L_0 = \partial L_1$. It is difficult to see
whether the isomorphisms obtained in this way are canonical or not. This is an
inherent weakness of our approach; there are alternatives which do not appear
to have this problem \cite{oh-new}. In any case, the question is irrelevant for
the purposes of the present paper.

A theorem of Floer \cite{floer89} says that $HF(L,L) \iso H^*(L;\Z/2)$ for any
$\theta$-exact Lagrangian submanifold $L$ (Floer stated the theorem for
$\partial L = \emptyset$ only, but his proof extends to the general case).
There is no general way of computing $HF(L_0,L_1)$ when $L_0 \neq L_1$, but
there are a number of ad hoc arguments which work in special situations. One
such argument will be the subject of the next two sections.

%
%
\subsection{Equivariant transversality}

Let $(M,\o,\alpha)$ and $\theta$ be as before. In addition we now assume that $M$ carries an
involution $\iota$ which preserves $\o$, $\theta$, and $\alpha$. The fixed point set $M^\iota$
is again a symplectic manifold with contact type boundary. Moreover, if $L \subset M$ is a
Lagrangian submanifold with $\iota(L) = L$, its fixed part $L^\iota = L \cap M^\iota$ is a
Lagrangian submanifold of $M^\iota$. We denote by $\JJ^\iota \subset \JJ$ the subspace of
those $\J = (J_t)_{0 \leq t \leq 1}$ such that $\iota^*J_t = J_t$ for all $t$.

Let $L_0',L_1'$ be a pair of $\theta$-exact and $\kappa$-compatible Lagrangian submanifolds of
$M$, which intersect transversally and satisfy $\partial L_0' \cap \partial L_1' = \emptyset$.
Assume that $\iota(L_j') = L_j'$ for $j = 0,1$. If one wants to use this symmetry property to
compute Floer cohomology, the following question arises: can one find a $\J \in \JJ^\iota$
which is regular? The answer is no in general, because of a simple phenomenon which we will
now describe. Take a $\J \in \JJ^\iota$, and let $u \in \moduli(\J)$ be a map whose image lies
in $M^\iota$. Then the spaces $\W^1_u$ and $\W^0_u$ carry natural $\Z/2$-actions, and
$D_{u,\J}$ is an equivariant operator. Assume that the invariant part
\[
D^{\Z/2}_{u,\J}: (\W^1_u)^{\Z/2} \longrightarrow (\W^0_u)^{\Z/2}
\]
is surjective, and that $\mathrm{index} \, D_{u,\J}^{\Z/2} > \mathrm{index} \, D_{u,\J}$. Then
$D_{u,\J}$ can obviously not be onto. This situation is stable under perturbation, which means
that for any $\J' \in \JJ^\iota$ sufficiently close to $\J$ one can find a $u' \in
\moduli(\J)$ with the same property as $u$. Hence $\JJ^\iota \cap \JJ^\reg(L_0',L'_1)$ is not
dense in $\JJ^\iota$. This is a general problem in equivariant transversality theory, and it
is not difficult to construct concrete examples. By taking the reasoning a little further, one
can find cases where $\JJ^\iota \cap \JJ^\reg(L'_0,L'_1) = \emptyset$. The aim of this section
is to show that the maps $u \in \moduli(\J)$ with $\im(u) \subset M^\iota$ are indeed the only
obstruction; in other words, equivariant transversality can be achieved everywhere except at
these solutions.

\begin{lemma} \label{th:regular}
Let $u \in \moduli(\J)$, for some $\J \in \JJ^\iota$, be a map which is not constant and such
that $\im(u) \not\subset M^\iota$. Let $x_-,x_+ \in L_0' \cap L_1'$ be the limits of $u$. Then
the subset $R(u) \subset \R \times (0;1)$ of those $(s,t)$ which satisfy
\[
du(s,t) \neq 0, \quad u(s,t) \notin u(\R \setminus \{s\},t) \cup \iota(u(\R,t)) \cup \{x_\pm\}
\cup \{\iota(x_\pm)\}
\]
is open and dense.
\end{lemma}

\proof We will use the methods developed in \cite{floer-hofer-salamon94}. To
begin with, take a more general situation: let $M$ be an arbitrary manifold and
$(J_t)_{0 < t < 1}$ a smooth family of almost complex structures on it. We
consider smooth maps $v: \Sigma \longrightarrow M$, where $\Sigma$ is some
connected open subset of $\R \times (0;1)$, satisfying
\begin{equation} \label{eq:j-holo}
\partial v/\partial s + J_t(v) \partial v/\partial t = 0.
\end{equation}
Some basic properties of this equation are:
\begin{Jlist}
\item \label{item:critical}
If $v$ is a non-constant solution of \eqref{eq:j-holo} then $C(v) = \{(s,t) \in \Sigma
\suchthat dv(s,t) = 0\}$ is a discrete subset of $\Sigma$.
\item \label{item:coincidence}
Let $v_1,v_2$ be two solutions of \eqref{eq:j-holo} defined on the same set $\Sigma$, such
that $v_1 \not\equiv v_2$. Then $A(v_1,v_2) = \{ (s,t) \in \Sigma \suchthat v_1(s,t) =
v_2(s,t)\}$ is a discrete subset.
\item \label{item:preimage}
Let $v$ be a non-constant solution of \eqref{eq:j-holo}. For each $x \in M$, $v^{-1}(x)
\subset \Sigma$ is a discrete subset.
\item \label{item:shifted}
Let $v_1: \Sigma_1 \longrightarrow M$ and $v_2: \Sigma_2 \longrightarrow M$ be two solutions
of \eqref{eq:j-holo} such that $v_2$ is an embedding. Assume that for any $(s,t) \in \Sigma_1$
there is an $(s',t') \in \Sigma_2$ with $t = t'$ and $v_1(s,t) = v_2(s',t')$. Then there is a
$\sigma \in \R$ such that for all $(s,t) \in \Sigma_1$ we have $(s-\sigma,t) \in \Sigma_2$ and
$v_1(s,t) = v_2(s-\sigma,t)$.
\item \label{item:better-coincidence}
Let $v_1,v_2: \R \times (0;1) \longrightarrow M$ be two non-constant solutions of
\eqref{eq:j-holo}. Assume that $v_2$ is not a translate of $v_1$ in $s$-direction. Then for
any $\rho>0$ the subset $S_\rho(v_1,v_2) = \{ (s,t) \in \R \times (0;1) \suchthat v_1(s,t)
\notin v_2([-\rho;\rho] \times \{t\})\}$ is open and dense in $\R \times (0;1)$.
\end{Jlist}
\ref{item:critical} is \fhs{Corollary 2.3(ii)}. \ref{item:coincidence} is a form of the unique
continuation theorem \fhs{Proposition 3.1}, and \ref{item:preimage} follows from it by letting
$v_2 \equiv x$ be a constant map. \ref{item:shifted} is a variant of \fhs{Lemma 4.2}, and is
proved as follows. $\psi = v_2^{-1}v_1: \Sigma_1 \longrightarrow \Sigma_2$ is a holomorphic
map of the form $\psi(s,t) = (\psi_1(s,t),t)$. By looking at the derivative of $\psi$ one sees
that it must be a translation in $s$-direction.

The proof of the final property \ref{item:better-coincidence} is similar to
that of \fhs{Theorem 4.3}; however, the differences are sufficiently important
to deserve a detailed discussion. It is clear that $S_\rho(v_1,v_2)$ is open
for any $\rho>0$. Now assume that $(S,T) \in \R \times (0;1)$ is a point which
has a neighbourhood disjoint from $S_\rho(v_1,v_2)$. After possibly moving
$(S,T)$ slightly, one can assume that $v_1(S,T) \notin v_2(C(v_2))$. This can
be achieved because $C(v_2)$ is a countable set by \ref{item:critical}, and
hence $v_1^{-1}(v_2(C(v_2)))$ is countable by \ref{item:preimage}. Again
because of \ref{item:preimage}, there are only finitely many $S_1',\dots,S_k'
\in [-\rho;\rho]$ such that $v_1(S,T) = v_2(S_j',T)$. By choice of $(S,T)$, we
have $dv_2(S_j',T) \neq 0$ for all these $S_j'$. Choose an $\epsilon>0$ such
that the restriction of $v_2$ to the open disc $D_\epsilon(S_j',T)$ of radius
$\epsilon$ around $(S_j',T)$ is an embedding for any $j$.

{\bf Claim A:} {\em There is an $\epsilon' > 0$ with the following property: for each $(s,t)
\in D_{\epsilon'}(S,T)$ there is an $(s',t') \in \bigcup_{j=1}^k D_{\epsilon/2}(S_j',T)$ with
$t' = t$ and $v_1(s,t) = v_2(s',t)$.}

To prove this, take a sequence of points $(s_i,t_i)_{i \in \N}$ with limit $(S,T)$. For
sufficiently large $i$ we have $(s_i,t_i) \notin S_\rho(v_1,v_2)$. Hence there are $s'_i \in
[-\rho;\rho]$, $s_i' \neq s_i$, such that $v_2(s_i',t_i) = v_1(s_i,t_i)$. If the claim is
false, one can choose the $(s_i,t_i)$ and $s_i'$ in such a way that $(s_i',t_i) \notin
\bigcup_j D_{\epsilon/2}(S_j',T)$. Then, after passing to a subsequence, one obtains a limit
$s' \in [-\rho;\rho]$ with $v_2(s',T) = v_1(S,T)$ and $s' \neq S_1',\dots,S_k'$, which is a
contradiction.

{\bf Claim B:} {\em There is a nonempty connected open subset $\Sigma_1 \subset
D_{\epsilon'}(S,T)$ and a $j \in \{1,\dots,k\}$ with the following property: for any $(s,t)
\in \Sigma_1$ there is an $(s',t') \in D_\epsilon(S_j',T)$ with $t' = t$ and $v_1(s,t) =
v_2(s',t)$.}

Let $U_j \subset D_{\epsilon'}(S,T)$ be the closed subset of those $(s,t)$ such that $v_1(s,t)
= v_2(s',t)$ for some $(s',t) \in \overline{D_{\epsilon/2}(S_j',T)}$. Claim A says that $U_1
\cup \dots \cup U_k = D_{\epsilon'}(S,T)$; hence at least one of the $U_j$ must have nonempty
interior. Define $\Sigma_1$ to be some small open disc inside that $U_j$; this proves claim B.

Now consider the restrictions $v_1|\Sigma_1$ and $v_2|\Sigma_2$, where
$\Sigma_2 = D_{\epsilon}(S_j',T)$. By definition $v_2|\Sigma_2$ is an
embedding; applying \ref{item:shifted} shows that $v_1|\Sigma_1$ is a translate
of $v_2|\Sigma_2$ in $s$-direction. Because of unique continuation
\ref{item:coincidence} it follows that $v_1$ must be a translate of $v_2$,
which completes the proof of \ref{item:better-coincidence}.

After these preliminary considerations, we can now turn to the actual proof of Lemma
\ref{th:regular}. Set $v = \iota \circ u$. $v$ cannot be a translate of $u$ in $s$-direction.
To see this, assume that on the contrary $v(s,t) = u(s-\rho,t)$ for some $\rho \in \R$. The
case $\rho = 0$ is excluded by the assumptions, since it would imply that $\im(u) \subset
M^\iota$. If $\rho \neq 0$ then $u(s,t) = u(s-2\rho,t)$ which, because of the finiteness of
the energy, means that $u$ is constant; this is again excluded by the assumptions.

Clearly $R(u) = R_1(u) \cap R_2(u) \cap R_3(u)$ where
\begin{align*}
 R_1(u) &= \{(s,t) \in \R \times (0;1) \suchthat du(s,t) \neq 0,\;
 u(s,t) \notin u(\R \setminus \{s\},t), \; u(s,t) \neq x_{\pm}\}, \\
 R_2(u) &= \{(s,t) \in \R \times (0;1) \suchthat u(s,t) \neq \iota(x_\pm)\}, \\
 R_3(u) &= \{(s,t) \in \R \times (0;1) \suchthat u(s,t) \notin \iota(u(\R,t))\}.
\end{align*}
$R_1(u)$ is the set of regular points of $u$ as defined in \cite{floer-hofer-salamon94}, and
it is open and dense by a slight variation of \fhs{Theorem 4.3}; see also \cite{oh96c}.
$R_2(u)$ is open and dense by \ref{item:preimage} above. $R_3(u)$ is the intersection of
countably many sets $S_\rho(u,v)$, each of which is open and dense by
\ref{item:better-coincidence}. Baire's theorem now shows that $R(u)$ is dense. As for its
openness, it can be proved by an elementary argument, as in \fhs{Theorem 4.3}. \qed

\begin{prop} \label{th:transversality}
There is a $\J \in \JJ^\iota$ such that every $u \in \moduli(\J)$ whose image is not contained
in $M^\iota$ is regular. In fact, the subspace of such $\J$ is $\smooth$-dense in $\JJ^\iota$.
\end{prop}

\proof[Sketch of proof] It is convenient to recall first the proof of the basic transversality
result in Floer theory, Proposition \ref{th:fhs}. Let $\T_\J$ be the tangent space of $\JJ$ at
some point $\J$. This space consists of smooth families $\Y = (Y_t)_{0 \leq t \leq 1}$ of
sections $Y_t \in \smooth(\End(TM))$ satisfying some additional properties, which we will not
write down here. For $\J \in \JJ$ and $u \in \moduli(\J)$ one defines an operator
\[
\widetilde{D}_{u,\J}: \W^1_u \times \T_\J \longrightarrow \W^0_u, \quad
\widetilde{D}_{u,\J}(X,\Y) = D_{u,\J}(X) + Y_t(u)\partial u/\partial t.
\]
The main point in the proof of Proposition \ref{th:fhs} is to show that
$\widetilde{D}_{u,\J}$ is always onto. Assume that the contrary holds: then,
because the image is closed, there is a nonzero $Z \in L^q(u^*TM)$, with $1/p +
1/q = 1$, which satisfies a linear $\bar\partial$-type equation, hence is
smooth on $\R \times (0;1)$, and such that
\begin{equation} \label{eq:orthogonality}
 \int_{\R \times [0;1]} \o(Y_t(u) \partial u/\partial t, J_t(u) Z) \; ds\, dt = 0
 \quad \text{for all }\Y \in \T_\J.
\end{equation}
Let $x_\pm$ be the limits of $u$. Excluding the case of constant maps $u$, for which
$D_{u,\J}$ is already onto, one can find an $(s,t) \in \R \times (0;1)$ with $du(s,t) \neq 0$,
$u(s,t) \neq x_{\pm}$, and $u(s,t) \notin u(\R \setminus \{s\},t)$. In order to avoid problems
with the behaviour of $\Y$ near the boundary, it is useful to require that $u(s,t)$ lies
outside a fixed small neighbourhood of $\partial M$. Then, assuming that $Z(s,t) \neq 0$, one
can construct a $\Y \in \T_\J$ which is concentrated on a small neighbourhood of $u(s,t)$ and
which contradicts \eqref{eq:orthogonality}. Since the set of points $(s,t)$ with these
properties is open, it follows that $Z$ vanishes on some open subset. By the unique
continuation for solutions of linear $\bar\partial$-equations, $Z \equiv 0$, which is a
contradiction.

The proof of Proposition \ref{th:transversality} differs from this in that one considers only
equivariant $\J$. This means that the domain of $\widetilde{D}_{u,\J}$ must be restricted to
those $\Y$ which are equivariant. One considers a supposed $Z$ as before, and takes a point
$(s,t) \in R(u)$ as in Lemma \ref{th:regular}. The same local construction as before yields a
non-equivariant $\Y \in \T_\J$ which is concentrated near $u(s,t)$. One then makes this
equivariant by taking $Y_t + \iota^*Y_t$; the properties of $R(u)$ ensure that this averaged
element of the tangent space still leads to a contradiction with \eqref{eq:orthogonality}. The
details are straightforward, and we leave them to the reader. \qed

%
%
\subsection{A symmetry argument}

The notations $(M,\o,\alpha)$, $\theta$, $\kappa$ and $\iota$ are as in the previous section.
We begin with the most straightforward application of equivariant transversality.

\begin{lemma} \label{th:first-symmetry}
Let $L_0,L_1 \subset M$ be two $\theta$-exact Lagrangian submanifolds with $\partial L_0 \cap
\partial L_1 = \emptyset$ and $\iota(L_j) = L_j$. Assume that the following properties hold:
\begin{Slist}
\item \label{item:first-intersection}
the intersection $L_0 \cap L_1$ is transverse, and each intersection point lies in $M^\iota$;
\item \label{item:first-no-invariants}
there is no continuous map $v: [0;1]^2 \longrightarrow M^\iota$ with the following properties:
$v(0,t) = x_-$, $v(1,t) = x_+$ for all $t$, where $x_- \neq x_+$ are points of $L_0 \cap L_1$,
and $v(s,0) \in L_0$, $v(s,1) \in L_1$ for all $s$.
\end{Slist}
Then $\dim_{\Z/2} HF(L_0,L_1) = |L_0 \cap L_1|$.
\end{lemma}

\proof An inspection of the proof of Lemma \ref{th:compatible}\ref{item:compatible-one} shows
that one can find $\kappa$-compatible Lagrangian submanifolds $L_j'$, $j = 0,1$, which are
isotopic to $L_j$ by an exact {\em and $\iota$-equivariant} Lagrangian isotopy rel $\partial
M$. Moreover, this can be done in such a way that $L_j'$ is equal to $L_j$ outside a small
neighbourhood of $\partial M$, and such that $L_0' \cap L_1' = L_0 \cap L_1$. Then the
conditions \ref{item:first-intersection} and \ref{item:first-no-invariants} continue to hold
with $(L_0,L_1)$ replaced by $(L_0',L_1')$. Choose a $\J \in \JJ^\iota$ which has the property
stated in Proposition \ref{th:transversality}. We will show that $\J \in \JJ^\reg(L_0',L_1')$.
Assume that on the contrary there is a $u \in \moduli(\J)$ which is not regular. Note that $u$
cannot be constant, since the constant maps are regular for any choice of $\J$. By Proposition
\ref{th:transversality} $\im(u) \subset M^\iota$. Then condition
\ref{item:first-no-invariants} says that the two endpoints $\lim_{s \rightarrow \pm \infty}
u(s,\cdot)$ must necessarily agree. Because of the gradient flow interpretation of Floer's
equation (recall that the action functional is not multivalued in the situation which we are
considering), $u$ must be constant, which is a contradiction. This shows that $\J$ is indeed
regular.

Composition with $\iota$ defines an involution on $\moduli(\J)$. It is not difficult to see
that the index of $D_{u,\J}$ and $D_{\iota \circ u,\J}$ are the same; indeed, there are
obvious isomorphisms between the kernels resp.\ cokernels of these two operators. Moreover,
$u$ and $\iota \circ u$ have the same endpoints because $L_0' \cap L_1' \subset M^\iota$.
Hence the involution preserves each subset $\moduli_k(x_-,x_+;\J) \subset \moduli(\J)$.
Consider the induced involution $\bar{\iota}$ on the quotients $\moduli_1(x_-,x_+;\J)/\R$. A
fixed point of $\bar{\iota}$ would be a map $u \in \moduli(\J)$ such that $u(s,t) =
\iota(u(s-\sigma,t))$ for some $\sigma \in \R$. By the same argument as in the proof of Lemma
\ref{th:regular}, the only maps with this property are the constant ones, which do not lie in
$\moduli_1(x_-,x_+)$ since the associated operators have index zero. Hence $\bar{\iota}$ is
free for all $x_-,x_+$.

By definition, $HF(L_0,L_1)$ is the homology of $(CF(L_0',L_1'),d_{\J})$. Since each set
$\moduli_1(x_-,x_+;\J)/\R$ admits a free involution, $\nu(x_-,x_+;\J)$ and hence $d_{\J}$ are
zero, so that $\dim HF(L_0,L_1) = \dim CF(L_0',L_1') = |L_0' \cap L_1'| = |L_0 \cap L_1|$.
\qed

\begin{prop} \label{th:symmetry}
Let $L_0,L_1 \subset M$ be two $\theta$-exact Lagrangian submanifolds with $\partial L_0 \cap
\partial L_1 = \emptyset$ and $\iota(L_j) = L_j$. Assume that the following properties hold:
\begin{Sprimelist}
\item \label{item:clean}
the intersection $N = L_0 \cap L_1$ is clean, and there is a $\iota$-invariant Morse function
$h$ on $N$ whose critical points are precisely the points of $N^\iota = N \cap M^\iota$;
\item \label{item:no-invariants}
same as \ref{item:first-no-invariants} above.
\end{Sprimelist}
Then $\dim_{\Z/2} HF(L_0,L_1) = |N^\iota| = \dim_{\Z/2} H^*(N;\Z/2)$.
\end{prop}

\proof Clean intersection means that $N$ is a smooth manifold and $TN = (TL_0|N) \cap
(TL_1|N)$. We will now describe a local model for clean intersections, due to Weinstein
\cite[Theorem 4.3]{weinstein73}. There are a neighbourhood $U \subset L_0$ of $N$, a function
$f \in \smooth(U,\R)$ with $df^{-1}(0) = N$ and which is nondegenerate in the sense of Bott, a
neighbourhood $V \subset T^*U$ of the zero-section, and a symplectic embedding
\[
\psi: V \longrightarrow M \setminus \partial M
\]
with $\psi|U = \id$ and $\psi^{-1}(L_1) = \Gamma_{df} \cap V$. Here $T^*U$ carries the
standard symplectic structure, $U$ is considered to be embedded in $V$ as the zero-section,
and $\Gamma_{df} \subset T^*U$ is the graph of $df$. Weinstein's construction can easily be
adapted to take into account the presence of a finite symmetry group. In our situation this
means that one can take $U,f$ to be invariant under $\iota|L_0$, $V$ to be invariant under the
induced action of $\Z/2$ on $T^*U$, and $\psi$ to be an equivariant embedding.

Let $h$ be a function as in \ref{item:clean}. One can find a $\iota$-invariant function $g \in
\smooth_c(U,\R)$ such that $g|N = h$. An elementary argument shows that for all sufficiently
small $t>0$, $f + t\,g$ is a Morse function whose critical point set is $N \cap dg^{-1}(0) =
dh^{-1}(0) = N^\iota$. Define a Lagrangian submanifold $L'_1 \subset M$ by setting $L'_1 \cap
(M \setminus \im\,\psi) = L_1$ and $\psi^{-1}(L'_1) = \Gamma_{df + t\,dg} \cap V$ for some
small $t>0$. $L'_1$ is obviously $\iota$-invariant, and can be deformed back to $L_1$ by an
exact Lagrangian isotopy rel $\partial M$. Moreover, $L_0 \cap L'_1$ consists of the critical
points of $f + tg$, hence is equal to $N^\iota$.

The $\iota$-invariant parts of $L_0,L_1$ are Lagrangian submanifolds of $M^\iota$ which
intersect cleanly in $N^\iota$. By assumption \ref{item:clean} $N^\iota$ is a finite set, so
that $L_0^\iota$ and $L_1^\iota$ actually intersect transversally. Since $(L_1')^\iota$ is a
$C^1$-small perturbation of $(L_1)^\iota$, there is a homeomorphism of $M^\iota$ which sends
$L_0^\iota$ to itself and carries $L_1^\iota$ to $(L_1')^\iota$. As a consequence of this,
\ref{item:no-invariants} continues to hold with $(L_0,L_1)$ replaced by $(L_0,L'_1)$. We have
now proved that the pair $(L_0,L_1')$ satisfies all the assumptions of Lemma
\ref{th:first-symmetry}. Using that Lemma and Proposition \ref{th:isotopy-invariance} one
obtains
\[
\dim HF(L_0,L_1) = \dim HF(L_0,L_1') = |L_0 \cap L_1'| = |N^\iota|.
\]
It remains to explain why $|N^\iota|$ is equal to the dimension of $H^*(N;\Z/2)$. This is in
fact a consequence of the finite-dimensional (Morse theory) analogue of Lemma
\ref{th:first-symmetry}. Consider the gradient flow of $h$ with respect to some
$\Z/2$-invariant metric on $N$. Since $N^\iota$ is a finite set, there are no $\Z/2$-invariant
gradient flow lines. Therefore one can perturb the metric equivariantly so that the moduli
spaces of gradient flow lines become regular. Then $\iota$ induces a free $\Z/2$-action on
these moduli spaces, which implies that the differential in the Morse cohomology complex (with
coefficients in $\Z/2$) for $h$ is zero. Hence $\dim H^*(N;\Z/2) = |dh^{-1}(0)| = |N^\iota|$.
\qed

%
%
\subsection{The grading on Floer cohomology\label{subsec:grading}}

The basic idea in this section is due to Kontsevich \cite{kontsevich94}; for a
detailed exposition see \cite{seidel99}. Let $(M,\o,\alpha,\theta)$ be as in
Section \ref{subsec:floer}, with $\dim M = 2n$; in addition assume that $2
c_1(M,\o)$ is zero.

Let $\LL \rightarrow M$ be the natural fibre bundle whose fibres $\LL_x$ are
the Lagrangian Grassmannians (the manifolds of all linear Lagrangian subspaces)
of $T_xM$. Recall that $\pi_1(\LL_x) \iso \Z$ and that $H^1(\LL_x;\Z)$ has a
canonical generator, the Maslov class $C_x$. Because of our assumption on the
first Chern class, there is a cohomology class on $\LL$ which restricts to
$C_x$ for any $x$. Correspondingly there is an infinite cyclic covering $\tLL
\longrightarrow \LL$ such that each restriction $\tLL_x \longrightarrow \LL_x$
is the universal covering. Fix one such covering (there may be many
non-isomorphic ones in general, so that there is some nontrivial choice to be
made here) and denote the $\Z$-action on it by $\chi_{\tLL}$. For any
Lagrangian submanifold $L$ in $(M,\o,\alpha)$ there is a canonical section
$s_L: L \longrightarrow \LL$ given by $s_L(x) = T_xL$. A {\em grading of $L$}
is a lift $\tL$ of $s_L$ to $\tLL$. Clearly, a grading exists iff $s_L^*\tLL$
is a trivial covering of $L$, so that the obstruction lies in $H^1(L)$. Pairs
$(L,\tL)$ are called {\em graded Lagrangian submanifolds}; we will usually
write $\tL$ instead of $(L,\tL)$. $\chi_{\tLL}$ defines a $\Z$-action on the
set of graded Lagrangian submanifolds. There is also a notion of isotopy, whose
definition is obvious.

An equivalent formulation of the theory goes as follows. Let $J$ be an $\o$-compatible almost
complex structure on $M$. By assumption the bicanonical bundle $K^{\otimes 2}$ is trivial; any
nowhere zero section $\Theta$ of it determines a map
\[
\delta_{\LL}: \LL \longrightarrow \C^*/\Rg,
\]
defined by $\delta_{\LL}(\R e_1 \oplus \dots \oplus \R e_n) = \Theta((e_1 \wedge \dots \wedge
e_n)^{\otimes 2})$ for any family of $n$ orthonormal vectors $e_1,\dots,e_n \in T_xM$ spanning
a Lagrangian subspace. The pullback of the universal cover $\R \longrightarrow \C^*/\Rg$, $\xi
\mapsto \exp(2\pi i \xi)$ by $\delta_{\LL}$ is a covering $\tLL$ of the kind considered above,
and one can show that all such coverings (up to isomorphism) can be obtained in this way. From
this point of view, a grading of a Lagrangian submanifold $L$ is just a map $\tilde{L}: L
\longrightarrow \R$ which is a lift of $L \longrightarrow \C^*/\Rg$, $x \mapsto
\delta_{\LL}(T_xL)$; and the $\Z$-action on the set of gradings is by adding constant
functions.

\begin{example}
Suppose that $M$ is the intersection of a complex hypersurface $g^{-1}(0) \subset \C^{n+1}$
with some ball in $\C^{n+1}$, and that $\o$ is the restriction of the standard symplectic
form. In this case one can take $J$ to be the standard complex structure, and $\Theta$ to be
the square of the complex $n$-form $\det_{\C}(\overline{dg_x} \wedge \dots)|TM$. Then, if
$e_1,\dots,e_n$ is an orthonormal basis of $T_xL$, one has
\begin{equation} \label{eq:explicit-deltal}
\delta_{\LL}(T_xL) = \textstyle{\det_{\C}}(\overline{dg_x} \wedge e_1 \wedge \dots \wedge
e_n)^2.
\end{equation}
\end{example}

Let $(L_0,\tL_0)$ and $(L_1,\tL_1)$ be a pair of graded Lagrangian submanifolds. To any point
$x \in L_0 \cap L_1$ one can associate an absolute Maslov index $\muabs(\tL_0,\tL_1;x) \in
\Z$. To define it, take a path $\tilde{\lambda}_0: [0;1] \longrightarrow \widetilde{\LL}_x$
with endpoints $\tilde{\lambda}_0(j) = \tL_j(x)$, $j = 0,1$. The projection of
$\tilde{\lambda}$ to $\LL$ is a path $\lambda_0: [0;1] \longrightarrow \LL_x$ with
$\lambda_0(j) = T_xL_j$. Let $\lambda_1: [0;1] \longrightarrow \LL_x$ be the constant path
$\lambda_1(t) = T_xL_1$, and set
\[
\muabs(\tL_0,\tL_1;x) = \mupaths(\lambda_0,\lambda_1) + \half(n - \dim(T_xL_0 \cap T_xL_1));
\]
here $\mupaths$ is the `Maslov index for paths' defined in \cite{robbin-salamon93}. The
absolute Maslov index is well-defined, essentially because $\lambda_0$ is unique up to
homotopy rel endpoints. It depends on the gradings through the formula
\begin{equation} \label{eq:shifting}
\muabs(\chi_{\widetilde{\LL}}(r_0)(\tL_0),
 \chi_{\widetilde{\LL}}(r_1)(\tL_1);x) =
 \muabs(\tL_0,\tL_1;x) + r_1 - r_0.
\end{equation}
For later use, we record a situation in which the index is particularly easy to
compute. For any $t \in [0;1]$ one can define a quadratic form $B_t$ on $\lambda_0(t) \cap
\lambda_1(t)$ by setting $B_t(v) = (d/d\tau)_{\tau = t} \o(v,w(\tau))$, where $w$ is a smooth
path in $T_xM$ (defined for $\tau$ near $t$) with $w(t) = v$, $w(\tau) \in \lambda_0(\tau)$
for all $\tau$. In particular, taking $t = 1$ yields a quadratic form $B_1$ defined on
$T_xL_1$.

\begin{lemma} \label{th:compute-index}
Assume that $\lambda_0(t) \cap \lambda_1(t) = \Lambda$ is the same for all $0 \leq t < 1$, and
that the quadratic form $B_1$ is nonpositive with nullspace precisely equal to $\Lambda$. Then
$\muabs(\tL_0,\tL_1;x) = 0$.
\end{lemma}

We omit the proof, which is an immediate consequence of the basic properties of the Maslov
index for paths.

Now let $L_0,L_1$ be two $\theta$-exact Lagrangian submanifolds with $\partial L_0 \cap
\partial L_1 = \emptyset$, and $\tL_0,\tL_1$ gradings of them. To define $HF(L_0,L_1)$ one has
to deform the given $L_j$ to $\kappa$-compatible Lagrangian submanifolds $L_j'$ as in Lemma
\ref{th:perturb}. The $L_j'$ inherit preferred gradings $\tL_j'$ from those of $L_j$ (this is
done by lifting the isotopy from $L_j$ to $L_j'$ to an isotopy of graded Lagrangian
submanifolds). Introduce a $\Z$-grading on $CF(L_0',L_1')$ by setting
\[
CF^r(L_0',L_1') = \bigoplus_{\muabs(\tL_0',\tL_1';x) = r} \Z/2\gen{x}.
\]
Then $d_\J$, for any $\J \in \JJ^\reg(L_0',L_1')$, has degree one, so that one gets a
$\Z$-grading on $HF(L_0',L_1';\J)$. One can show that the continuation isomorphisms between
these groups for various choices of $L_0',L_1'$ and $\J$ have degree zero. Hence there is a
well-defined graded Floer cohomology group $HF^*(\tL_0,\tL_1)$, which satisfies
\[
HF^*(\chi_{\widetilde{\LL}}(r_0)\tL_0,
 \chi_{\widetilde{\LL}}(r_1)\tL_1) =
 HF^{*+r_0-r_1}(\tL_0,\tL_1).
\]
In view of the definition of Floer cohomology, the whole discussion carries over to pairs of
Lagrangian submanifolds with $\partial L_0 = \partial L_1$.

We will now explain the relation between this formalism and the traditional {\em relative
grading} on Floer cohomology, which we have used in the statement of Corollary
\ref{th:comparison}. The relative grading is defined under the assumption that a certain class
$m_{L_0,L_1} \in H^1(\P(L_0,L_1);\Z)$ vanishes; then $HF(L_0,L_1)$ splits as a direct sum of
groups corresponding to the different connected components of $\P(L_0,L_1)$, and for each of
these summands there is a $\Z$-grading which is unique up to a constant shift; see e.g.
\cite{viterbo87}, \cite{floer88e}. If the manifolds $L_0,L_1$ admit gradings, for some
covering $\tLL$, then $m_{L_0,L_1}$ must necessarily be zero, and the absolute grading
$HF^*(\tL_0,\tL_1)$ with respect to any $\tL_0,\tL_1$ is also one of the possible choices of
relative grading. The converse statement is false: there are more choices of relative gradings
than those which occur as $HF^*(\tL_0,\tL_1)$.

\begin{prop} \label{th:graded-symmetry}
Let $L_0,L_1 \subset M$ be two Lagrangian submanifolds as in Proposition \ref{th:symmetry},
and assume that they have gradings $\tL_0,\tL_1$ with respect to some covering $\tLL$. Then
there is a graded isomorphism
\[
HF^*(\tL_0,\tL_1) \iso \bigoplus_{k=1}^d H^{*-\mu_k}(N_k;\Z/2).
\]
Here $N_1,\dots,N_d$ are the connected components of $N = L_0 \cap L_1$; and $\mu_k$ is the
index $\muabs(\tL_0,\tL_1;x)$ at some point $x \in N_k$ (this is independent of the choice of
$x$ because the intersection $L_0 \cap L_1$ is clean).
\end{prop}

The proof is essentially the same as that of Proposition \ref{th:symmetry}. There is just one
additional computation, which describes the change in the Maslov index when $L_0,L_1$ are
perturbed to make the intersection transverse. Since this is a well-known issue, which has
been addressed e.g.\ in \cite{pozniak}, we will not discuss it here.

\section{The symplectic geometry of the $(A_m)$-singularities}

We will first define the general notion of symplectic monodromy map for deformations of
isolated hypersurface singularities. After that we consider the Milnor fibres of the
singularities of type $(A_m)$ in more detail. Concretely, following a suggestion of Donaldson,
we will set up a map
\[
 \left(
 \begin{matrix} \text{curves on the disc with} \\ \text{$(m+1)$ marked points} \end{matrix}
 \right)
 \longmapsto
 \left(
 \begin{matrix} \text{Lagrangian submani-} \\ \text{folds of the Milnor fibre} \end{matrix}
 \right)
\]
Here ``curve'' is taken in the sense of Section \ref{subsec:curves}. At least
for those curves which meet at least one marked point, the Floer cohomology of
the associated Lagrangian submanifolds is well-defined and recovers the
geometric intersection number of the curves. This quickly leads to a proof of
our main symplectic results, Theorem \ref{th:gin-two} and Corollary
\ref{th:faithful-two}. We will also carry out a refined version of the same
argument, which involves Floer cohomology as a graded group. A comparison
between this and the results of Section \ref{subsec:bigrdim} yields Corollary
\ref{th:comparison}.

%
%
\subsection{More basic symplectic geometry\label{subsec:families}}

Throughout this section, $(M,\o,\alpha)$ is a fixed symplectic manifold with contact type
boundary, such that the relative symplectic class $[\o,\alpha]$ is zero.

Let $W$ be a connected manifold and $\widehat{W} \subset W$ a nonempty submanifold. A {\em
symplectic fibration} over $(W,\widehat{W})$ consists of the following data:
\begin{Flist}
\item \label{item:f1}
a proper smooth fibration $\pi: E \longrightarrow W$ whose fibres $E_w = \pi^{-1}(w)$ are
manifolds with boundary. If $W$ also has a boundary, $E$ is a manifold with corners; in that
case we adopt the convention that $\partial E$ refers just to the boundary in fibre direction,
which is the union of the boundaries $\partial E_w$ of the fibres;
\item \label{item:f2}
the structure of a symplectic manifold with contact type boundary on each fibre. More
precisely, we want differential forms $\o_w \in \Omega^2(E_w)$ and $\alpha_w \in
\Omega^1(\partial E_w)$, depending smoothly on $w \in W$ and satisfying the obvious
conditions. We also require that $[\o_w,\alpha_w] \in H^2(E_w,\partial E_w;\R)$ is zero for
all $w$;
\item \label{item:f3}
a trivialization of the family of contact manifolds $(\partial E_w,\alpha_w)_{w \in W}$. This
is a fibrewise diffeomorphism $\tau: \partial M \times W \longrightarrow \partial E$, such
that $\tau_w: (\partial M,\alpha) \longrightarrow (\partial E_w,\alpha_w)$ is a contact
diffeomorphism for every $w \in W$;
\item \label{item:f4}
a trivialization of $E|\widehat{W}$, compatible with all the given structure. This is a
fibrewise diffeomorphism $\eta: M \times \widehat{W} \longrightarrow E|\widehat{W}$ such that
$\eta_w^*\o_w = \o$ and $(\eta_w|\partial M)^*\alpha_w = \alpha$ for each $w \in \widehat{W}$,
and which agrees with $\tau$ on $\partial M \times \widehat{W}$.
\end{Flist}

Two symplectic fibrations over $(W,\widehat{W})$ are {\em cobordant} if they can be joined by
a symplectic fibration over $(W \times [0;1], \widehat{W} \times [0;1])$. A symplectic
fibration is called {\em strict} if $\tau_w^*\alpha_w = \alpha$ for all $w$ (for a general
symplectic fibration, this holds only as an equality of contact structures).

What we have defined is actually a family of symplectic manifolds with contact type boundary,
equipped with certain additional data. Still, we prefer to use the name symplectic fibration
for the sake of brevity; for the same reason we usually write $E$ rather than
$(E,\pi,(\o_w),(\alpha_w),\tau,\eta)$.

\begin{example} \label{ex:interval}
Let $\phi$ be a map in $\Symp(M,\partial M,\o)$. Consider the projection $\pi: M \times [0;1]
\longrightarrow [0;1]$. Define $\eta_0: M \longrightarrow \pi^{-1}(0)$ to be the identity, and
$\eta_1: M \longrightarrow \pi^{-1}(1)$ to be $\phi^{-1}$. This, with the remaining data
chosen in the obvious way, defines a strict symplectic fibration over $([0;1],\{0;1\})$, which
we denote by $E_\phi$.
\end{example}

\begin{prop} \label{th:monodromy}
One can associate to any symplectic fibration $E$ over $(W,\widehat{W})$ a map
$\rho_s^E: \pi_1(W,\widehat{W}) \longrightarrow \pi_0(\Symp(M,\partial M,\o))$,
such that
\begin{Alist}
\item \label{item:naturality}
$\rho_s^E$ is natural with respect to the pullback of $E$ by smooth maps on the base;
\item \label{item:cobordism}
cobordant symplectic fibrations have the same maps $\rho_s^E$;
\item \label{item:interval}
if $E = E_\phi$ is as in Example \ref{ex:interval} and $[\beta] \in \pi_1(W,\widehat{W})$ is
the class of the identity map $\beta = \id_{[0;1]}$, then $\rho_s^E([\beta]) = [\phi]$.
\end{Alist}
Conversely, these properties determine the assignment $E \mapsto \rho_s^E$ uniquely.
\end{prop}

The nontrivial issue here is that different fibres of a symplectic fibration
may not be symplectically isomorphic; not even their volumes need to be the
same. We will therefore define $\rho_s^E$ first for strict symplectic
fibrations, where the problem does not occur, and then extend the definition
using cobordisms. An alternative approach is outlined in Remark
\ref{th:alternative} below.

Take a symplectic fibration $E$ over $(W,\widehat{W})$, and let $\psi \in
\smooth(\partial M \times W,\R)$ be the function defined by $\tau_w^*\alpha_w =
e^{\psi_w}\alpha$; by assumption, it vanishes on $\partial M \times
\widehat{W}$. Choose $\theta \in \Omega^1(M)$ with $\theta|\partial M =
\alpha$, $d\theta = \o$, and let $\kappa: \Rleq \times \partial M \rightarrow
M$ be the embedding determined by the corresponding Liouville vector field.
Similarly choose a smooth family $\theta_w \in \Omega^1(E_w)$, $w \in W$, with
$\theta_w|\partial E_w = \alpha_w$ and $d\theta_w = \o_w$ for all $w$, and such
that $\eta_w^*\theta_w = \theta$ whenever $w \in \widehat{W}$. This determines
an embedding of $\R^{\leq 0} \times \partial E$ into $E$ which fibers over $W$.
By combining this with $\tau$, one obtains an embedding $\gamma: \Rleq \times
\partial M \times W \rightarrow E$ which satisfies $\gamma_w^*\o =
d(e^{r+\psi_w}\alpha)$, and such that $\eta_w^{-1} \circ \gamma_w = \kappa$ for
$w \in \widehat{W}$.

\begin{lemma} \label{th:fibre-bundle}
Any strict symplectic fibration $E$ over $(W,\widehat{W})$ has the structure of
a fibre bundle with structure group $\Symp(M,\partial M,\o)$, compatible with
$\eta$.
\end{lemma}

\proof More precisely, we will show that for any point of $W$ there is a
neighbourhood $U$ and a local trivialization $\eta_U: M \times U \rightarrow
E|U$ which agrees with $\eta$ on $M \times (U \cap \widehat{W})$, such that
\[
\eta_U(\kappa(r,x),w) = \gamma(r,x,w)
\]
for all $(r,x,w)$ in a neighbourhood of $\{0\} \times \partial M \times U
\subset \Rleq \times \partial M \times U$. The transition maps between two
local trivializations with this property obviously take values in
$\Symp(M,\partial M,\o)$, providing the desired structure.

The first step is to find a two-form $\Omega \in \Omega^2(E)$ with the
following properties:
\begin{Olist}
\item \label{item:o1}
$\Omega|E_w = \o_w$,
\item \label{item:o2}
$d\Omega = 0$,
\item \label{item:o3}
$\gamma^*\Omega \in \Omega^2(\R^{\leq 0} \times \partial M \times W)$ agrees
with the pullback of $d(e^r\alpha) \in \Omega^2(\R^{\leq 0} \times \partial M)$
on some neighbourhood of $\{0\} \times \partial M \times W$,
\item \label{item:o4}
$\eta^*\Omega \in \Omega^2(M \times \widehat{W})$ is the pullback of $\omega$.
\end{Olist}
This is not difficult. One takes a one-form $\Theta$ on $E$ with $\Theta|E_w =
\theta_w$, such that $\gamma^*\Theta$ is equal to the pullback of $e^r\alpha$
near $\{0\} \times \partial M \times W$, and with $\eta^*\Theta$ the pullback
of $\theta$ to $M \times \widehat{W}$. Such a one-form exists because the
fibration is strict. Then $\Omega = d\Theta$ satisfies the conditions above.

Define, for any $x \in E$, a subspace $H_x = \{X \in T_xE \suchthat \O(X,Y) = 0
\text{ for all } Y \in T_xE \text{ with } D\pi(Y) = 0\}$. The subbundle $H
\subset TE$ formed by the $H_x$ is complementary to $\ker\,D\pi \subset TE$ by
\ref{item:o1} and is tangent to $\partial E$ by \ref{item:o3}; hence it defines
a connection on $\pi: E \longrightarrow W$. The parallel transport maps of this
connection are symplectic by \ref{item:o2} (the reader can find a more
extensive discussion of symplectic connections in \cite[Chapter
6]{mcduff-salamon96}). Moreover, with respect to the map $\gamma$, they are
trivial near $\partial E$ by \ref{item:o3}. Using those maps, one easily
obtains the desired local trivializations. \qed

As a special case of this construction, one sees that any strict symplectic
fibration over $([0;1],\{0;1\})$ is isomorphic to one of those in Example
\ref{ex:interval}.

\begin{lemma} \label{th:strict}
Any symplectic fibration is cobordant to a strict one.
\end{lemma}

\proof Take a symplectic fibration $E$, and define $\psi$ and $\gamma$ as
before. Let $c: W \longrightarrow \R^{\geq 0}$ be given by $c(w) = \max
(\{\psi_w(x) \suchthat x \in \partial M\} \cup \{0\})$. Choose some
$\epsilon>0$ and a function $\xi \in \smooth(\Rleq \times W,\R)$, such that the
restriction $\xi_w$ has the following properties for each $w \in W$: $\xi_w(r)
= 0$ for $r \leq -c(w)-2\epsilon$, $\xi_w(r) = 1$ for $r \geq -\epsilon$, and
$0 \leq \xi_w'(r) < c(w)^{-1}$ for $r \in \R^{\leq 0}$. For $(w,s) \in W \times
[0;1]$ define $\tilde{\o}_{w,s} \in \Omega^2(E_w)$ by setting $\tilde{\o}_{w,s}
= \o_w$ outside $\gamma_w([-c(w)-2\epsilon;0] \times \partial M)$ and
\begin{equation} \label{eq:pullback}
\gamma_w^*\tilde{\o}_{w,s} = d(e^{r + \psi_w - s\xi_w(r)\psi_w}\alpha).
\end{equation}
$\tilde{\o}_{w,s}$ is a symplectic form on $E_w$. In fact, the top exterior
power of the two-form on the r.h.s\ of \eqref{eq:pullback} is
$(1-s\xi'_w(r)\psi_w)\, dr \wedge \alpha \wedge m (d\alpha)^{m-1}$, and this is
everywhere nonzero by assumption on $\xi'_w$. Define $\tilde{\alpha}_{w,s} \in
\Omega^1(\partial E_w)$ by $\tau_w^*\tilde{\alpha}_{w,s} =
\exp((1-s)\psi_w)\alpha$. In the same way one can construct one-forms
$\tilde{\theta}_{w,s}$ on $E_w$ with $d\tilde{\theta}_{w,s} = \tilde{\o}_{w,s}$
and $\tilde{\theta}_{w,s}|\partial E_w = \tilde{\alpha}_{w,s}$, which shows
that the relative cohomology class $[\tilde{\o}_{w,s},\tilde{\alpha}_{w,s}]$ is
zero. $(E_w,\tilde{\o}_{w,s},\tilde{\alpha}_{w,s})$ is a symplectic manifold
with contact type boundary for all $(w,s) \in W \times [0;1]$. For $s = 0$
these are the original fibres $(E_w,\o_w,\alpha_w)$; the same holds for
arbitrary $s$ when $w \in \widehat{W}$, since then $\psi_w = 0$. Note also that
for $s = 1$, $\tau_w^*\tilde{\alpha}_{w,1} = \alpha$.  Finally, by putting
these forms onto $\widetilde{E} = E \times [0;1] \rightarrow W \times [0;1]$
and taking $\tilde{\tau},\tilde{\eta}$ to be induced from $\tau,\eta$ in the
obvious way, one obtains a cobordism between $E$ and the strict symplectic
fibration $\widetilde{E}_1 = \widetilde{E}| W \times \{1\}$. \qed

\proof[Proof of Proposition \ref{th:monodromy}] For a strict symplectic
fibration $E$ one defines $\rho_s^E$ in terms of the fibre bundle structure
from Lemma \ref{th:fibre-bundle}. To show that this is independent of all
choices, it is sufficient to observe that any two fibre bundle structures on
$E$ obtained from different choices of $\Omega$ can be joined by a fibre bundle
structure on $E \times [0;1]$. As mentioned earlier, the definition is then
extended to general symplectic fibrations by using Lemma \ref{th:strict}. There
is again a uniqueness part of the proof, which consists in showing that two
strict symplectic fibrations which are cobordant are also cobordant by a strict
cobordism; this is a relative version of Lemma \ref{th:strict}. The properties
\ref{item:naturality}--\ref{item:interval} follow immediately from the
definition; and the fact that they characterize $\rho_s^E$ is obvious. \qed

Our construction implies that $\rho_s^E$ is multiplicative in the following
sense: if $\beta_0,\beta_1$ are two paths $([0;1],\{0;1\}) \longrightarrow
(W,\widehat{W})$ with $\beta_0(1) = \beta_1(0)$, then $\rho_s^E([\beta_1 \circ
\beta_0]) = \rho_s^E([\beta_1]) \rho_s^E([\beta_0])$. In the special case where
$\widehat{W} = \{\hat{w}\}$ is a single point, one obtains a group homomorphism
$\pi_1(W,\hat{w}) \longrightarrow \pi_0(\Symp(M,\partial M,\o))$ (because of
the way in which we write the composition of paths, our fundamental groups are
the opposites of the usual ones). As a final observation, note that $\rho_s^E$
depends only on the homotopy class (rel $\widehat{W}$) of the contact
trivialization $\tau$; the reason is that changing $\tau$ within such a class
gives a fibration which is cobordant to the original one.

\begin{remark} \label{th:alternative}
An alternative way of proving Proposition \ref{th:monodromy} is to attach an
infinite cone to $M$, forming a noncompact symplectic manifold
$(\overline{M},\bar{\o}) = (M,\o) \cup_{\partial M} (\R^{\geq 0} \times
\partial M, d(e^r\alpha))$. Because $[\o,\alpha] = 0$, the inclusion
$\Symp(M,\partial M,\o) \hookrightarrow \Symp^c(\overline{M},\bar{\o})$ is a
weak homotopy equivalence. Given a symplectic fibration $E$, one similarly
attaches an infinite cone to each fibre. The resulting family of non-compact
symplectic manifolds admits a natural structure of a fibre bundle with the
group $\Symp^c(\overline{M},\bar{\o})$ as structure group, and this defines a
map
\[
\pi_1(W,\widehat{W}) \longrightarrow \pi_0(\Symp^c(\overline{M},\bar{\o})) \iso
\pi_0(\Symp(M,\partial M,\o)).
\]
\end{remark}

\begin{lemma} \label{th:relative-monodromy-one}
Let $E$ be a symplectic fibration over $(W,\widehat{W}) = ([0;1],\{0;1\})$. Let $F \subset E
\setminus \partial E$ be a compact submanifold which intersects each fibre $E_w$
transversally, such that $F_w = F \cap E_w \subset E_w$ is a closed Lagrangian submanifold 
for all $w$. Define Lagrangian submanifolds $L_0,L_1 \subset M$
by $L_w = \eta_w^{-1}(F_w)$. Let $\phi \in \Symp(M,\partial M,\o)$ be a map which represents
$\rho_s([\beta])$, where $\beta = \id_{[0;1]}$. Then $L_1$ is Lagrangian isotopic to
$\phi(L_0)$.
\end{lemma}

\proof The result is obvious when $E$ is strict; as remarked after Lemma \ref{th:fibre-bundle}, 
we may then assume that $E = E_\phi$ is as in Example \ref{ex:interval}, and 
the submanifolds $F_w \subset E_w = M$ will form a Lagrangian isotopy between $L_0$ and
$\phi^{-1}(L_1)$.

Now consider the general case; for this we adopt the notation from the proof of Lemma
\ref{th:strict}. The vector fields dual to the one-forms $\theta_w$ define a semiflow
$(\lambda_s)_{s \leq 0}$ on $E$, whose restriction to each fibre is Liouville (rescales the
symplectic form exponentially). Define
\[
\widetilde{F} = \lambda_s(F) \times [0;1] \subset \widetilde{E},
\]
where $s \in \R$ is a number smaller than $-c(w)-2\epsilon$ for all $w \in W$. The restriction
$\widetilde{F}_{w,s} = \widetilde{F} \cap \widetilde{E}_{w,s}$ is a closed Lagrangian
submanifold in $(\widetilde{E}_{w,s},\tilde{\o}_{w,s},\tilde{\alpha}_{w,s})$ for all $(w,s)
\in W \times [0;1]$. By restricting to the subset $W \times \{1\}$, on which $\widetilde{E}$
becomes strict, and using the observation made above concerning the strict case together with
the cobordism invariance of $\rho_s^E$, one finds that $\phi(\eta_0^{-1}(\lambda_s(F_0)))$ is
Lagrangian isotopic to $\eta_1^{-1}(\lambda_s(F_1))$. Moreover, the flow $\lambda_s$ provides
symplectic isotopies between $\eta_w^{-1}(\lambda_s(F_w))$ and $L_w$ for $w = 0,1$. \qed

\begin{lemma} \label{th:relative-monodromy-two}
Let $E$ be a symplectic fibration over $(W,\widehat{W}) = ([0;1],\{0;1\})$. Let $F \subset E$
be a compact submanifold with which intersects each fibre $E_w$ transversally, and such that
$F_w = F \cap E_w \subset E_w$ is a Lagrangian submanifold with nonempty boundary 
for all $w$. Assume that the path in the space of all Legendrian submanifolds of $(\partial
M,\alpha)$ defined by $w \mapsto \tau_w^{-1}(\partial F_w)$ is closed and contractible. Then,
for $L_0,L_1$ and $\phi$ as in the previous Lemma, $L_1$ is isotopic to $\phi(L_0)$ by a
Lagrangian isotopy rel $\partial M$.
\end{lemma}

This is the analogue of Lemma \ref{th:relative-monodromy-one} for Lagrangian submanifolds with
nonempty boundary. The proof uses Lemma \ref{th:compatible} and the same ideas as before;
we omit it.

%
%
\subsection{Singularities}

To begin with, some notation: $B^{2k}(r)$, $\barB^{2k}(r)$, $S^{2k-1}(r)$ denote the open
ball, the closed ball, and the sphere of radius $r>0$ around the origin in $\C^k$;
$\theta_{\C^k} = i/4 \sum_j (x_j d\bar{x}_j - \bar{x}_j dx_j) \in \Omega^1(\C^k)$ and
$\omega_{\C^k} = d\theta_{\C^k} \in \Omega^2(\C^k)$ are the standard forms.

Let $g \in \C[x_0,\dots,x_n]$ be a polynomial with $g(0) = 0$ and $dg(0) = 0$, such that $x =
0$ is an isolated critical point. A theorem of Milnor says that $H_0 = g^{-1}(0)$ intersects
$S^{2n+1}(\epsilon)$ transversally for all sufficiently small $\epsilon>0$. Fix an $\epsilon$
in that range. Assume that we are given a family $(g_w)_{w = (w_0,\dots,w_m) \in \C^{m+1}}$ of
polynomials of the form $g_w(x) = g(x) + \sum_{j=0}^m w_j\tilde{g}_j(x)$, where $\tilde{g}_0
\equiv 1$. The intersection of $H_w = g_w^{-1}(0)$ with $S^{2n+1}(\epsilon)$ remains
transverse for all sufficiently small $w$. Fix a $\delta>0$ such that this is the case for all
$|w| < \delta$.

\begin{lemma} \label{th:gray}
The smooth family of contact manifolds $(H_w \cap S^{2n+1}(\epsilon),\theta_w|H_w \cap
S^{2n+1}(\epsilon))$, $w \in B^{2m+2}(\delta)$, admits a trivialization, which is unique up to
homotopy.
\end{lemma}

This is an immediate consequence of Gray's theorem, because $B^{2m+2}(\delta)$ is
contractible. Let $W \subset B^{2m+2}(\delta)$ be the subset of those $w$ such that $H_w \cap
\barB^{2n+2}(\epsilon)$ contains no singular point of $H_w$. It is open and connected, because
its complement is a complex hypersurface. For each $w \in W$ equip $E_w = H_w \cap
\barB^{2n+2}(\epsilon)$ with the forms $\o_w = \o_{\C^{n+1}}|E_w$ and $\theta_w =
\theta_{\C^{n+1}}|E_w$, and set $\alpha_w = \theta_w|\partial E_w$. Then each
$(E_w,\o_w,\alpha_w)$ is a symplectic manifold with contact boundary, and the relative
symplectic class $[\o_w,\alpha_w]$ is zero because of the existence of $\theta_w$. Consider
the smooth fibration
\[
\pi: E = \{ (x,w) \in \C^{n+1} \times W \suchthat |x| \leq \epsilon, \; g_w(x) = 0\}
\longrightarrow W
\]
whose fibres are the manifolds $E_w$. Choose some base point $\hat{w} \in W$ and set
$(M,\o,\alpha) = (E_{\hat{w}}, \o_{\hat{w}}, \alpha_{\hat{w}})$. By restriction, any
trivialization as in Lemma \ref{th:gray} defines a contact trivialization $\tau: \partial M
\times W \longrightarrow \partial E$. Take $\eta: M \longrightarrow E_{\hat{w}}$ to be the
identity map. What we have defined is a symplectic fibration $E$ over $(W,\{\hat{w}\})$.
By Proposition \ref{th:monodromy} one can associate to it a homomorphism
\[
\rho_s = \rho_s^E: \pi_1(W,\hat{w}) \longrightarrow \pi_0(\Symp(M,\partial M,\o)).
\]
In the usual terminology, $E$ is the Milnor fibration of the deformation $(H_w)$ of the
singularity $0 \in H_0$, and the manifolds $E_w$ are the Milnor fibres. We call $\rho_s$ the
symplectic monodromy map associated to the deformation.

We need to discuss briefly how $\rho_s$ depends on the various parameters involved. $\tau$ is
unique up to homotopy, so that $\rho_s$ is independent of it. Changing $\delta$ affects the
space $W$, but the fundamental group remains the same for all sufficiently small $\delta$, and
we will assume from now on that $\delta$ has been chosen in that range. The dependence on the
base point $\hat{w}$ is a slightly more complicated issue. Different fibres of the Milnor
fibration are not necessarily symplectically isomorphic, but they become isomorphic after
attaching an infinite cone, as in Remark \ref{th:alternative}. This means that, given a path
in $W$ from $w_0$ to $w_1$, one can identify
\[
\pi_0(\Symp(E_{w_0},\partial E_{w_0},\o_{w_0})) \iso \pi_0(\Symp(E_{w_1},\partial
E_{w_1},\o_{w_1})),
\]
and this fits into a commutative diagram with the corresponding isomorphism $\pi_1(W,w_0) \iso
\pi_1(W,w_1)$ and with the symplectic monodromy maps at these base points. For a similar
reason, making $\epsilon$
smaller changes $M$ but does not affect $\pi_0(\Symp(M,\partial M,\o))$ or $\rho_s$. 
Finally (but we will not need this) one can
choose $(g_w)$ to be a miniversal deformation of $g$, and then the resulting symplectic
monodromy map is really an invariant of the singularity.

\begin{remark}
Apart from the present paper, symplectic monodromy has been studied only for the one-parameter
deformations $g_w(x) = g(x) + w$. In that case $\pi_1(W,\hat{w}) = \Z$, so that it is
sufficient to consider the class $[\phi] = \rho_s(1) \in \pi_0(\Symp(M,\partial M,\o))$. When
$g$ is the ordinary singularity $g(x) = x_0^2 + \dots + x_n^2$, $M$ is isomorphic to a
neighbourhood of the zero-section in $T^*S^n$, and $\phi$ can be chosen to be the generalized
Dehn twist along the zero-section in the sense of \cite{seidel98b}, \cite{seidel99}; those
papers show that $[\phi]$ has infinite order in all dimensions $n \geq 1$. We should say that
the relevance of symplectic geometry in this example was first pointed out by Arnol'd
\cite{arnold95}.

More generally, in \cite{seidel99} it is proved that $[\phi]$ has infinite order for all
weighted homogeneous singularities such that the sum of the weights is $\neq 1$ (the
definition of symplectic monodromy in \cite{seidel99} differs slightly from that here, but it
can be shown that the outcome is the same). The general question, {\em does $[\phi]$ have
infinite order for all nontrivial isolated hypersurface singularities?}, is open.
\end{remark}

%
%
\subsection{From curves to Lagrangian submanifolds}

We continue to use the notation introduced in the previous section. The $n$-dimensional
singularity of type $(A_m)$, for $m,n \geq 1$, is given by $g(x_0,\dots,x_n) = x_0^2 + x_1^2 +
\dots + x_{n-1}^2 + x_n^{m+1}$. We consider the $(m+1)$-parameter deformation $g_w(x) = x_0^2
+ \dots + x_{n-1}^2 + h_w(x_n)$, where $h_w(z) = z^{m+1} + w_mz^m + \dots + w_1z + w_0 \in
\C[z]$ (this is slightly larger than the usual miniversal deformation, but that makes no
difference as far as the monodromy map is concerned). Since $g$ is weighted homogeneous, one
can take $\epsilon = 1$. In contrast, we will make use of the right to choose $\delta$ small.
For $w = (w_0,\dots,w_m) \in B^{2m+2}(\delta)$ set $D_w = \{z \in \C \suchthat |h_w(z)| +
|z|^2 \leq 1\}$, and $\Delta_w = h_w^{-1}(0) \subset \C$.

\begin{lemma} \label{th:convex}
Provided that $\delta$ has been chosen sufficiently small, the following properties hold for
any $w \in B^{2m+2}(\delta)$:
\begin{Plist}
\item \label{item:convex}
$D_w \subset \C$ is an compact convex subset with smooth boundary, which contains $B^2(\half)$
and is contained in $B^2(2)$.
\item \label{item:positivity}
Let $\beta: [0;1] \longrightarrow \partial D_w$ be an embedded path which moves in positive
sense with respect to the obvious orientation. Then
\[
\textstyle{\frac{d}{dt} \arg h_w(\beta(t)) > 0, \quad \frac{d}{dt} \arg \beta(t) > 0}
\]
for all $t \in [0;1]$, where $\arg = \im\,\log: \C^* \longrightarrow \R/2\pi\Z$ is the
argument function.
\item \label{item:delta-w}
$\Delta_w$ is contained in $B^2(\frac{1}{3})$.
\end{Plist}
\end{lemma}

\proof For $w = 0$, $h_0(z) = z^{m+1}$, and $\partial D_0 = \{|z|^{m+1} + |z|^2 = 1\}$ is a
circle centered at the origin, so the conditions are certainly satisfied. Perturbing $w$
slightly will change $h_w$ and $\partial D_w$ only by a small amount in any $C^k$-topology,
and this implies the first two parts. The proof of \ref{item:delta-w} is even simpler. \qed

\begin{lemma}
$W \subset B^{2m+2}(\delta)$ is the subset of those $w$ such that $h_w$ has no multiple zeros.
\end{lemma}

This follows immediately from property \ref{item:delta-w} above. It implies in particular that
for $w \in W$, $\Delta_w \subset D_w \setminus \partial D_w$ is a set of $(m+1)$ points.

We will now study the symplectic geometry of the Milnor fibres $E_w = H_w \cap \barB^{2n+2}(1)
= \{ |x| \leq 1, \; x_0^2 + \dots x_{n-1}^2 + h_w(x_n) = 0\}$ equipped with the forms $\o_w,
\theta_w, \alpha_w$. Consider the projection $H_w \longrightarrow \C$, $x \longmapsto x_n$,
whose fibres are the affine quadrics $Q_{w,z} = \{x_0^2 + \dots + x_{n-1}^2 = -h_w(z), \; x_n
=z\}$. $Q_{w,z}$ is smooth for all $z \notin \Delta_w$, and actually symplectically isomorphic
to $T^*S^{n-1}$. If one chooses a symplectic isomorphism which respects the obvious
$O(n)$-actions, the zero-section $S^{n-1} \subset T^*S^{n-1}$ corresponds to $\Sigma_{w,z} =
\sqrt{-h_w(z)}S^{n-1} \times \{z\} \subset Q_{w,z}$, which is the subset of those $x =
(x_0,\dots,x_{n-1},x_n)$ such that $x_n = z$, $|x_0|^2 + \dots + |x_{n-1}|^2 = |h_w(z)|$, and
$x_i \in \sqrt{-h_w(z)}\R$ for $i = 0,\dots,n-1$. For $z \in \Delta_w$, $Q_{w,z}$ is
homeomorphic to $T^*S^{n-1}$ with the zero-section collapsed to a point, that point being
precisely $\Sigma_{w,z}$. Moreover, by definition of $D_w$,
\[
\Sigma_{w,z} \text{ is}
 \begin{cases}
 \text{contained in $E_w \setminus \partial E_w$}
 & \text{if $z \in D_w \setminus \partial D_w$},\\
 \text{contained in $\partial E_w$}
 & \text{if $z \in \partial D_w$},\\
 \text{disjoint from $E_w$}
 & \text{otherwise.}
 \end{cases}
\]
For $w \in W$ and $c$ a curve in $(D_w,\Delta_w)$ as defined in Section \ref{subsec:curves},
set
\[
L_{w,c} = \bigsqcup_{z \in c} \Sigma_{w,z} \subset E_w.
\]

\begin{lemma} \label{th:donaldson}
\begin{theoremlist}
\item \label{item:donaldson-one}
If $c$ meets $\Delta_w$ then $L_{w,c}$ is a $\theta_w$-exact Lagrangian submanifold of
$(E_w,\o_w,\alpha_w)$.
\item \label{item:donaldson-two}
Let $c_0,c_1$ be two curves in $(D_w,\Delta_w)$ which meet $\Delta_w$ and such that $c_0
\isotopic c_1$. Then $L_{w,c_0}, L_{w,c_1}$ are isotopic by an exact Lagrangian isotopy rel
$\partial E_w$.
\end{theoremlist}
\end{lemma}

\proof \ref{item:donaldson-one} Assume first that $c$ can be parametrized by a smooth
embedding $\gamma: [0;1] \longrightarrow D_w$ with $\gamma^{-1}(\Delta_w) = \{0\}$ and
$\gamma^{-1}(\partial D_w) = \{1\}$. Because $\gamma(0)$ is a simple zero of $h_w$, one can
write $h_w(\gamma(t)) = t\,k(t)$ for some $k \in \smooth([0;1],\C^*)$. Choose a smooth square
root $\sqrt{-k}$, and define a map from the unit ball $\barB^n \subset \R^n$ to $E_w$ by
\begin{equation} \label{eq:parametrization}
y \longmapsto (y \sqrt{-k(|y|)^2}, \gamma(|y|^2)).
\end{equation}
This is a smooth embedding with image $L_{w,c}$. Since $\gamma$ meets $\partial D_w$
transversally, the embedding intersects $\partial E_w$ transversally. Moreover, the pullback
of $\theta_w$ under \eqref{eq:parametrization} is of the form $\psi(|y|^2)d|y|^2$ for some
$\psi \in \smooth([0;1],\R)$. This implies that $\o_w|L_{w,c} = 0$ and $\theta_w|\partial
L_{w,c} = 0$, which means that $L_{w,c}$ is indeed a Lagrangian submanifold.

In the other case, when $c$ joins two points of $\Delta_w$, $L_{w,c}$ is a Lagrangian
submanifold diffeomorphic to the $n$-sphere; this is proved in the same way as before, only
that now one splits $c$ into two parts, and covers $L_{w,c}$ by two smooth charts.

The $\theta_w$-exactness of $L_{w,c}$ is obvious for all $n \geq 2$, since then
$H^1(L_{w,c},\partial L_{w,c};\R) = 0.$ For $n = 1$ what one has to prove is that
\[\int_{L_{w,c}} \theta_w = 0\] for some orientation of $L_{w,c}$. Now the involution
$\iota(x_0,x_1) = (-x_0,x_1)$ preserves $\theta_w$ and reverses the orientation of $L_{w,c}$,
as one can see by looking at the tangent space at any point $(0,z)$, $z \in c \cap \Delta_w$;
this implies the desired result. Part \ref{item:donaldson-two} of the Lemma follows from
\ref{item:donaldson-one} together with Lemma \ref{th:exact-isotopies}. \qed

$L_{w,c}$ is Lagrangian even if $c$ does not meet $\Delta_w$, but then it is no longer
$\theta_w$-exact; which is why we do not use these submanifolds.

\begin{lemma} \label{th:positive-isotopies}
Let $c$ be a curve in $(D_w,\Delta_w)$ which joins a point of $\Delta_w$ with a point of
$\partial D_w$. Take a vector field on $\partial D_w$ which is nonvanishing and positively
oriented. Extend it to a vector field $Z$ on $D_w$ which vanishes near $\Delta_w$, and let
$(f_t)$ be the flow of $Z$. Then the Lagrangian isotopy $t \mapsto L_{w,f_t(c)}$, $0 \leq t
\leq 1$, is exact and also positive in the sense of Section \ref{subsec:symplectic}.
\end{lemma}

\proof The exactness follows again from Lemma \ref{th:exact-isotopies}. Define a smooth path
$\beta: [0;1] \longrightarrow \partial D_w$ by requiring that $\beta(t)$ be the unique point
of $f_t(c) \cap \partial D_w$. Choose some square root of $-h_w \circ \beta$, and consider
\begin{equation} \label{eq:d-para}
S^{n-1}  \times [0;1] \longrightarrow \partial E_w, \quad (y,t) \longmapsto
(\sqrt{-h_w(\beta(t))}y,\beta(t)).
\end{equation}
This maps $S^{n-1} \times \{t\}$ to $\partial L_{w,f_t(c)}$. A computation shows that the
pullback of $\alpha_w$ under \eqref{eq:d-para} is $\psi(t)\,dt$, where $\psi(t) =
|h_w(\beta(t))|\frac{d}{dt}(\arg h_w(\beta(t))) + \half |\beta(t)|^2 \frac{d}{dt}(\arg
\beta(t))$ is a positive function by \ref{item:positivity} in Lemma \ref{th:convex}. This is
one of the definitions of a positive isotopy. \qed

The intersections of the submanifolds $L_{w,c}$ are governed by those of the underlying
curves. In fact $L_{w,c_0} \cap L_{w,c_1} = \bigsqcup_{z \in c_0 \cap c_1} \Sigma_{w,z}$
consists of an $(n-1)$-sphere for any point of $(c_0 \cap c_1) \setminus \Delta_w$ together
with an isolated point for each point of $c_0 \cap c_1 \cap \Delta_w$. The next Lemma
translates this elementary fact into a statement about Floer cohomology.

\begin{lemma} \label{th:floer-cohomology}
Let $c_0,c_1$ be two curves in $(D_w,\Delta_w)$, each of which 
meets $\Delta_w$. Then the dimension of
$HF(L_{w,c_0},L_{w,c_1})$ is $2\,I(c_0,c_1)$.
\end{lemma}

\proof We consider first the case when $c_0 \cap c_1 \cap \partial D_w = \emptyset$ and $c_0
\not \isotopic c_1$. $I(c_0,c_1)$ depends only on the isotopy classes of $c_0,c_1$; the same
is true of $HF(L_{w,c_0},L_{w,c_1})$, by Lemma \ref{th:donaldson}\ref{item:donaldson-two} and
Proposition \ref{th:isotopy-invariance}. Hence we may assume that $c_0,c_1$ have minimal
intersection.

A straightforward computation of the tangent spaces shows that $L_{w,c_0},L_{w,c_1}$ have
clean intersection. Moreover, since $c_0 \cap c_1 \cap \partial D_w = \emptyset$, $L_{w,c_0}
\cap L_{w,c_1} \cap \partial E_w = \emptyset$. Consider the involution $\iota(x_0,\dots,x_n) =
(-x_0,\dots,-x_{n-2},x_{n-1},x_n)$ on $E_w$, which preserves $\o_w$ and $\theta_w$. Each
connected component $\Sigma_{w,z}$ of $L_{w,c_0} \cap L_{w,c_1}$ carries a $\iota$-invariant
Morse function $h_z$ whose critical point set is $\{(0,\dots,\pm \sqrt{-h_w(z)},z)\} =
\Sigma_{w,z} \cap M^\iota$. Namely, for $z \in \Delta_w$ one takes $h_z$ to be constant, and
otherwise $h_z(x) = (-h_w(z))^{-1/2}x_{n-1}$ for some choice of square root. Hence
$L_{w,c_0},L_{w,c_1}$ satisfy condition \ref{item:clean} in Proposition \ref{th:symmetry}.

The fixed point set $M^\iota$ is precisely the double branched cover of $D_w$ from Section
\ref{subsec:branched-cover}, and $L_{w,c_j} \cap M^\iota \subset M^\iota$ is the preimage of
$c_j$ under the covering projection. Lemma \ref{th:lifting-curves} shows that $L_{w,c_0} \cap
M^\iota$ and $L_{w,c_1} \cap M^\iota$ have minimal intersection as curves in $M^\iota$. Hence,
using Lemma \ref{th:no-maps}, it follows that $L_{w,c_0},L_{w,c_1}$ also satisfy condition
\ref{item:no-invariants} in Proposition \ref{th:symmetry}. Applying that Proposition shows
that $\dim HF(L_{w,c_0},L_{w,c_1}) = \dim H^*(L_{w,c_0} \cap L_{w,c_1};\Z/2) = 2|(c_0 \cap
c_1) \setminus \Delta_w| + |c_0 \cap c_1 \cap \Delta_w| = 2\,I(c_0,c_1)$.

Next, consider the case when $c_0 \cap c_1 \cap \partial D_w = \emptyset$ and $c_0 \isotopic
c_1$. This means that $c_0$ must be a path joining two points of $\Delta_w$; hence $I(c_0,c_1)
= I(c_0,c_0) = 1$. On the other hand, since $L_{w,c_0}$ and $L_{w,c_1}$ are isotopic by an
exact Lagrangian isotopy, one has $HF(L_{w,c_0},L_{w,c_1}) \iso HF(L_{w,c_0},L_{w,c_0}) \iso
H^*(L_{w,c_0};\Z/2) \iso (\Z/2)^2$ by Floer's theorem \cite{floer89}.

The final case is $c_0 \cap c_1 \cap \partial D_w \neq \emptyset$. Then $c_0,c_1$ have each
one point on the boundary $\partial D_w$, which is the same for both, so that $\partial
L_{w,c_0} = \partial L_{w,c_1}$. To compute $HF(L_{w,c_0},L_{w,c_1})$ one has to perturb
$L_{w,c_0}$ by an exact and positive Lagrangian isotopy. One can take the isotopy to be
$L_{w,f_t(c_0)}$ for some flow $(f_t)$ on $D_w$ as in Lemma \ref{th:positive-isotopies}. By
definition, $I(c_0,c_1) = I(f_t(c_0),c_1)$ for small $t>0$. This reduces the Floer cohomology
computation to the case which we have already proved. \qed

We are now ready to bring in the symplectic monodromy map. Fix a base point $\hat{w} \in W$.
Write $(M,\o,\theta,\alpha) = (E_{\hat{w}},\o_{\hat{w}},\theta_{\hat{w}},\alpha_{\hat{w}})$,
$(D,\Delta) = (D_{\hat{w}},\Delta_{\hat{w}})$, $h = h_{\hat{w}}$, and $L_c = L_{\hat{w},c}$
for any curve $c$ in $(D,\Delta)$.

One can embed $W$ as an open subset into $\Conf_{m+1}(D \setminus \partial D)$ by sending $w
\mapsto \Delta_w$. The image of this embedding contains all configurations of points which lie
close to the origin, and it is mapped into itself by the maps which shrink configurations by
some factor $0<r<1$. From this one can easily see that the embedding is a weak homotopy
equivalence. By combining this with the observation made in Section \ref{subsec:disc} one
obtains canonical isomorphisms
\begin{equation} \label{eq:canonical-isos}
\pi_1(W,\hat{w}) \iso \pi_1(\Conf_{m+1}(D),\Delta) \iso \pi_0(\Diffeo),
\end{equation}
where $\Diffeo = \Diff(D,\partial D;\Delta)$.

\begin{lemma} \label{th:curves-and-monodromy}
Let $f \in \Diffeo$ and $\phi \in \Symp(M,\partial M,\o)$ be related in the following way:
$[f] \in \pi_0(\Diffeo)$ corresponds to a class in $\pi_1(W,\hat{w})$ whose image under
$\rho_s$ is $[\phi] \in \pi_0(\Symp(M,\partial M,\o))$. Then, for any curve $c$ in
$(D,\Delta)$ which meets $\Delta$, $L_{f(c)}$ and $\phi(L_c)$ are isotopic by a Lagrangian
isotopy rel $\partial M$.
\end{lemma}

\proof Consider first the case when $c \cap \partial D = \emptyset$. Since the result depends
only on the isotopy class of $c$, we may assume that $c \subset B^2(1/2)$. Similarly, the
result depends only on the isotopy class of $f$, and we may assume that $f$ is the identity
outside $B^2(1/2) \subset D$. Take a smooth path $\beta: ([0;1],\{0;1\}) \longrightarrow
(W,\hat{w})$ whose class in $\pi_1(W,\hat{w})$ corresponds to $[f] \in \pi_0(\Diffeo)$ under
the canonical isomorphism. This means that there is an isotopy $(f_t)_{0 \leq t \leq 1}$ in
$\Diff(D,\partial D)$ with $f_0 = \id$, $f_1 = f$, and $f_t(\Delta) = \Delta_{\beta(t)}$ for
all $t$. We may assume that each $f_t$ is equal to the identity outside $B^2(1/2)$. Then
$f_t(c)$ is a curve in $(D_{\beta(t)},\Delta_{\beta(t)})$ for all $t \in [0;1]$. Moreover, the
Lagrangian submanifolds $L_{\beta(t),f_t(c)}$ depend smoothly on $t$; by this we mean that they form a
submanifold $F$ of the pullback $\beta^*E \longrightarrow [0;1]$, as in Lemma
\ref{th:relative-monodromy-one}. By using that Lemma and the naturality of the maps $\rho_s^E$
under pullbacks, it follows that $L_{f(c)}$ is Lagrangian isotopic to $\phi(L_c)$.

The case when $c \cap \partial D \neq \emptyset$ is not very different. We may assume that the
intersection $c \cap (D \setminus B^2(1/2))$ consists of a straight radial piece $\R^{>0}\zeta
\cap (D \setminus B^2(1/2))$, for some $\zeta \in S^1$. We take a smooth path $\beta$ and an
isotopy $(f_t)$ as in the previous case. Then $f_t(c) \cap D_{\beta(t)}$ is a curve in
$(D_{\beta(t)}, \Delta_{\beta(t)})$ for all $0 \leq t \leq 1$. The rest is as before except
that in order to apply Lemma \ref{th:relative-monodromy-two}, one has to show that the path of
Legendrian submanifolds of $(\partial M,\alpha)$ given by
\[
t \mapsto \tau_{\beta(t)}^{-1}(\partial L_{\beta(t),f_t(c) \cap D_{\beta(t)}})
\]
is closed and contractible (in the space of Legendrian submanifolds). The closedness is
obvious because $L_{f(c)}$ and $L_c$ have the same boundary. In order to prove
contractibility, observe that for any $w \in B^{2m+2}(\delta)$, one can define a Legendrian
submanifold $\Lambda_w$ in $(H_w \cap S^{2n+1}(1), \alpha_w |H_w \cap S^{2n+1}(1))$ by
\[
\Lambda_w = \Big(\bigsqcup_{1/2 \leq s \leq 2} \Sigma_{w,s\zeta} \Big) \cap S^{2n+1}(1).
\]
By using a relative version of Lemma \ref{th:gray}, one sees that the trivialization $\tau$
can be chosen in such a way that $\tau_w^{-1}(\Lambda_w) \subset \partial M$ is the same for
all $w \in W$. Now our path of Legendrian submanifolds consists of $\Lambda_{\beta(t)}$, so
that its contractibility is obvious for that specific choice of $\tau$. Since any two $\tau$
are homotopic, it follows that contractibility holds for an arbitrary choice. \qed

\begin{lemma} \label{th:refinement}
Assume that the symplectic automorphism $\phi$ in Lemma \ref{th:curves-and-monodromy} has been
chosen $\theta$-exact. Then the isotopy between $L_{f(c)}$ and $\phi(L_c)$ in Lemma
\ref{th:curves-and-monodromy} can be made exact.
\end{lemma}

\proof For $n>1$ there is nothing to prove: first of all, $H^1(M,\partial M;\R)$ is zero, so
that any choice of $\phi$ is $\theta$-exact; and secondly, $H^1(L_c,\partial L_c;\R) = 0$, so
that any Lagrangian isotopy between such submanifolds is exact.

Therefore we now assume that $n = 1$. Since $\phi$ is $\theta$-exact, both $\phi(L_c)$ and
$L_{f(c)}$ are $\theta$-exact Lagrangian submanifolds. With one exception, it is true that
$L_c$ represents a nontrivial class in $H_1(M,\partial M)$, so that the restriction map
$H^1(M,\partial M;\R) \longrightarrow H^1(L_c,\partial L_c;\R)$ is onto; then one can apply
Lemma \ref{th:exist-exact}\ref{item:deform-to-exact-isotopy} to prove the desired result.

The one exception is when $(m,n) = (1,1)$ and $c$ is the unique (up to isotopy) curve joining
the two points of $\Delta$. In this particular case the Milnor fibre $(M,\o)$ is isomorphic to
a cylinder $([-s_0;s_0] \times S^1, ds \wedge dt)$, and $L_c$ is a simple closed curve which
divides $M$ into two pieces of equal area; to verify this last fact, it is sufficient to note
that the involution $\iota$ exchanges the two pieces. Clearly $\phi(L_c)$ also divides $M$
into two pieces of equal area; and this is sufficient to prove that $\phi(L_c)$ is isotopic to
$L_c$ by an exact isotopy. On the other hand, $f(c) \isotopic c$ for any $f \in \Diffeo$, so
that $L_c$ is exact isotopic to $L_{f(c)}$ by Lemma
\ref{th:donaldson}\ref{item:donaldson-two}. \qed

Now choose a set of basic curves $b_0,\dots,b_m$ in $(D,\Delta)$, and let $L_0 = L_{b_0},
\dots, L_m = L_{b_m}$ be the corresponding Lagrangian submanifolds of $(M,\o,\alpha)$; from
the discussion above, it follows that these do indeed satisfy \eqref{eq:a-chain}. 
As explained in Section \ref{subsec:disc}, the choice of $b_0,\dots,b_m$
allows one to identify each of the canonically isomorphic groups \eqref{eq:canonical-isos} with
$B_{m+1}$. In particular, the symplectic monodromy can now be thought of as a map $\rho_s :
B_{m+1} \longrightarrow \pi_0(\Symp(M,\partial M,\o))$.

\proof[Proof of Theorem \ref{th:gin-two}] From the various isomorphisms which we have
introduced, it follows that $\phi = \phi_\sigma$ and $f = f_\sigma$ are related as in Lemmas
\ref{th:curves-and-monodromy} and \ref{th:refinement}. Using those two results together with
Lemma \ref{th:floer-cohomology}, one finds that $\dim HF(L_i,\phi(L_j)) = \dim
HF(L_{b_i},L_{f(b_j)}) = 2\,I(b_i,f(b_j))$. \qed

\proof[Proof of Corollary \ref{th:faithful-two}] Let $\sigma \in B_{m+1}$ be an element whose
image under the monodromy map is the trivial class in $\pi_0(\Symp(M,\partial M,\o))$. Then
one can take $\phi_\sigma = \id_M$ in Theorem \ref{th:gin-two}, so that
\[
I(b_i,f_\sigma(b_j)) = \half\dim HF(L_i,\phi_\sigma(L_j)) = \half\dim HF(L_i,L_j) = I(b_i,b_j)
\]
for all $0 \leq i,j \leq m$. The same holds for $f_\sigma^2$, since $\rho_s(\sigma^2)$ is
again trivial. In view of Lemma \ref{th:detect-identity}, it follows that $\sigma$ must be the
trivial element. \qed

%
%
\subsection{Bigraded curves and graded Lagrangian submanifolds\label{subsec:both-gradings}}

We use the notation $(D,\Delta)$, $h$, $L_c$, and $(M,\o,\alpha,\theta)$ as in the previous
section. As explained in Section \ref{subsec:bigraded-curves}, the polynomial $h$ can be used
to define a map $\delta_P$ which classifies the covering $\tP$ of $P$, and one can express the
notion of bigraded curves in terms of this map. On the other hand, $M$ is the intersection of
a complex hypersurface $H_{\hat{w}}$ with the unit ball $B^{2n+2}(1)$. As explained in Section
\ref{subsec:grading}, the polynomial $g_{\hat{w}}$ which defines $H_{\hat{w}}$ determines a
map $\delta_{\LL}: \LL \longrightarrow \C^*/\R^{>0}$; and one can use this map to define a
covering $\tLL$ and the corresponding notion of graded Lagrangian submanifold of $M$.

\begin{lemma} \label{th:gradings-and-bigradings}
Let $c$ be a curve in $(D,\Delta)$ which meets $\Delta$. Any bigrading $\tc$ of $c$ determines
in a canonical way a grading $\tL_c$ of $L_c$. Moreover, if $\tL_c$ is the grading associated
to $\tc$, the grading associated to $\Pchi(r_1,r_2)\tc$ is $\chi_{\tLL}(r_1+nr_2)\tL_c$.
\end{lemma}

\proof Let $\gamma: (0;1) \longrightarrow D \setminus \Delta$ be an embedding which
parametrizes an open subset of $c$, and $x \in M$ a point of $\Sigma_{\hat{w},\gamma(t)}
\times \{\gamma(t)\} \subset L_c$. A straightforward computation, starting from
\eqref{eq:explicit-deltal}, shows that
\begin{equation} \label{eq:computation}
\delta_{\LL}(T_xL_c) = (-h(\gamma(t)))^{n-2}\gamma'(t)^2 \in \C^*/\Rg.
\end{equation}
Now assume that we have a bigrading of $c$, which is a map $\tc = (\tc_1,\tc_2): c \setminus
\Delta \longrightarrow \R^2$ with specific properties. Define
\[
\tL_c: L \setminus \bigcup_{z \in c \cap \Delta} \Sigma_{\hat{w},z} \times \{z\}
\longrightarrow \R
\]
by $\tL_c(x) = \tc_1(\gamma(t)) + n\tc_2(\gamma(t))$ for a point $x$ as before. A comparison
of \eqref{eq:computation} with \eqref{eq:explicit-deltal} and \eqref{eq:explicit-deltap} shows
that $\tL_c$ is a grading of $L$, defined everywhere except at finitely many points. It
remains to prove that $\tL_c$ can be extended continuously over those missing points. If $n
\geq 2$, this is true for general reasons. Namely, whenever has a infinite cyclic covering of
an $n$-dimensional manifold (which in our case is $s_L^*\tLL$) and a continuous section of it
defined outside a finite subset, the values at the missing points can always be filled in
continuously. In the remaining case $n = 1$, any one of the missing points divides $L_c$
locally into two components, and one has to check that the limits of $\tL_c$ from both
directions (which certainly exist) are the same. But this follows from the fact that $\tL_c$
is invariant under the involution $\iota$.

The last sentence of the Lemma is obvious from the construction. \qed

\begin{lemma} \label{th:two-indices}
Let $c_0,c_1$ be two curves in $(D,\Delta)$ which meet $\Delta$ and which intersect
transversally. Let $\tc_0,\tc_1$ be two bigradings, and $(r_1,r_2) = \mubigr(\tc_0,\tc_1;z)$
be the local index at a point $z \in c_0 \cap c_1$. Let $\tL_{c_0}, \tL_{c_1}$ be the
corresponding gradings of $L_{c_0},L_{c_1}$. Then the absolute Maslov index at any point $x
\in \Sigma_{\hat{w},z} \subset L_{c_0} \cap L_{c_1}$ is
\[
\muabs(\tL_{c_0},\tL_{c_1};x) = r_1 + nr_2.
\]
\end{lemma}

\proof In view of the last part of Lemma \ref{th:gradings-and-bigradings} it is sufficient to
prove this for $(r_1,r_2) = (0,0)$. Take an isotopy $(f_t)_{0 \leq t \leq 1}$ in
$\Diff(D,\partial D;\Delta)$ with $f_t(z) = z$ for all $t$, $f_0 = \id$, and $f_1(c_0) = c_1$
in some neighbourhood of $z$. We also assume that $Df_t(T_zc_0)$ rotates clockwise with
positive speed during the isotopy, and that the total angle by which is moves is less than
$\pi$ (see Figure\ \ref{fig:turning}). Let $\tL_{f_t(c_0)}$ be the famiy of graded Lagrangian
submanifolds corresponding to the bigraded curves $\tilde{f}_t(\tc_0)$.
\includefigure{turning}{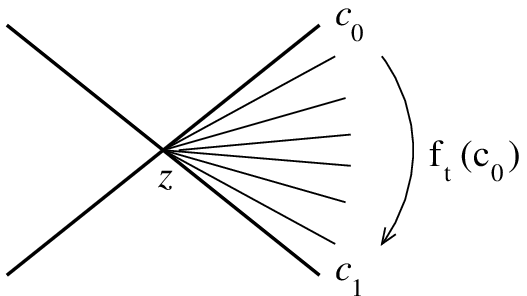}{hb}%

By definition of the local index, $\tilde{f}_1(\tc_0)$ is equal to $\tc_1$ near $z$. Hence
$\tL_{f_t(c_0)}$ agrees with $\tL_{c_1}$ near $\Sigma_{\hat{w},z}$. It follows that for any $x
\in \Sigma_{\hat{w},z}$ one can use the paths $\lambda_0,\lambda_1: [0;1] \longrightarrow
\LL_x$, $\lambda_0(t) = T_xL_{f_t(c_0)}$, $\lambda_1(t) = T_xL_{c_1}$, to compute the absolute
Maslov index. Now, it is not difficult to see that these two paths satisfy the conditions of
Lemma \ref{th:compute-index}, so that the absolute Maslov index is in fact zero. \qed

\begin{lemma} \label{th:graded-floer-cohomology}
Let $(c_0,\tc_0)$ and $(c_1,\tc_1)$ be two bigraded curves in $(D,\Delta)$, such that both
$c_0$ and $c_1$ intersect $\Delta$. Let $\tL_{c_0},\tL_{c_1}$ be the corresponding graded
Lagrangian submanifolds of $(M,\o,\alpha)$. Then the Poincar{\'e} polynomial of the graded
Floer cohomology is obtained by setting $q_1 = q$, $q_2 = q^n$ in the bigraded intersection
number:
\[
\sum_{r \in \Z} q^r \dim HF^r(\tL_{c_0},\tL_{c_1}) = \Ibigr(\tc_0,\tc_1)_{q_1 = q,\, q_2 =
q^n}.
\]
\end{lemma}

\proof This is essentially the same argument as for Lemma \ref{th:floer-cohomology}. Consider
the case when $c_0 \cap c_1 \cap \partial D = \emptyset$ and $c_0 \not\isotopic c_1$. We may
assume that $c_0,c_1$ have minimal intersection. Using Proposition \ref{th:graded-symmetry},
and Lemma \ref{th:two-indices} to compute the relevant absolute Maslov indices, one finds that
the following holds: a point $z \in c_0 \cap c_1$ with local index $(\mu_1(z),\mu_2(z)) \in
\Z^2$ contributes
\[
\begin{cases}
 q^{\mu_1(z) + n\mu_2(z)}(1 + q^{n-1}) & \text{if } z \notin \Delta, \\
 q^{\mu_1(z) + n\mu_2(z)} & \text{if } z \in \Delta \\
\end{cases}
\]
to the Poincar{\'e} polynomial of the Floer cohomology. The factor $(1 + q^{n-1})$ in the
first case comes from the ordinary cohomology of the $(n-1)$-sphere $\Sigma_{\hat{w},z}$. It
is now straightforward to compare this with the definition of the bigraded intersection number
$\Ibigr(\tc_0,\tc_1)$.

As in Lemma \ref{th:floer-cohomology}, the case $c_0 \isotopic c_1$ can be checked by a simple
computation, and that case $c_0 \cap c_1 \cap \partial D \neq \emptyset$ follows from the
previously considered one. \qed

\proof[Proof of Corollary \ref{th:comparison}] One can forget about graded Lagrangian
submanifolds and consider Lemma \ref{th:graded-floer-cohomology} as a statement about the
relative grading on $HF^*(L_{c_0},L_{c_1})$. Using this, and arguing as in the proof of
Theorem \ref{th:gin-two}, one sees that (for some choice of relative grading on the left hand
side)
\[
\sum_{r \in \Z} q^r \dim HF^r(L_i,\phi(L_j)) = \Ibigr(b_i,f(b_j))_{q_1 = q,\, q_2 = q^n}.
\]
The proof is completed by comparing this with Proposition \ref{th:dimequals}. \qed

\providecommand{\bysame}{\leavevmode\hbox to3em{\hrulefill}\thinspace}

\end{document}